\documentclass[a4paper,11pt,reqno]{amsart}
\usepackage[a4paper,hscale=0.7,vscale=0.75,centering]{geometry}

\usepackage{fancyhdr}
\usepackage{amsthm}
\usepackage{hyperref}
\usepackage{amssymb}
\usepackage{latexsym}
\usepackage{amsmath}
\usepackage{amsfonts}
\usepackage{mathabx}
\usepackage[dvips]{graphicx}
\usepackage[utf8]{inputenc}

\usepackage[LGR,T1]{fontenc}%
\usepackage[francais,english]{babel}
\usepackage{url}
\usepackage{enumerate}
\usepackage{hyperref}
\usepackage[usenames,dvipsnames]{color}

\def \RR{\mathbb R}
\def \urpi{\underline{\varpi}}
\def \EE{\mathbb E}
\def \PP{\mathbb P}

\def \FF{\mathbb F}
\providecommand{\U}[1]{\protect\rule{.1in}{.1in}}

\newtheorem{theorem}{Theorem}

\newtheorem{corollary}[theorem]{Corollary}

\newtheorem{example}[theorem]{Example}

\newtheorem{lemma}[theorem]{Lemma}

\newtheorem{proposition}[theorem]{Proposition}
\newtheorem{remark}[theorem]{Remark}

\renewenvironment{proof}[1][Proof]{\noindent\textbf{#1.}}{}

\newcommand{\cD}{\mathcal{D}}
\newcommand{\cE}{\mathcal{E}}
\newcommand{\cF}{\mathcal{F}}

\newcommand{\cP}{\mathcal{P}}

\newcommand{\cR}{\mathcal{R}}
\newcommand{\cS}{\mathcal{S}}

\newcommand{\cU}{\mathcal{U}}

\newcommand{\cW}{\mathcal{W}}

\newcommand{\R}{\mathbb{R}}

\def \proof{{\noindent \bf Proof. }}
\def \eproof{\hbox{ }\hfill$\Box$}

\newcommand{\ud}{\mathrm{d}}

\newcommand{\HYP}[1]
    {\ensuremath{({H#1} ) }}

\newcommand{\set}[1]
    {\ensuremath{\{ #1 \}}}
    
\newcommand{\HP}[1] %L2DPDT sur 0,T
    {\ensuremath{\mathscr{H}^{#1}}}

\newcommand{\esp}[1]{\ensuremath{\mathbb{E} \, \!\! \left[#1\right] }}

\newcommand{\EFp}[2]
    {\ensuremath{% raccourci esperance conditionnelle discretisation
     \mathbb{E}_{#1}\, \!\!\left[#2\right] }}
     
\newcommand{\ti}[1]{t_{i #1}}

\renewcommand{\Xi}[1]{X_{i #1}}

\newcommand{\whesp}[1]
    {\ensuremath{%
     \widehat{\mathbb{E}}\!\left[#1\right] }}

\newcommand{\pic}[2]{\ensuremath{ 
\text{\texttt{picard[}}#1\text{\texttt{](}}#2\text{\texttt{)}}
}}

\newcommand{\sol}[2]{\ensuremath{ 
\text{\texttt{solver[}}#1\text{\texttt{](}}#2\text{\texttt{)}}
}}

\newcommand{\tsol}[4]{\ensuremath{ 
\overline{\text{\texttt{solver}}}\text{\texttt{[}}#1\text{\texttt{](}}#2\text{\texttt{,}}#3\text{\texttt{,}}#4\text{\texttt{)}}
}}

\newcommand{\tX}[1]{\tilde{X}^{#1}}
\newcommand{\tY}[1]{\tilde{Y}^{#1}}
\newcommand{\tZ}[1]{\tilde{Z}^{#1}}

\newcommand{\Xd}{\bar{X}}
\newcommand{\Ud}{\bar{U}}

\newcommand{\Vd}{\bar{V}}

\newcommand{\Yd}{\bar{Y}}

\newcommand{\Zd}{\bar{Z}}

\newcommand{\ccesp}[1]{\ensuremath{\hat{\mathbb{E}} \; \!\! \left[#1\right] }}

\newcommand{\cc}[1]{\ensuremath{ \langle #1 \rangle} } %copie d'une v.a.
\newcommand{\law}[1]{\ensuremath{ [#1]} } %loi d'une v.a. 

\newcommand{\Tr}[1]{\ensuremath{{\rm Tr}[#1] } }

\newcommand{\norm}[1]{\left \|{#1} \right \|_2} 
\newcommand{\normalpha}[1]{\bigl \|{#1} \bigr \|_{2\alpha}} %\vvvert
\newcommand{\normalphas}[1]{\|{#1}  \|_{2\alpha}}
\newcommand{\bd}{\bar{\delta}}

\newcommand{\normalphadeuxB}[1]{\Bigl \|{#1} \Bigr \|_{\alpha}}

\begin{document}

\title[Numerical Method for FBSDEs of McKean-Vlasov Type]{Numerical Method for FBSDEs of McKean-Vlasov Type}

\author[Jean-Fran\c{c}ois Chassagneux, Dan Crisan and Fran\c{c}ois Delarue]{Jean-Fran\c{c}ois Chassagneux$^*$, Dan Crisan$^{\dagger}$ and Fran\c{c}ois Delarue$^{\ddagger}$ }

\maketitle

\renewcommand{\thefootnote}{\fnsymbol{footnote}}
 \footnotetext[1]{Laboratoire de Probabilit\'es et Mod\`eles al\'eatoires, Universit\'e Paris Diderot. {\sf chassagneux@math.univ-paris-diderot.fr}}
 \footnotetext[2]{Department of Mathematics, Imperial College London. {\sf d.crisan@imperial.ac.uk}}
 \footnotetext[3]{Laboratoire Jean-Alexandre Dieudonn\'e, 
 Universit\'e de Nice Sophia-Antipolis. {\sf delarue@unice.fr}}

 \renewcommand{\thefootnote}{\arabic{footnote}}
 
\begin{abstract}
This paper is dedicated to the presentation and the analysis of a numerical scheme for forward-backward SDEs of the McKean-Vlasov type, or equivalently for solutions to PDEs on the Wasserstein space. Because of the mean field structure of the equation, earlier methods for classical forward-backward systems fail. 
The scheme is based on a variation of the method of continuation. The principle is to implement recursively local Picard iterations on small time intervals. 

We establish a bound for the rate of convergence under the assumption that the decoupling field of the forward-bakward SDE (or equivalently the solution of the PDE) satisfies mild regularity conditions. We also provide numerical illustrations.
\end{abstract}

%\tableofcontents

\section{Introduction}

In this paper, we investigate a probabilistic numerical method to approximate
the solution of the following non-local PDE 
\begin{equation}
  \label{eq pde decoupling field}
\begin{split}
 &\partial_t \cU(t,x,\mu) 
+ b(x,\cU(t,x,\mu),\nu) \cdot \partial_x \cU(t,x,\mu) 
\\
&\hspace{120pt}+ \frac12 \Tr{ \partial_{xx}^2 \cU(t,x,\mu)a(x,\mu)}  
+ f\bigl(x,\cU(t,x,\mu),\partial_x \cU(t,x,\mu) \sigma(x,\mu),\nu\bigr) 
\\
&\hspace{5pt}+ \int_{\R^d}\partial_\mu \cU(t,x,\mu)(\upsilon) 
\cdot b(\upsilon,\cU(t,\upsilon,\nu),\nu)\ud \mu(\upsilon)  
 + \int_{\R^d} \frac12 \Tr{\partial_{x}\partial_\mu \cU(t,x,\mu)(\upsilon) a(\upsilon,\mu)} \ud \mu(\upsilon)
 = 0 \; ,
\end{split}
\end{equation}
for $(t,x,\mu) \in [0,T)\times\R^d\times \cP_2(\R^d)$ with the terminal condition $\cU(T,\cdot)=g(\cdot)$, where $\nu$ is a notation for the image of the probability measure $\mu$ by the mapping $\RR^d \ni x \mapsto 
(x,\cU(t,x,\mu)) \in \RR^{2d}$. Above, $a(x,\mu)=[\sigma \sigma^\dagger](x,\mu)$. The set $\cP_2(\R^d)$ is the set of probability measures with a finite second-order moment, endowed with the Wasserstein distance i.e.
\begin{align*}
\cW_2(\mu,\mu') := \inf_{\pi} 
\left(
\int_{\RR^d \times \RR^d} \vert x-x' \vert^2 \ud \pi(x,x')
\right)^{\frac12}
\;,
\end{align*}
for  $(\mu,\mu') \in \cP_2(\R^d) \times \cP_2(\R^d)$, the infimum being taken over the probability distributions $\pi$ on $\RR^d \times \RR^d$ whose marginals on $\RR^d$ are respectively $\mu$ and $\mu'$. 
\vspace{5pt}

Whilst the first two lines in 
\eqref{eq pde decoupling field}
form a classical non-linear parabolic equations, the 
last two terms are non-standard. Not only are they non-local, in the sense that the solution or its derivatives are computed at points $\upsilon$ different from $x$, but also they involve derivatives in the argument $\mu$, which lives in a space of probability measures. In this regard, 
the notation $\partial_{\mu} \cU(t,x,\mu)(\upsilon)$ denotes the so-called \textit{Wasserstein derivative} of the function $\cU$ in the direction of the measure, computed at point
$(t,x,\mu)$ and taken at the \textit{continuous coordinate} $\upsilon$. 
We provide below a short reminder of the construction of this derivative, as introduced by Lions, see \cite{Cardaliaguet} or \cite[Chap. 5]{CarmonaDelarue_book_I}.  

These PDEs arise in the study of large population stochastic control
problems, either of mean field game type, see for instance 
\cite{Cardaliaguet,CardaliaguetDelarueLasryLions,chacri15,Lionsvideo}
or 
\cite[Chap. 12]{CarmonaDelarue_book_II}
and the references therein,
or of mean field control type, see for instance
 \cite{BensoussanFrehseYam2,BensoussanFrehseYam3,chacri15,PhamWei1}. In both cases, $\cU$ plays 
 the role of a value function or, when the above equation is replaced by a system of equations of the same form, the gradient of the value function. 
 Generally speaking, these types of equations are known as ``master equations''. We refer to the aforementioned papers and monographes for a complete overview of the subject, in which existence and uniqueness of classical or viscosity solutions have been studied. In particular, in our previous paper \cite{chacri15}, we tackled classical solutions by connecting $\cU$ with a system of fully coupled Forward-Backward Stochastic Differential Equations of the McKean-Vlasov type (MKV FBSDE), for which $\cU$ plays the role of a decoupling field. We also refer to 
 \cite[Chap. 12]{CarmonaDelarue_book_II} for a similar approach. 
 
% More precisely, it is introduced there as the decoupling field of this equation (see equation \cite{eq decoupling}
%below for a clear explanation of the name).
In the current paper, we build on this link to design our numerical method.
\vskip 5pt

The connection between $\cU$ and FBSDEs may be stated as follows. Basically, $\cU$ may be written as $\cU(t,x,\mu) = Y^{t,x,\mu}_t$ for all $(t,x,\mu) \in [0,T]\times \R^d \times \cP_2(\R^d)$,
where $Y^{t,x,\mu}$ together with $(X^{t,x,\mu},Z^{t,x,\mu})$ solves the following standard FBSDE:
\begin{align}
 X^{t,x,\mu}_s &= x + \int_t^s b(X^{t,x,\mu}_{r},Y^{t,x,\mu}_{r},\law{X^{t,\xi}_r,Y^{t,\xi}_r}) \ud r + \int_t^s \sigma(X^{t,x,\mu}_r,\law{X^{t,x,\mu}_r}) \ud W_r
\label{eq Xtxmu}
 \\
 Y^{t,x,\mu}_s &= g(X^{t,x,\mu}_{T},\law{X_T^{t,\xi}}) + \int_s^T f(X^{t,x,\mu}_{r},Y^{t,x,\mu}_{r},Z^{t,x,\mu}_r,\law{X^{t,\xi}_r,Y^{t,\xi}_r})  \ud r - \int_s^T Z^{t,x,\mu}_r \cdot \ud W_r \, ,
\label{eq Ytxmu}
 \end{align}
 which is parametrized by the law of the following MKV FBSDE: 
\begin{align}
 X^{t,\xi}_s &= \xi + \int_t^s b(X^{t,\xi}_{r},Y^{t,\xi}_{r},\law{X^{t,\xi}_r,Y^{t,\xi}_r}) \ud r + \int_t^s \sigma(X^{t,\xi}_r,\law{X^{t,\xi}_{r}}) \ud W_r
 \label{eq Xtxi}
 \\
 Y^{t,\xi}_s &= g(X^{t,\xi}_{T},\law{X_T^{t,\xi}}) + \int_s^T f(X^{t,\xi}_{r},Y^{t,\xi}_{r},Z^{t,\xi}_r,\law{X^{t,\xi}_r,Y^{t,\xi}_r})  \ud r - \int_s^T Z^{t,\xi}_r \cdot \ud W_r \, ,
 \label{eq Ytxi}
 \end{align}
 where $(W_{t})_{0 \le t \le T}$ is a Brownian motion and $\xi$ has $\mu$ as distribution. 
 In the previous equations and in the sequel, we use the notation $\law{\theta}$ for the law of a random variable $\theta$. In particular, in the above, we have that $\law{\xi} = \mu$.
% \textcolor{magenta}{
 %comment: The form of the coefficients is not completely general, but reflects the main difficulties one can encounter.
 %}
 So, 
to obtain an approximation of $\cU(t,x,\mu)$ given by the initial value of  \eqref{eq Ytxmu}, our strategy is to
approximate the system \eqref{eq Xtxi}-\eqref{eq Ytxi} as its solution appears in the coefficients of
\eqref{eq Xtxmu}-\eqref{eq Ytxmu}.
In this regard, our approach is probabilistic. 

Actually, our paper is not the first one to address the numerical approximation of equations of the type 
\eqref{eq pde decoupling field} by means of a probabilistic approach.
In its PhD dissertation, Alanko \cite{Alanko} develops a numerical method for mean field games based upon a Picard iteration: Given the proxy 
for the equilibrium distribution of the population (which is represented by the mean field component in the above FBSDE), one solves for the value function by approximating the solution of the (standard) BSDE associated with the control problem; given the solution of the BSDE, we then get a new proxy for the equilibrium distribution and so on... Up to a Girsanov transformation, the BSDE associated with the control problem coincides with the backward equation in the above FBSDEs.
In \cite{Alanko}, the Girsanov transformation is indeed used to decouple the 
forward and backward equations and it is the keystone of the paper to 
address the numerical impact of the change of measure onto the mean field component. 
 Under our setting, this method would more or less consist in solving for the backward equation given a proxy for the forward equation and then in iterating, which is what we call the \textit{Picard method} for the FBSDE system. Unfortunately, 
convergence of the Picard iterations is a difficult issue, as the convergence is known in small time only, see the numerical examples in Section \ref{se numerics} below. It is indeed well-known that Picard theorem only applies in small time for fully coupled problems. In this regard, it must be stressed that our system  \eqref{eq Xtxi}-\eqref{eq Ytxi}
is somehow doubly coupled, once in the variable $x$ and once in the variable $\mu$, which explains why a change measure does not permit to decouple it entirely.   

The goal of our paper is precisely to go further and to propose an algorithm whose convergence is known on any interval of a given length (observe that the convergence is not studied in \cite{Alanko}). 
In the classical case, this question has been addressed by several authors, among which \cite{delmen06,delmen08} and \cite{BenderZhang}, but all these methods rely on the Markov structure of the problem. Here, the Markov property is true but at the price of regarding the entire $\RR^d \times \cP_{2}(\RR^d)$ as state space: The fact that the second component is infinite dimensional makes intractable the complexity of these approaches.  
To avoid any similar problem, we use a pathwise approach for the forward component; it consists in iterating successively the \textit{Picard method} on small intervals, all the Picard iterations being implemented with 
a tree
approximation of the Brownian motion. This strategy is inspired from the method of continuation, the parameter in the continuation argument being the time length $T$ itself. The advantage for working on a tree is twofold: 
as we said, we completely bypass any Markov argument; also, we get, not only, an approximation of the system \eqref{eq Xtxi}-\eqref{eq Ytxi}
but also, for free, an approximation of the system  
\eqref{eq Xtxmu}-\eqref{eq Ytxmu}, which ``lives'' on a subtree obtained by conditioning on the initial root. 
We prove that the method is convergent and provide a rate of convergence for it. Numerical examples are given in Section \ref{se numerics}. 
Of course, the complexity remains pretty high in comparison with the methods developed in the classical non McKean-Vlasov case. This should not come as a surprise since, as we already emphasized, the problem is somehow infinite dimensional. 

We refer the interested reader to 
the following papers for various numerical methods, based upon finite differences or variational approaches, for mean field games:
\cite{AchdouCapuzzo-Dolcetta,AchdouCamilliCapuzzo-Dolcetta2,AchdouPorretta}
and 
\cite{BenamouCarlier,LachapelleSalomonTurinici,Gueant_numerical}. {Recently, a Markov chain approximation method was {also} suggested in \cite{baybud16}.}
\vskip 5pt

%\textcolor{magenta}
{
The paper is organized as follows. 
The method for the system 
\eqref{eq Xtxi}-\eqref{eq Ytxi} is exposed in Section \ref{se algo MKVFBSDE 1}. 
The convergence is addressed in Section \ref{se algo MKVFBSDE analysis}.  In Section \ref{se numerics}, we explain how to compute in practice $\cU(t,x,\mu)$ (and thus approximate 
\eqref{eq Xtxmu}-\eqref{eq Ytxmu})
from the approximation of the sole \eqref{eq Xtxi}-\eqref{eq Ytxi} and we present some numerical results validating empirically the convergence results obtained in Section \ref{se algo MKVFBSDE analysis}. We collect in the appendix some key results for the convergence analysis.
}

\vspace{5pt}
%\textcolor{red}{say a word on system HJB-Fokker Planck, its approx and the link with what we have
%$u(t,x)=U(t,x,\mu_t)$ and where $\mu_t$ is the flow of measure for \eqref{eq Xtxi}}

%%%%%%%

\section{A new algorithm for coupled forward backward systems}
\label{se algo MKVFBSDE 1}
As announced right above, we will focus on the approximation of 
the following type of McKean-Vlasov forward-backward stochastic differential equation:
\begin{equation}
\label{eq fbsde}
\begin{split}
&dX_{t} = b\bigl(X_{t},Y_{t},\law{X_t,Y_t}\bigr) dt + \sigma\bigl(X_t,\law{X_{t}}\bigr)dW_{t},
\\
&dY_{t} = -f\bigl(X_{t},Y_{t},Z_t,\law{X_t,Y_t}\bigr)\ud t + Z_{t} \cdot dW_{t}, \quad t \in [0,T] \ , 
\\
&Y_{T} = g\bigl(X_{T},\law{X_T}\bigr)\,\;\text{ and } \; X_0 = \xi \, ,
\end{split}
\end{equation}
for some time horizon $T>0$. 
Throughout the analysis, the equation is regarded on a complete filtered probability space $(\Omega,\cF,\FF,\PP)$, equipped with a $d$-dimensional $\FF$-Brownian motion $(W_{t})_{0 \leq t \leq T}$. To simplify, we assume that the state process $(X_{t})_{0 \leq t \leq T}$ is of the same dimension. 
The process $(Y_{t})_{0 \le t \le T}$ is $1$-dimensional. As a result, 
$(Z_{t})_{0 \le t \le T}$ is $d$-dimensional. 

In \eqref{eq fbsde}, the three processes $(X_{t})_{0 \le t \le T}$, $(Y_{t})_{0 \le t \le T}$ and $(Z_{t})_{0 \le t \le T}$ are required to be 
$\FF$-progressively measurable. Both $(X_{t})_{0 \le t \le T}$
and $(Y_{t})_{0 \le t \le T}$ have continuous trajectories. Generally speaking, the initial condition $X_{0}$ is assumed to be square-integrable, but 
at some point, we will assume that $X_{0}$ belongs to $L^p(\Omega,\cF,\PP;\RR^d)$, for some $p >2$. Accordingly, $(X_{t})_{0 \le t \le T}$, $(Y_{t})_{0 \le t \le T}$ and $(Z_{t})_{0 \le t \leq T}$ must satisfy:
\begin{equation*}
\vvvert (X,Y,Z) \vvvert_{[0,T]}
:= 
\EE \biggl[ \sup_{0 \le t \le T} \bigl( \vert X_{t} \vert^2 + \vert Y_{t} \vert^2 \bigr) + \int_{0}^T \vert Z_{t} \vert^2 dt \biggr]^{1/2} < \infty. 
\end{equation*}

The domains and codomains of the coefficients are defined accordingly. 
The assumption that $\sigma$ is assumed to be independent of the variable $y$ is consistent with the global solvability results that exist in the literature for equations like \eqref{eq fbsde}. For instance, it covers cases coming from optimization theory for large mean field interacting particle systems. We refer to our previous paper \cite{chacri15} for a complete overview on the subject, together with the references
\cite{BensoussanFrehseYam,Cardaliaguet,CarmonaDelarue_book_I,CarmonaDelarue_book_II,CarmonaDelarueLachapelle}. In light of the examples tackled in 
\cite{chacri15}, the fact that $b$ is independent of $z$ may actually seem more restrictive, as it excludes cases when the forward-backward system of the McKean-Vlasov type is used to represent the value function of the underlying optimization problem. It is indeed a well-known fact that, with or without McKean-Vlasov interaction, the value function of a standard optimization problem may be represented as the backward component of a standard FBSDE with a drift term depending upon the $z$ variable. This says that, in order to tackle the aforementioned optimization problems of the mean field type by means of the numerical method investigated in this paper, one must apply the algorithm exposed below to the Pontryagin system. The latter one is indeed of the form \eqref{eq fbsde}, provided that 
$Y$ is allowed to be multi-dimensional. (Below, we just focus on the one-dimensional case, but the adaptation is straightforward.)

In fact, our choice for assuming $b$ to be independent of $z$ should not come as a surprise. The same assumption appears in the papers \cite{delmen06,delmen08} dedicated to the numerical analysis of standard FBSDEs, which will serve us as a benchmark throughout the text. See however Remark \ref{re z meme pas peur}.

Finally, the fact that the coefficients are time-homogeneous is for convenience only. 
\vskip 5pt

As a key ingredient in our analysis, we use the following representation result given in e.g. Proposition 2.2 in \cite{chacri15},
\begin{align} \label{eq link U-Y}
Y^\xi_t := \cU(t,X^\xi_t,\law{X_t^\xi})\;,
\end{align}
where $\cU : [0,T] \times \RR^d \times \cP_{2}(\RR^d) \rightarrow \RR$ is assumed to be the classical solution, in the sense of \cite[Definition 2.6]{chacri15}, to \eqref{eq pde decoupling field}. In this regard, 
the derivative with respect to the measure argument is defined according to Lions' approach to the Wasserstein derivative. In short, the \textit{lifting} $\hat{\cU}$ of $\cU$ to $L^2(\Omega,\cF_{0},\PP;\RR^d)$, which we define by
\begin{equation*}
\hat{\cU}(t,x,\xi) = \cU(t,x,[\xi]), \qquad t \in [0,T], \ x \in \RR^d, \ \xi \in L^2(\Omega,\cF_{0},\PP;\RR^d),
\end{equation*}
is assumed to be Fr\'echet differentiable. Of course, this makes sense as long as
the space $(\Omega,\cF_{0},\PP)$ is rich enough so that, for any 
$\mu \in \cP_{2}(\RR^d)$, there exists a random variable 
$\xi \in L^2(\Omega,\cF_{0},\PP;\RR^d)$ such that $\xi \sim \mu$. So, in the sequel, $(\Omega,\cF_{0},\PP)$ is assumed to be atomless, which makes it rich enough. 
A crucial point with Lions' approach to Wasserstein differential calculus is that the Fr\'echet derivative of $\hat{\cU}$ in the third variable, which can be identified with a square-integrable random variable, may be represented at point $(t,x,\xi)$ 
as $\partial_{\mu} \cU(t,x,[\xi])(\xi)$ for a mapping $\partial_{\mu} \cU(t,x,\mu)(\cdot) : \RR^d \ni v \mapsto \partial_{\mu} \cU(t,x,\mu)(v) \in \RR^d$. 
This latter function plays the role of Wasserstein derivative of $\cU$ in the measure argument. To define a classical solution, it is then required that 
$\RR^d \ni v \mapsto \partial_{\mu} \cU(t,x,\mu)(v)$ is differentiable, both 
 $\partial_{\mu} \cU$ and 
 $\partial_{v} \partial_{\mu} \cU$ being required to be continuous at any point $(t,x,\mu,v)$ such that $v$ is in the support of $\mu$.

%\vspace{5pt}
%We assume that \\
%- \HYP{1}  there exists a decoupling field $U$ for the above equation satisfying $U(t,X_t,\law{X_t})=Y_t$ 
%and ``classical solution'' to
%\begin{align}
% \partial_t U(t,x,\mu) &+ b(U(t,x,\mu),\nu) \partial_x U(t,x,\mu) + \frac12 \Tr{ \partial_{xx}^2 U(t,x,\mu)\sigma\sigma(x)^{\dagger}}  
% \nonumber
% \\
% &+ f(\partial_x U(t,x,\mu) \sigma(x)) + \int_{\R^d}\partial_\mu U(t,x,\mu)(\upsilon) b(U(t,\upsilon,\nu),\nu)\ud \mu(\upsilon)  \nonumber
% \\
% & \quad\quad\quad+ \int_{\R^d} \frac12 \Tr{\partial_{x}\partial_\mu U(t,x,\mu)(\upsilon) \sigma\sigma^{\dagger}(\upsilon)} \ud \mu(\upsilon)
% = 0 \;   \label{eq pde decoupling field}
%\end{align}
%with the terminal condition $u(T,\cdot)=g(\cdot)$ and where $\nu = \law{(\xi,U(t,\xi,\mu))}$ for $\law{\xi} = \mu$.

\vspace{5pt}
\noindent \textbf{Assumptions.} Our analysis requires some minimal regularity assumptions on the coefficients $b$, $\sigma$, $f$ and the function $\cU$. As for the coefficients functions, we assume that there exists a constant $\Lambda \geq 0$ such that:
\vspace{5pt}

\noindent - \HYP{0}: The functions $b$, $\sigma$, $f$ and $g$ are $\Lambda$-Lipschitz continuous in all the variables, the space $\cP_{2}(\RR^d)$ being equipped with the Wasserstein distance $\cW_{2}$. Moreover, the function $\sigma$ is bounded by $\Lambda$.

\vspace{5pt}
We now state the main assumptions on $\cU$, see Remark \ref{re main ass} for comments.
\vspace{5pt}

\noindent - \HYP{1}: 
for any $t \in [0,T]$ and $\xi \in L^2(\Omega,\cF_{t},\PP;\RR^d)$, 
the McKean-Vlasov forward-backward system \eqref{eq fbsde}
set on $[t,T]$ instead of $[0,T]$ with $X_{t}=\xi$ as initial condition at time $t$ has a unique solution $(X_{s}^{t,\xi},Y_{s}^{t,\xi},Z_{s}^{t,\xi})_{t \leq s \leq T}$; in parallel, 
$\cU$ is the classical solution, in the sense of \cite[Definition 2.6]{chacri15}, to \eqref{eq pde decoupling field}; and $\cU$
and its  derivatives satisfy
\begin{align}
&|\cU(t,x,\mu) - \cU(t,x,\mu')| + |\partial_x \cU(t,x,\mu) - \partial_x \cU(t,x,\mu')| \le \Lambda\cW_2(\mu,\mu') \, ,
\label{eq cont grad mu}
\\
& |\partial_x \cU(t,x,{\mu}) |+\norm{\partial_\mu \cU(t,x,\law{\xi})(\xi)} \le \Lambda \; , \label{eq bound first deriv}
\\
&  |\partial^2_{xx} \cU(t,x,{\mu}) |+\norm{\partial_\upsilon \partial_\mu \cU(t,x,\law{\xi})(\xi)} \le \Lambda \; , \label{eq bound second deriv}
\\ 
&{
\text{and } |\partial^2_{xx} \cU(t,x,{\mu}) - \partial^2_{xx} \cU(t,x',{\mu})| \le \Lambda |x-x'|\;, \label{eq lip second deriv}
}
\end{align}
for $(t,x,x',\xi) \in [0,T] \times \RR^{d} \times \RR^{d} \times L^2(\Omega,\cF_{0},\PP;\RR^d)$ and $\mu,\mu' \in \cP_{2}(\RR^d)$. 
Also, we require that
\begin{equation}
\begin{split}
&|\cU(t+h,x,\law{\xi}) - \cU(t,x,\law{\xi})|
+
|\partial_x \cU(t+h,x,\law{\xi}) - \partial_{x}\cU(t,x,\law{\xi})| \le \Lambda h^\frac12 \bigl( 1+ \vert x \vert + \| \xi \|_{2} \bigr) \, ,
\end{split}
\label{eq cont grad t}
\end{equation}
and for all $h\in[0,T)$, $(t,x) \in [0,T-h]\times\R^d$,
$\xi \in L^2(\Omega,\cF_{0},\PP;\RR^d)$ and $v,v' \in \RR^d$,
\begin{align}\label{eq control dvdmu}
 |\partial_\upsilon \partial_\mu\cU(t,x,\law{\xi})(\upsilon)
 &-\partial_\upsilon \partial_\mu \cU(t,x,\law{\xi})(\upsilon')|
 \le
\Lambda \set{1 + |\upsilon|^{2\alpha} + |\upsilon'|^{2\alpha} + \norm{\xi}^{2\alpha}}^\frac12|\upsilon-\upsilon'|\;,
\end{align}
for some $\alpha >0$. 

\begin{remark}\label{re main ass}
In \cite{chacri15}, it is shown that, under some conditions on the coefficients $b$, $f$ and $\sigma$, the PDE \eqref{eq pde decoupling field} has indeed a unique classical solution which satisfies the assumption \HYP{1}.
\begin{enumerate}
\item Estimate \eqref{eq control dvdmu} is obtained by combining Definition 2.6 and Proposition 3.9 in \cite{chacri15}.
A major difficulty in the analysis provided below is the fact that $\alpha$ may be larger than 1, in which case the Lipschitz bound for the second order derivative is super-linear. This problem is proper to the McKean-Vlasov structure of the equation and does not manifest in the classical setting, 
compare for instance with  \cite{delmen06,delmen08}. Below, we tackle two cases: the case when $\alpha \leq 1$, which has been investigated in 
\cite{CardaliaguetDelarueLasryLions}
and \cite[Chap. 12]{CarmonaDelarue_book_II} under stronger conditions on the coefficients, and the case when 
$\alpha >1$ but 
$\cU$ is bounded. 
\item Estimates \eqref{eq cont grad mu}-\eqref{eq cont grad t} are required to control the convergence error when the coefficients ($b$ or $f$)
depend on $Z$. 
\begin{enumerate}
 
\item The estimate \eqref{eq cont grad mu} can be retrieved from the computations made in \cite{chacri15}. 
See the comments at the bottom of page 60, near equation $(4.58)$.
\item  The estimate \eqref{eq cont grad t} comes from the theory of FBSDEs (without McKean-Vlasov interaction).
Indeed, using the Lipschitz property of $\cU$ and $\partial_{x} \cU$ in the variable $\mu$, it suffices to prove 
\begin{equation*}
\begin{split}
&|\cU(t+h,x,\law{X^{t,\xi}_{t+h}}) - \cU(t,x,\law{X^{t,\xi}_{t}})|
+
|\partial_x \cU(t+h,x,\law{X^{t,\xi}_{t+h}}) - \partial_{x}\cU(t,x,\law{X^{t,\xi}_{t}})| 
\\
&\le \Lambda h^\frac12\bigl( 1+ \vert x \vert + \| \xi \|_{2} \bigr) \, .
\end{split}
\end{equation*}
As stated in Proposition 2.2 in \cite{chacri15}, for $\xi \sim \mu$, $\cU(s,x,\law{X^{t,\xi}_{s}}) = u_{t,\mu}(s,x)$ where $u_{t,\mu}$
is solution to a quasi-linear PDE.  
%(\textcolor{red}{Achtung! unif ellipticity assumed there but might be useless... redo the computation in small time? etc.})
Then the estimate \eqref{eq cont grad t} follows from 
standard results on non-linear PDEs, see
e.g. Theorem 2.1 in \cite{delmen06}.
\end{enumerate}
\end{enumerate}
%
%\item For later use, we stress that the notion of classical solution that we use here \cite[Definition 2.6]{chacri15} imposes
%some key bounds on the derivatives of the function $U$, namely
%
%This also implies that $U$ is Lipschitz continuous in the following sense
%\begin{align}\label{eq U Lipschitz}
%|U(t,x,\mu)-U(t,x',\mu')| \le L(|x-x'| + \cW_2(\mu,\mu')) \;,
%\end{align}
%for all $t \in [0,T]$, $(x,x') \in \R^{2d}$ and $(\mu,\mu') \in \cP_2(\R^d) \times \cP_2(\R^d)$.
In comparison with the assumption used in \cite{delmen06}, 
the condition 
\HYP{1} is more demanding. In 
\cite{delmen06}, there is no need for assuming 
the second-order derivative to be Lipschitz in space. This follows from the fact that, here, we approximate the Brownian increments 
by random variables taking a small number of values, whilst in 
\cite{delmen06}, the Brownian increments are approximated by a quantization grid
with a larger number of points. 
In this regard, our approach is closer to the strategy implemented in 
\cite{delmen08}.
\end{remark}

%\textcolor{blue}{T: $X \rightarrow (t,X)$, $Y \rightarrow (t,X,Y)$ and $Z \rightarrow (t,X,Y,Z)$ still ok. }

\subsection{Description}

The goal of the numerical method exposed in the paper is to approximate 
$\cU$. The starting point is the formula \eqref{eq fbsde} and, quite naturally, the strategy is to approximate the process $(X^\xi,Y^\xi,Z^\xi):=(X^{0,\xi},Y^{0,\xi},Z^{0,\xi})$. 

Generally speaking, this approach raises a major difficulty, as it requires to handle the strongly coupled forward-backward structure of \eqref{eq fbsde}.  Indeed, theoretical solutions to \eqref{eq fbsde} may be constructed by means of basic Picard iterations but in small time only, which comes in contrast with similar results for decoupled forward or backward equations for which Picard iterations converge on any finite time horizon. In the papers \cite{delmen06,delmen08} --which deal with the non McKean-Vlasov case--, this difficulty is bypassed by approximating \textit{the decoupling field} $\cU$ at the nodes of a time-space grid. Obviously, this strategy is hopeless in the McKean-Vlasov setting as the state variable is infinite dimensional; discretizing it on a grid would be of a non-tractable complexity. 
This observation is the main rationale for the approach exposed below. 

Our method is a variation of the so-called \textit{method of continuation}. In full generality, it consists in increasing step by step the coupling parameter between the forward and backward equations. Of course, the intuition is that, for a given time length $T$, 
the Picard scheme should converge for
very small values of the coupling parameter. The goal is then to insert the approximation computed for a small coupling parameter into the scheme used to compute a numerical solution for a higher value of the coupling parameter.  
Below, we adapt this idea, but we directly regard $T$ itself as a coupling parameter. So we increase $T$ step and by step and, on each step, we make use of a Picard iteration based on the approximations obtained at the previous steps. 

This naturally motivates the introduction of an equidistant grid $\Re = \set{r_0=0, \dots, r_N = T}$ of the time interval $[0,T]$, with
$r_k = k\delta$ and $\delta = \frac{T}{N}$ for $N \ge 2$. In the following we shall consider that $\delta$ is ``small enough''
and state more precisely what it means in the main results, see Theorem  \ref{th error picard} and Theorem \ref{th propagation error II}.

For $0 \le k \le N-1$, we consider  intervals $I_k = [r_{k},T]$ and on each interval, the following
FBSDE, for $\xi \in L^2(\cF_{r_{k}})$ (which is a shorter notation for 
$L^2(\Omega,\cF_{r_{k}},\PP;\RR^d)$):
\begin{align} \label{eq main step k}
 X_t &=  \xi + \int_{r_k}^t b\bigl(X_{s},Y_s,\law{X_s,Y_s}\bigr)\ud s + \int_{r_k}^t \sigma(X_s,\law{X_{s}})\ud W_s\;,\;\\ %X_{(N-k)\delta} = \xi \in L^2(\cF_{(N-k)\delta})\\
 Y_t &= g\bigl(X_T,\law{X_T}\bigr) + \int_{t}^T f\bigl(X_{s},Y_{s},Z_s,\law{X_s,Y_s}\bigr)\ud s -\int_t^TZ_s \cdot \ud W_s. \label{eq main step k y}
\end{align}
%For $0 \le k \le N-1$, we also consider the interval $[r_{k},r_{k+1}]$
% and the following FBSDE, for $\xi \in L^2(\cF_{r_k})$ and $\chi \in L^2(\cF_{r_{k+1}})$,
% \begin{align} \label{eq de fbsde tsol}
% \left \{
% \begin{array}{rcl}
%  \tilde{X}_t &= &\xi + \int_{r_k}^t b(\tilde{Y}_s,\law{\tilde{X}_s,\tilde{Y}_s}) \ud s + \int_{r_k}^t \sigma(\tilde{X}_s)\ud W_s %\quad \text{ and } \quad \tilde{X}_{(N-k-1)\delta}=  
%  \\
%  \tilde{Y}_t &= & \chi + \int_{t}^{r_{k+1}} f(\tilde{Z}_s,\law{\tX{}_s,\tY{}_s})\ud s -\int_t^T\tilde{Z}_s \ud W_s
%  \end{array}
%  \right.\,.
% \end{align}

%\vspace{5pt} 
%Our goal is to get an approximation of $U(0,x,\mu)$ for some $x \in \R^d$
%\textcolor{red}{BELONGING TO THE SUPPORT OF $\mu$}. To obtain it, we will
% compute an approximation of $U(0,\xi,\mu)$ where $\xi$ is a random variable with distribution $\mu$.
% We will see in Section \ref{subse approx U(t,x,mu)}, how it  allows us to obtain approximation of $U(0,x,\mu)$ in practice. 
 
 \paragraph{Picard iterations.}  We need to compute backwards the value of $\cU(r_{k},\xi,\law{\xi})$ for some $\xi \in L^2(\cF_{r_{k}})$, $0 \le k \le N-2$. 
 We are then going to solve the FBSDE \eqref{eq main step k}-\eqref{eq main step k y} on the interval $I_{k}$. As explained above, the difficulty
 is the arbitrariness of $T$: When $k$ is large, $I_{k}$ is of a small length, but this becomes false as $k$ decreases. Fortunately, we 
 can rewrite the forward-backward system on a smaller interval at the price of changing the terminal boundary condition. Indeed, from $\HYP{1}$, we know that
 $(X_{s}^{r_{k},\xi},Y_{s}^{r_{k},\xi},Z_{s}^{r_{k},\xi})_{r_{k} \le s \leq r_{k+1}}$ solves:
% there exits a unique solution which can be expressed in term of the decoupling field $U$ in the following way:
 \begin{align*}
  \left \{
 \begin{array}{rcl}
 X_t &=&  \xi + \int_{r_{k}}^t b\bigl(X_{s},Y_s,\law{X_s,Y_s}\bigr)\ud s + \int_{r_k}^t \sigma\bigl(X_s,\law{X_{s}}\bigr)\ud W_s \; ,\;\\ %X_{(N-k)\delta} = \xi \in L^2(\cF_{(N-k)\delta})\\
 Y_t &=& \cU\bigl(r_{k+1}, X_{r_{k+1}},\law{X_{r_{k+1}}}\bigr) 
 + \int_{t}^{r_{k+1}} f\bigl(X_{s},Y_{s},Z_s,\law{X_s,Y_s}\bigr)\ud s - \int_t^{r_{k+1}} Z_s \cdot \ud W_s \; ,
   \end{array}
  \right.
   \end{align*}
  for $t \in [r_{k},r_{k+1}]$.  
   
If $\delta$ is small enough, a natural approach is to introduce a Picard iteration scheme to 
approximate the solution of the above equation. To do so, one can  implement the following recursion (with respect to the index $j$):
 \begin{align*}
  \left \{
 \begin{array}{rcl}
 X^{j}_t &=&  \xi + \int_{r_{k}}^t b\bigl(X^j_{s},Y^{j}_s,\law{X^j_s,Y^j_s}\bigr)\ud s +  \int_{r_k}^t \sigma\bigl(X^j_s,
\law{X^j_s} 
 \bigr)\ud W_s\;,\;\\ %X_{(N-k)\delta} = \xi \in L^2(\cF_{(N-k)\delta})\\
 Y^{j}_t &=& \cU\bigl(r_{k+1}, X^{j-1}_{r_{k+1}},\law{X^{j-1}_{r_{k+1}}}\bigr) 
  + \int_{t}^{r_{k+1}} f\bigl(X_{s}^{j-1},Y_{s}^j,Z^j_s,\law{X^{j-1}_s,Y^j_s}\bigr)\ud s -\int_{t}^{r_{k+1}}Z^j_s \cdot \ud W_s
   \end{array}
  \right.\,.
   \end{align*}
 with $(X_{s}^0 = \xi
  + \int_{r_{k}}^t b\bigl(X^0_{s},0,\law{X^0_s,0}\bigr)\ud s +  \int_{r_k}^t \sigma\bigl(X^0_s,
\law{X^0_s} 
 \bigr)\ud W_s )_{r_{k} \leq s \leq r_{k+1}}$ and $(Y^0_{s} = 0)_{r_{k} \le s \le r_{k+1}}$. It is known that, for $\delta$ small enough, $(X^j,Y^j,Z^j) \rightarrow_{j \rightarrow \infty} (X,Y,Z)$, in the sense that
$\vvvert (X^j - X,Y^j - Y,Z^j - Z) \vvvert_{[r_{k},r_{k+1}]}
\rightarrow_{j \rightarrow \infty} 0$.

 \vspace{5pt}
 \noindent But in practice we
 will encounter three main difficulties. 
 \begin{enumerate}
 \item The procedure has to be stopped after a given number of iterations $J$. 
 %\footnote{\textcolor{red}{Actually it could
 %depend on $k$ as well, but we do not consider this option here.}}. 
 \item The above Picard iteration assumes the perfect knowledge of the map $\cU$ at time $r_k$, but
 $\cU$ is exactly what we want to compute... 
 \item The solution has to be discretized in time and space.
 \end{enumerate}
 
\noindent \textbf{Ideal recursion.} 
We first discuss 1) and 2) above. 
The main idea is to use a recursive algorithm (with a new recursion, but on the time parameter). 
\vspace{5pt}
\\
Namely, for $k \le N-1$, we assume that we are given a solver which computes 
\begin{align}\label{eq error solver}
\sol{k+1}{\xi} =  \cU(r_{k+1},\xi,\law{\xi}) + \epsilon^{k+1}(\xi)\,,
\end{align}
where $\epsilon$ is an error made, for any $\xi \in L^2(\cF_{r_{k+1}})$. 
We shall sometimes refer to $\sol{k+1}{\cdot}$ as ``the solver at level $k+1$''.
\vskip 5pt

\noindent Taking these observations into account,
we first define an ideal solver, which assumes that 
each Picard iteration in the approximation of the solution of the forward-backward system can be perfectly computed. We denote it by 
$\pic{}{}$. Accordingly, we identify (for the time being)  
$\sol{k+1}{}$ with 
$\pic{k+1}{}$. Given $\pic{k+1}{}$, $\pic{k}{}$ is defined as follows.
 \begin{align}\label{eq algo solver no error}
  \left \{
 \begin{array}{rcl}
 \tX{k,j}_t &=&  \xi + \int_{r_{k}}^t b\bigl(\tX{k,j}_{s},\tY{k,j}_s,\law{\tX{k,j}_s,\tY{k,j}_s}\bigr)\ud s 
+ \int_{r_k}^t 
 \sigma\bigl(\tX{k,j}_s,
 \law{\tX{k,j}_s}
 \bigr) \ud W_s \;,\;\\ %X_{(N-k)\delta} = \xi \in L^2(\cF_{(N-k)\delta})\\
 \tY{k,j}_t &=& \pic{k+1}{\tX{k,j-1}_{r_{k+1}}} 
 - \int_t^{r_{k+1}} \tZ{k,j}_s \cdot \ud W_s
\\
 &&\hspace{15pt}+ \int_{t}^{r_{k+1}} f\bigl(\tX{k,j-1}_{s},\tY{k,j}_{s},\tZ{k,j}_s,\law{\tX{k,j-1}_s,\tY{k,j}_s}\bigr)\ud s \,,
   \end{array}
  \right.
   \end{align}
for $j \ge 1$ and with 
\begin{center}
$\Bigl(\tilde{X}^{k,0}_{s}= 
\xi
 + \int_{r_{k}}^t b(X^{k,0}_{s},0,\law{X^{k,0}_s,0})\ud s +  \int_{r_k}^t \sigma(X^{k,0}_s,
\law{X^{k,0}_s} )\ud W_s
\Bigr)_{r_{k} \le s \le r_{k+1}}$,
\end{center}
and 
$(\tilde{Y}^{k,0}_{s}=0)_{r_{k} \le s \le r_{k+1}}$. We then define 
\begin{align*}
 \pic{k}{\xi} := Y^{k,J}_{r_{k}} \text{ and } \epsilon^{k}(\xi) :=  Y^{k,J}_{r_{k}} - \cU(r_{k},\xi,\law{\xi})\;,
\end{align*}
where $J \geq 1$ is the number of Picard iterations. 
%\textcolor{red}{say a word about the last solver}\\

At level $N-1$, which is the last level for our recursive algorithm, the Picard iteration scheme is given by
 \begin{align}\label{eq algo solver last level}
  \left \{
 \begin{array}{rcl}
 \tX{N-1,j}_t &=&  \xi + \int_{r_{N-1}}^t b\bigl(\tX{N-1,j}_{s},\tY{N-1,j}_s,\law{\tX{N-1,j}_s,\tY{N-1,j}_s}\bigr)\ud s 
 \\ 
 &&\hspace{15pt}
 + \int^t_{r_{N-1}}\sigma\bigl(\tX{N-1,j}_s,
 \law{\tX{N-1,j}_s}
 \bigr) \ud W_s\;,\;\\ %X_{(N-k)\delta} = \xi \in L^2(\cF_{(N-k)\delta})\\
 \tY{N-1,j}_t &=& g(\tX{N-1,j-1}_T,\law{\tX{N-1,j-1}_T})
  - \int_t^{T} \tZ{N-1,j}_s \cdot \ud W_s\,
 \\ 
 &&\hspace{15pt}+ \int_t^Tf\bigl(\tX{N-1,j-1}_{s},\tY{N-1,j}_{s},\tZ{N-1,j}_s,\law{\tX{N-1,j-1}_{s},\tY{N-1,j}_{s}}\bigr)\ud s\, .
   \end{array}
  \right.
   \end{align}
Here, the terminal condition $g$ is known and the error comes from the fact that the Picard iteration is stopped. It is then natural to set, for $\xi \in L^2(\cF_{T})$,% and $\eta \in L^2(\cF_{T})$.
\begin{align}
%\sol{N-1}{\xi} & =  U(r_{N-1},\xi,\law{\xi}) + \epsilon^{N-1}(\xi)\,, \label{eq error solver last level}
%\\
\pic{N}{\xi} & = g(\xi,\law{\xi}) \text{ and } \epsilon^N(\xi) = 0\,. \label{eq error solver very last level}
\end{align}

\paragraph{Practical implemention.}
As already noticed in 3) above, it is not possible to solve the backward and forward equations in \eqref{eq algo solver no error} perfectly, even though 
the system is decoupled. Hence, we need to introduce
an approximation that can be implemented in practice. Given a continuous adapted input process $\mathfrak{X}=(\mathfrak{X}_{s})_{r_{k} \le s \le r_{k+1}}$ such that $\EE[\sup_{r_{k} \le s \le r_{k+1}} \vert 
{\mathfrak X}_{s} \vert^2]< \infty$
 and $\eta \in L^2(\Omega,\cF_{r_{k+1}},\PP;\RR)$, we thus would like to solve
 \begin{align*} 
 \left \{
 \begin{array}{rcl}
  \tilde{X}_t &= &\mathfrak{X}_{r_k} + \int_{r_k}^t b\bigl(\tilde{X}_{s},\tilde{Y}_s, \law{\tilde{X}_s,\tilde{Y}_s}\bigr) \ud s + \int_{r_k}^t \sigma\bigl(\tilde{X}_s,\law{\tilde{X}_s} 
  \bigr)\ud W_s %\quad \text{ and } \quad \tilde{X}_{(N-k-1)\delta}=  
  \\
  \tilde{Y}_t &= & \eta + \int_{t}^{r_{k+1}} f\bigl(\mathfrak{X}_{s},\tilde{Y}_{s},\tilde{Z}_s,\law{\mathfrak{X}_{s},\tilde{Y}_s}\bigr)\ud s -\int_t^{r_{k+1}} \tilde{Z}_s \cdot \ud W_s\;,
  \end{array}
  \right.
 \end{align*}
 for $t \in [r_{k},r_{k+1}]$. 
 
Let $\pi$ be a discrete time grid of $[0,T]$ such that $\Re \subset \pi$,
\begin{align}\label{eq de disc grid pi}
\pi := \set{t_0 := 0 < \dots < t_n := T} \,\text{ and }\, |\pi| := \max_{i<n}(t_{i+1}-t_i).
\end{align}
For $0 \le k \le N-1$, we note $\pi^k := \set{t \in \pi \,|\, r_k \le t \le r_{k+1}}$
and for later use, we define the indices $(j_k)_{0 \le k\le N}$ as follows
\begin{align*}
\pi^k = \set{t_{j_k} := r_k < \dots<t_i< \dots < r_{k+1} =:t_{{j_{k+1}}}}\,,
\end{align*}
for all $k<N$.
So, instead of a perfect solver for an iteration of the Picard scheme \eqref{eq algo solver no error}, we assume that we are given a numerical solver, denoted by $\tsol{k}{\bar{\mathfrak{X}}}{{\eta}}{f}$, which computes
an approximation  of the process $(\tilde X_{s},\tilde Y_{s},\tilde Z_{s})_{r_{k} \le s \le r_{k+1}}$ on $\pi^k$
for a discretization
$(\bar{\mathfrak{X}}_{t})_{t \in \pi^k}$ of the time 
continuous process $({\mathfrak{X}}_{s})_{r_{k} \le s \le r_{k+1}}$. The \textit{output} is denoted by $(\Xd{}_t,\Yd{}_t,\Zd{}_t)_{t \in \pi^k}$. In parallel, we call \textit{input} the triplet formed by the random variable ${\eta}$, the discrete-time process $(\bar{\mathfrak{X}}_t)_{t \in \pi^k}$ and the driver $f$ of the backward equation. 
In short, the output is what the numerical solver returns after one iteration in the Picard scheme when the discrete input is $(\eta,\bar{\mathfrak X},f)$. 
Pay attention that, in contrast with $b$ and $\sigma$, we shall allow $f$ to vary; this is the rationale for regarding it as an input. However, when the value of $f$ is clear, we shall just regard the input as the pair $(\eta, (\bar{\mathfrak{X}}_t)_{t \in \pi^k})$.

\vspace{5pt} The full convergence analysis, including the discretization error, will be discussed in the next section in the following two cases: first for a generic (or abstract) solver $\tsol{}{}{}$ and second for an explicit solver, as given in the example below. 

\begin{example} \label{ex binomial tree}
This example is the prototype of the solver \textrm{\rm $\tsol{}{}{}$}. 
We consider an approximation of the Brownian motion obtained by
quantization of the Brownian increments. At every time $t \in \pi$, 
we denote by $\bar{W}_{t}$ the value at time $t$ of the discretized Brownian motion. It may expressed as
\begin{equation*}
\bar{W}_{t_{i}} := \sum_{j=0}^{i-1} \Delta \bar{W}_{j},
\end{equation*}
where 
\begin{equation*}
\Delta \bar{W}_{j} := h_{j}^{\frac12} \varpi_{j}, \quad
\varpi_{j}
:=  \Gamma_{d} \Bigl( h_{j}^{-\frac12} \bigl( W_{t_{j+1}}-W_{t_{j}}\bigr) \Bigr),
\end{equation*}
$\Gamma_{d}$ mapping $\RR^d$ onto a finite grid of $\RR^d$. 
Importantly, $\Gamma_{d}$ is assumed to be 
bounded by $\Lambda$ and 
each $\varpi_{j}$ is assumed to be centered and to have the identity matrix as covariance matrix. Of course, this is true if
$\Gamma_{d}$ is of the form
\begin{equation*}
\Gamma_{d}\bigl(w_{1},\cdots,w_{d}\bigr) := 
\bigl( \Gamma_{1}(w_{1}),\cdots,\Gamma_{1}(w_{d})
\bigr), \quad (w_{1},\cdots,w_{d}) \in \RR^d, 
\end{equation*}  
where $\Gamma_{1}$ is a bounded odd function from $\RR$ onto a finite subset of $\RR$ with a normalized second order moment under the standard Gaussian measure. In practice, $\Gamma_{d}$ is intended to take a small number of values. Of course, the typical example is the so-called \textit{binomial approximation}, in which case
$\Gamma_{1}$ is the sign function.

On each interval $[r_k,r_{k+1}]$, given a discrete-time input process $\bar{\mathfrak{X}}$ and a terminal condition
${\eta}$, we thus implement the following scheme (below, $\EE_{t_{i}}$ is the conditional expectation given $\cF_{t_{i}}$):
\begin{enumerate}
 \item For the backward component:
 \begin{enumerate}
  \item Set as terminal condition,  $(\Yd_{t_{j_{k+1}}},\Zd_{t_{j_{k+1}}}) = (\eta,0 )$.
  \item For $j_k \le i < j_{k+1}$, compute recursively
  \begin{align*}
&\Yd_{\ti{}}  =  \EFp{\ti{}}{ \Yd_{\ti{+1}} + (t_{i+1}-t_i)f\bigl(\bar{\mathfrak X}_{t_{i}},\Yd_{\ti{}},\Zd_{\ti{}},\law{\bar{\mathfrak{X}}_{\ti{}},\Yd_{\ti{}}}\bigr) } \, , 
\quad
&\Zd{}_{\ti{}} = \EFp{\ti{}}{\frac{\Delta \bar{W}_i }{\ti{+1}-\ti{}} \Yd{}_{\ti{+1}}}\,.%\; t_{j_k} \le \ti{} < t_{j_{k+1}}\,.
\end{align*}
%where $\Delta \bar{W}_i=\bar W_{t_{i+1}} - \bar{W}_{t_{i}}$.
 \end{enumerate}

 \item For the forward component:
 \begin{enumerate}
\item Set as initial condition, $\Xd_{t_{j_{k}}} = \bar{\mathfrak{X}}_{r_k}$.
\item For $j_k < i  \le j_{k+1}$, compute recursively
 \end{enumerate}
\begin{align*}
 \Xd_{\ti{+1}} & = \Xd_{\ti{}} +b\bigl(\Xd_{\ti{}},\Yd_{\ti{}},\law{\Xd_{\ti{}},\Yd_{\ti{}}}\bigr)(\ti{+1}-\ti{})+ \sigma\bigl(\Xd_{\ti{}},\law{\Xd_{\ti{}}}\bigr)\Delta \bar{W}_i\;.
\end{align*}
 
\end{enumerate}
% 
% \begin{align*}
%  \Xd_{t_{j_{k}}} &= \bar{\mathfrak{X}}_{r_k}\, \text{ and } \Yd_{t_{j_{k+1}}} = \bar{\chi} ,
% \\
% \Xd_{\ti{+1}} & = \Xd_{\ti{}} +b(\Yd_{\ti{}},\law{\Xd_{\ti{}},\Yd_{\ti{}}})h+ \sigma(\Xd_{\ti{}})\Delta \bar{W}_i \,,
% %\; t_{j_k} \le \ti{} < t_{j_{k+1}}\,, 
% \\
% \Yd_{\ti{}} & =  \dEFp{\ti{}}{ \Yd_{\ti{+1}} + (t_{i+1}-t_i)f(\Zd_{\ti{}},\law{\bar{\mathfrak{X}}_{\ti{}},\Yd_{\ti{}}}) } \text{ and } 
% \Zd{}_{\ti{}} = \dEFp{\ti{}}{\frac{\Delta \bar{W}_i }{\ti{+1}-\ti{}} \Yd{}_{\ti{+1}}}\,,\; t_{j_k} \le \ti{} < t_{j_{k+1}}\,.
% \end{align*}
\end{example}

\vspace{5pt}

\paragraph{Full algorithm  for $\sol{}{}$.} Using $\tsol{}{}{}$, for each level, we can now give a completely implementable algorithm for $\sol{}{}$. Its description is as follows.

\vspace{5pt}

The value $\sol{k}{\xi}$, i.e. the value of the solver at level $k$ with initial condition $\xi \in L^2(\cF_{r_k})$, is obtained through:
\begin{enumerate}
 \item Initialize the backward component at $\bar{Y}^{k,0}_t=0$ for $t \in \pi_k$
 and regard
 $(\bar{X}^{k,0}_{t})_{t \in \pi_{k}}$
 as the forward component of 
 $\tsol{k}{\xi}{0}{0}$ 
 \item for $ 1 \le j \le J$ 
 \begin{enumerate}
 \item compute $\bar{Y}^{k,j}_{r_{k+1}}=\sol{k+1}{\bar{X}^{k,j-1}_{r_{k+1}}}$.
 \item compute $(\bar{X}^{k,j},\bar{Y}^{k,j},\bar{Z}^{k,j})=\tsol{k}{\bar{X}^{k,j-1}}{\bar{Y}^{k,j}_{r_{k+1}}}{f}$
 \end{enumerate}
 \item return $\bar{Y}^{k,J}_{r_{k+1}}$.
\end{enumerate}
Following 
\eqref{eq error solver very last level}, we let
\begin{align}
\sol{N}{\xi} & = g(\xi,\law{\xi}) \,. \label{eq error solver very last level solver}
\end{align}

We first explain the initialization
step. The basic idea is to set the backward component 
to $0$ and then to solve the forward component as an approximation of the autonomous McKean-Vlasov diffusion process in which the backward entry is null. Of course, this may be solved by means of any standard method, but to make the notation shorten, we felt better to regard the underlying solver 
as a specific case of a forward-backward solver with null coefficients in the backward equation.  We specify in the analysis below the conditions that this initial solver 
$\tsol{}{}{0}{0}$ must satisfy.

It is also worth noting that each Picard iteration used to define the solver at level $k$ calls the solver at level $k+1$. This is a typical feature of the way the continuation method manifests from the algorithmic point of view. 
In particular, the total complexity is of order $O(J^N {\mathfrak K})$, where 
${\mathfrak K}$ is the complexity of the 
solver
$\tsol{}{}{}$. In this regard, it must be stressed that, for a given length $T$, $N$ is fixed, regardless of the time step $\vert \pi \vert$. Also, 
$J$ is intended to be rather small as the Picard iterations are expected to converge geometrically fast, see the numerical examples in Section 
\ref{se numerics} in which we choose $J=5$. However, it must be noticed that 
the complexity increases exponentially fast when $T$ tends to $\infty$, which is obviously the main drawback of this method. Again, we refer to 
Section \ref{se numerics} for numerical illustrations. 
\vskip 5pt

\paragraph{Useful notations.} 
Throughout the paper, $\| \cdot \|_{p}$ denotes the $L^p$ norm on $(\Omega,\cF,\PP)$. Also, $(\hat{\Omega},\hat{\cF},\hat{\PP})$
stands for a copy of $({\Omega},{\cF},{\PP})$. It is especially useful to represent the Lions' derivative of a function of a probability measure and to distinguish the (somewhat artificial) space used for representing these derivatives from the (physical) space carrying the Wiener process. For a random variable 
$X$ defined on $(\Omega,\cF,\PP)$, we shall denote by $\cc{X}$ its copy on $(\hat{\Omega},\hat{\cF},\hat{\PP})$.
\vskip 5pt

\noindent We shall use the notations 
$C_{\Lambda},c_{\Lambda}$ for constants 
only depending on $\Lambda$ (and possibly on the dimension as well). They are allowed to increase from line to line. We shall use the notation $C$ for constants not depending upon the discretization parameters. Again, they are allowed to increase from line to line. In most of the proofs, we shall just write $C$ for $C_{\Lambda}$, even if we use the more precise notation $C_{\Lambda}$ in the corresponding statement.

\subsection{A first analysis with no discretization error}%- no error on $\tsol{}{}{}$}

\vspace{2mm}
To conclude this section, we  want to understand how the error propagates through the solvers used at different levels in the ideal case where 
the Picard iteration in 
\eqref{eq algo solver no error} can be perfectly computed or equivalently 
when the solver is given by $\sol{k}{}{}=\pic{k}{}${}.
%
%
%. On each level, the error is due to the discretization procedure and the fact that the Picard iteration have to be stopped. We will first consider that no error is made on
%$\tsol{k}{\cdot}{\cdot}$ or equivalently that \eqref{eq algo solver no error} is computed perfectly.
%\\ \vspace{2mm}
For $j \ge 1$, we then denote by $(\tX{k,j},\tY{k,j},\tZ{k,j})$, the solution on $[r_k,r_{k+1}]$ of \eqref{eq algo solver no error}.
% with the initial condition $\xi$
% and the terminal condition $Y^{j-1,k}_{(N-k)\delta}$.

The main result of the section, see Theorem \ref{th error picard}, is an upper bound for the error when we use $\pic{\cdot}{\cdot}$ to approximate $\cU$.
The proof of this theorem requires the following proposition, which gives a local error estimate for each level.

\begin{proposition} \label{pr loc error picard}
Let us define, for $j \in \{1,\cdots,J\}$, $k \in \{1,\cdots,N-1\}$, 
$$\Delta_k^j := 
\Bigl\| \sup_{t \in [r_k,r_{k+1}]}\bigl(\tY{k,j}_t - \cU(t,\tX{k,j}_t,\law{\tX{k,j}_t}) \bigr)
\Bigr\|_2$$
then, there exist constants $C_{\Lambda},c_{\Lambda}$ such that, for 
$\bd := C_{\Lambda} \delta < c_{\Lambda}$, 
\begin{align}\label{eq pr local error}
 \Delta_k^j \le \bd^j \Delta^0_k + \sum_{\ell=1}^{j} \bd^{\ell-1} e^{\bd}
 %\norm{ \epsilon^{k+1}(\tX{k,j-\ell}_{r_{k+1}}) }
 \bigl\|\epsilon^{k+1}(\tX{k,j-\ell}_{r_{k+1}}) \bigr\|_2 \, .
\end{align}
\end{proposition}
We recall that $\epsilon^{k}(\xi)$ stands for the error term:
\begin{align*}
\epsilon^{k}(\xi) =   \pic{k}{\xi} - \cU(r_{k},\xi,\law{\xi}) \; , 
\quad \textrm{with} \ \epsilon^{N}(\xi) = 0 \; .
\end{align*}

\begin{remark}\label{re z meme pas peur}
A careful inspection of the proof shows that, whenever $\sigma$ depends on $Y$ or $b$ depends on $Z$, the same result holds true but with a constant 
$C_{\Lambda}$ depending on $N$. As $N$ is fixed in practice, this might still suffice to complete the analysis of the discretization scheme in that more general setting.
\end{remark}

\proof We suppose that the full algorithm is initialized at some level $k \in \{0,\cdots,N-1\}$,  with an initial condition $\xi \in L^2(\cF_{r_k})$. As the value of the index $k$ is fixed throughout the proof, we will drop it in the notations 
$(\tilde X^{k,j},\tilde Y^{k,j},\tilde Z^{k,j})$ and $\Delta_{k}^j$.  
%and set, for $X \in L^2$,
%$$ \norm{X} := \esp{\left| X \right|^2 \Big| \cF_{r_k}}^\frac12 \,.$$
\vspace{5pt}
\\
Applying Ito's formula for functions of a measure argument, see
 \cite{buclij14, chacri15}, we have
\begin{align*}
\ud \cU(t,\tX{j}_t,\law{\tX{j}_t}) 
 & = \bigg(  b(\tX{j}_{t},\tY{j}_t,\law{\tX{j}_t,\tY{j}_t}) \cdot \partial_x \cU(t,\tX{j}_t,\law{\tX{j}_t})
 \\
&\hspace{15pt} 
  + \frac12 {\rm Tr}
 \bigl[ a\bigl(\tX{j}_t,\law{\tX{j}_t}\bigr) \partial^2_{xx} \cU(t,\tX{j}_t,\law{\tX{j}_t}) \bigr]
 \\
 &\hspace{15pt} + \hat{\EE} \Bigl[ b(\cc{\tX{j}_t},\cc{\tY{j}_t},\law{\tX{j}_t,\tY{j}_t}) 	\cdot \partial_\mu \cU(t,\tX{j}_t,\law{\tX{j}_t}) \Bigr] 
\\
&\hspace{15pt}   + \hat{\EE} \Bigl[ \frac12 {\rm Tr} \bigl[ 
a\bigl(\cc{\tX{j}_t},\law{\tX{j}_t}\bigr) \partial_\upsilon \partial_\mu \cU(t,\tX{j}_t,\law{\tX{j}_t}) \bigr]  
\Bigr] + \partial_t \cU(t,\tX{j}_t,\law{\tX{j}_t}) \biggr)\ud t 
\\
&\hspace{5pt}+ \partial_x \cU(t,\tX{j}_t,\law{\tX{j}_t}) \cdot \bigl( \sigma\bigl(\tX{j}_t,\law{\tX{j}_t}\bigr) \ud W_t \bigr) \, .
 \end{align*}
 Expressing the integral in \eqref{eq pde decoupling field} as expectations on $(\hat{\Omega},\hat{\cF},\hat{\PP})$ and combining with \eqref{eq pde decoupling field} and \eqref{eq algo solver no error}, we obtain
\begin{align*} %\label{eq analysis I 1}
 \ud [\check{Y}^j_t - \tY{j}_t] 
 &=\Big(  \bigl\{b\bigl(\tX{j}_t,\tY{j}_t,\law{\tX{j}_t,\tY{j}_t}\bigr) - b\bigl(\tX{j}_t,\check{Y}^j_t,\law{\tX{j}_t,\check{Y}^j_t}\bigr)\bigr\}
 \cdot 
 \partial_x \cU(t,\tX{j}_t,\law{\tX{j}_t})
  %\nonumber 
\\
&\hspace{5pt} + \whesp{
\bigl\{b\bigl(\cc{\tX{j}_t},\cc{\tY{j}_t},\law{\tX{j}_t,\tY{j}_t}\bigr) - b\bigl(\cc{\tX{j}_t},\cc{\check{Y}^j_t},\law{\tX{j}_t,\check{Y}^j_t}\bigr)\bigr\}
\cdot 
\partial_\mu \cU(t,\tX{j}_t,\law{\tX{j}_t})
 }
%\nonumber
 \\
 &\hspace{5pt} + f\bigl(\tX{j}_t,\tY{j}_t,\tZ{j}_t,\law{\tX{j}_t,\tY{j}_t}\bigr) - f\bigl(\tX{j}_t,\check{Y}^j_{t},\check{Z}^j_t,\law{\tX{j}_t,\check{Y}^j_t}\bigr) \Big) \ud t
+ [\check{Z}^j_t - \tZ{j}_t] \cdot \ud W_t\,, %\nonumber
\end{align*}
where $\check{Y}^j_t := \cU(t,\tX{j}_t,\law{\tX{j}_t}) $ and $\check{Z}^j_t := \partial_x u(t,\tX{j}_t,\law{\tX{j}_t}) \sigma(\tX{j}_t,\law{\tX{j}_t})$. Observe that this argument is reminiscent of the four-step scheme, 
see \cite{MaProtterYong}. 
\vspace{5pt}
\\
Using standard arguments from BSDE theory and \HYP{0}--\HYP{1}, we then compute
\begin{align*}
 \Delta^j &\le e^{C\delta} \bigl\| \cU(r_{k+1},\tX{j}_{r_{k+1}},\law{\tX{j}_{r_{k+1}}}) - \tY{j}_{r_{k+1}} \bigr\|_{2}
 \\
 &\le e^{C \delta} \Bigl( \bigl\| \epsilon^{k+1}(\tX{j-1}_{r_{k+1}}) 
 \bigr\|_{2} + \bigl\| \cU(r_{k+1},\tX{j}_{r_{k+1}},\law{\tX{j}_{r_{k+1}}}) - \cU(r_{k+1},\tX{j-1}_{r_{k+1}},\law{\tX{j-1}_{r_{k+1}}})
\bigr\|_{2} \Bigr),
\end{align*}
recalling $\tY{j}_{r_{k+1}} = \pic{k+1}{\tX{j-1}_{r_{k+1}}}$ and  \eqref{eq error solver}. Since $\cU$ is Lipschitz, we have
\begin{align}\label{eq analysis I dY}
 \Delta^j
 &\le e^{C \delta} \left( \bigl\| \epsilon^{k+1}(\tX{j-1}_{r_{k+1}}) \bigr\|_{2} + 2L \bigl\| 
 \tX{j}_{r_{k+1}} - \tX{j-1}_{r_{k+1}} \bigr\|_{2} \right)\;.
\end{align}
We also have that
\begin{align*}
 \tX{j}_t - \tX{j-1}_t  = \int_{r_k}^t \bigl\{ b\bigl(
 \tX{j}_{s},\tY{j}_s,\law{\tX{j}_t,\tY{j}_t}\bigr) & - b\bigl(\tX{j-1}_s,\tY{j-1}_s,\law{\tX{j-1}_t,\tY{j-1}_t}\bigr)\bigr\} \ud s
 \\
 &\hspace{-40pt} + \int_{r_k}^t \bigl\{ \sigma\bigl(\tX{j}_s,
\law{\tX{j}_s}  \bigr) - \sigma\bigl(\tX{j-1}_s,
\law{\tX{j-1}_s}
\bigr)\bigr\} \ud W_s \;.
\end{align*}
Using usual arguments (squaring, taking the sup, using B\"urkholder-Davis-Gundy inequality), we get, since $b$ and $\sigma$ are Lipschitz continuous,
\begin{align*}
&\Bigl\| \sup_{t\in [r_k,r_{k+1}] } \vert \tX{j}_t - \tX{j-1}_t \vert 
 \Bigr\|_{2}
 \le C \Bigl( \delta 
 \Bigl\| 
 \sup_{t\in [r_k,r_{k+1}] } \vert \tY{j}_t - \tY{j-1}_t \vert
 \Bigr\|_{2}
 + \delta^{\frac12} 
 \Bigl\| 
 \sup_{t\in [r_k,r_{k+1}] } \vert \tX{j}_t - \tX{j-1}_t \vert
 \Bigl\|_{2}  \Bigr) \, .
\end{align*}
Observing that
\begin{align*}
 |\tY{j}_s - \tY{j-1}_s| & \le |\tY{j}_s - \cU(s,\tX{j}_s,\law{\tX{j}_s})| + |\tY{j-1}_s - \cU(s,\tX{j-1}_s,\law{\tX{j-1}_s})| 
  \\
  &\hspace{15pt}+ \Lambda(|\tX{j-1}_s-\tX{j}_s|+ \|\tX{j-1}_s-\tX{j}_s \|_2) \, ,
\end{align*}
we obtain, for $\delta$ small enough,
\begin{align} \label{eq analysis I dX}
 \Bigl\| \sup_{t\in[r_k,r_{k+1}]} \vert \tX{j}_t - \tX{j-1}_t \vert \Bigr\|_{2} \le C \delta (\Delta^j + \Delta^{j-1})\, .
\end{align}
Combining the previous inequality with \eqref{eq analysis I dY}, we obtain, for $\delta$ small enough,
\begin{align*}
 \Delta^j &\le e^{C\delta}\bigl\|\epsilon^{k+1}(\tX{j-1}_{r_{k+1}}) \bigr\|_2 + C\delta \Delta^{j-1}\, ,
 \end{align*}
 which by induction leads to
 \begin{align*}
 \Delta^j 
 &\le (C\delta)^j \Delta^0 + \sum_{\ell=1}^{j} (C\delta)^{\ell-1} e^{C\delta}\bigl\|\epsilon^{k+1}(\tX{j-\ell}_{r_{k+1}}) \bigr\|_2 \, ,
\end{align*}
and concludes the proof.
\eproof

\vspace{10pt} 

We now state the main result of this section, which explains how
the local error induced by the fact that the Picard iteration is stopped at rank $J$ propagates through the various levels $k=N-1,\cdots,0$.

\begin{theorem}\label{th error picard}
We can find two constants $C_\Lambda,c_{\Lambda}>0$
and a continuous non-decreasing function $\mathfrak{B} : \RR_{+} \rightarrow \RR_{+}$ matching $0$ in $0$, only depending on $\Lambda$,  
 such that, for $\bd:=C_{\Lambda} \delta < \min(c_{\Lambda},1)$
 and $\beta \geq \mathfrak{B}(\bd)$ satisfying 
\begin{align} \label{eq main ass}
(J-1)  \Lambda  \bd^J \frac{ e^{\beta C_\Lambda T}}{e^{\beta \bd}-1}  \le 1 
\end{align}
where  $J$ is the number of Picard iterations in a period, it holds, for any period $k\in \{0,\cdots, N\}$ and $\xi \in L^2(\cF_{r_k})$, 
  {\rm
 \begin{align} \label{eq control picard only}
  \norm{\sol{k}{\xi} - \cU(r_k,\xi,\law{\xi})} \le \Lambda
  \frac{e^{\beta C_{\Lambda} T}}{\beta}  
   \bd^{J-1}\bigl(1+\norm{P_{r_{k},T}^\star(\xi)}\bigr) \; ,
 \end{align}
 }
 where $P_{r_{k},t}(\xi)$ is the solution at time $t$ of the stochastic differential equation
 \begin{equation*}
d X_{s}^0  = b\bigl(X_{s}^0,0,\law{X_{s}^0,0}\bigr) \ud s + 
\sigma \bigl( X_{s}^0,\law{X_{s}^0} \bigr) \ud W_{s} \, ,
\end{equation*}
with $X_{r_{k}}^0=\xi$ as initial condition, and 
$P_{r_{k},t}^\star(\xi) = \sup_{s \in [r_{k},t]} \vert P_{r_{k},s}(\xi) \vert$.  
%provided that $\bd$ is small enough.
\end{theorem}
Of course, it is absolutely straightforward 
to bound 
$\norm{P_{r_{k},T}^\star(\xi)}$ by $C  ( 1 + \norm{\xi})$ in 
\eqref{eq control picard only}. Theorem
\ref{th error picard} may be restated accordingly, but the form used in the statement is more faithful to the spirit of the proof.

\vskip 5pt

\proof
We prove the claim by an induction argument. 
We show below that for all $k \in \{0,\dots,N\}$,
\begin{align} \label{eq proof hyp rec}
 \norm{\epsilon^k(\xi)} = \norm{\sol{k}{\xi} - \cU(r_k,\xi,\law{\xi})} \le \theta_{k}\left(1+\norm{P_{r_{k},T}^\star(\xi)}\right) \, ,
\end{align}
where $(\theta_k)_{k=0,\cdots, N-1}$ is defined by the following backward induction: $\theta_N := 0$, recall \eqref{eq error solver very last level}, and for $k \in \{0,\cdots,N-1\}$,
\begin{align} \label{eq de theta_k}
\theta_k := \Lambda \bd^J + e^{\beta \bd} \theta_{k+1} \,, 
% \theta_k := \Lambda \bd^J + \gamma \theta_{k+1} + 2\gamma \theta_{k+1} %\bd e^{2 \gamma \theta_{k+1}(J-1) \bd}\left( \gamma(J-1) \theta_{k+1} + %\frac{\Lambda}{1-\bd}\right)
\end{align}
where $\beta$ is such that
\begin{align}\label{eq de beta}
 \left(\gamma +  \gamma \bd e^{\gamma \bd } (  \gamma  + \frac{\Lambda}{1-\bd}) \right) \le e^{\beta \bd} \;,\quad 
 \text{ with }
 \gamma := \frac{ e^{\bd} }{1-\bd} \;.
\end{align}
With this definition, 
we have, for all $k \in \{0,\cdots, N\}$, 
\begin{align} \label{eq bound theta}
\theta_k =
\Lambda \bd^J
\sum_{j=0}^{N-k-1} e^{j \beta \bar{\delta}}
\leq \Lambda \bd^J\frac{e^{\beta C_\Lambda T}}{e^{\beta \bd} - 1} \, ,
\end{align}
which gives the expected result. 
\vspace{5pt} 
\\
We now prove 
\eqref{eq proof hyp rec}. Observe that it is obviously true for the last step $N$. 
Assume now that it holds true at step $k+1$, for $k<N$, and that \eqref{eq bound theta} holds true for $\theta_{k+1}$.
Then, using \eqref{eq main ass}, we have
\begin{align} \label{eq control theta k+1}
 \theta_{k+1}  j \le 1, \quad \text{ for all } j \le J-1\;.
\end{align}
From Proposition \ref{pr loc error picard}, we have
\begin{align}
\Delta^j_k \le \bd^j \Delta^0_k  
+ \sum_{\ell=1}^j \bd^{\ell-1} e^{\bd}\|\epsilon^{k+1}(\tX{k,j-\ell}_{r_{k+1}})\|_2\,.
\end{align}
%note for later use that $\Delta^0_k \le \Lambda(1+\norm{\xi})$. 
Using the induction hypothesis \eqref{eq proof hyp rec}, we compute
\begin{align}\label{eq temp init}
\Delta^j_k \le \bd^j \Delta^0_k  
+ \frac{e^{\bd} }{1-\bd}\theta_{k+1} +  e^{\bd} \theta_{k+1} \sum_{\ell=0}^{j-1} \bd^{j-1-\ell} 
\bigl\|P_{r_{k+1},T}^\star\bigl(\tilde{X}^{k,\ell}_{r_{k+1}}\bigr) \bigr\|_2 \ . 
\end{align} 
We study the last sum. 
Observe that for $\ell \in \{1,\cdots,j-1\}$, 
\begin{equation*}
\begin{split}
\bigl\|P_{r_{k+1},T}^\star\bigl(
\tilde{X}^{k,\ell}_{r_{k+1}} \bigr) \bigr\|_2 &\le 
\bigl\| P_{r_{k+1},T}^\star \bigl( 
\tilde{X}^{k,0}_{r_{k+1}} \bigr) \bigr\|_2 
  + \sum_{i=1}^\ell \norm{
P_{r_{k+1},T}^\star
\bigl(
\tilde{X}^{k,i}_{r_{k+1}}\bigr)-
P_{r_{k+1},T}^\star\bigl(
\tilde{X}^{k,i-1}_{r_{k+1}} \bigr)}.
 \end{split}
 \end{equation*} 
 We observe that 
 $P_{r_{k+1},t}( 
\tilde{X}^{k,0}_{r_{k+1}}) = P_{r_{k},t}(\tilde{X}^{k,0}_{r_{k}})= P_{r_{k},t}(\xi)$, for $t \in [r_{k+1},T]$. 
Hence,  $P_{r_{k+1},T}^\star( 
\tilde{X}^{k,0}_{r_{k+1}})
\leq 
P_{r_{k},T}^\star(\xi)$.
Also, it is well-checked that there exists a constant $C_{\Lambda}$
such that each $P_{t,T}^\star$
is $C_{\Lambda}$-Lipschitz continuous from $L^2(\cF_{t})$ into 
$L^2(\cF_{T})$. 
  Then,
\begin{align*}
 &\sum_{\ell=0}^{j-1} \bd^{j-1-\ell} \bigl\|P_{r_{k},T}^\star\bigl( \tilde{X}^{k,\ell}_{r_{k+1}} \bigr) \bigr\|_{2} 
 \le C_{\Lambda} \sum_{\ell=1}^{j-1} \bd^{j-1-\ell} \sum_{i=1}^\ell \norm{\tilde{X}^{k,i}_{r_{k+1}} -\tilde{X}^{k,i-1}_{r_{k+1}} }
 + \sum_{\ell=0}^{j-1}\bd^\ell \bigl\|P_{r_{k},T}^\star(\xi)\bigr\|_2 \, .
\end{align*}
Using \eqref{eq analysis I dX} in the proof of Proposition \ref{pr loc error picard} and changing the definition of $\bd$, we obtain
\begin{equation}
\label{autre equation 31 decembre}
\begin{split}
& \sum_{\ell=0}^{j-1} \bd^{j-1-\ell} \bigl\|
P_{r_{k},T}^\star\bigl( \tilde{X}^{k,\ell}_{r_{k+1}}\bigr) \bigr\|_2 
\le 
 \bd  \sum_{i=1}^{j-1}  (\Delta^i_k + \Delta^{i-1}_k )\sum_{\ell=i}^{j-1} \bd^{j-1-\ell} 
 + \sum_{\ell=0}^{j-1}\bd^\ell \bigl\|P_{r_{k},T}^\star(\xi)
 \bigr\|_2 \, . 
\end{split}
\end{equation}
Observing that, for all $i \le j-1$, 
$ 
\sum_{\ell=i}^{j-1} \bd^{j-1-\ell}  
%= \sum_{\ell=1}^{j-i} (\bd)^{\ell} 
\le \frac{1}{1-\bd},
$
we  get
\begin{align}
 \sum_{\ell=0}^{j-1} \bd^{j-\ell} \bigl\|
 P_{r_{k},T}^\star\bigl(\tilde{X}^{k,\ell}_{r_{k+1}}\bigr) \bigr\|_2 
 \le 
 \frac{ 2\bd }{1-\bd} \cS^{j-1}_k + \frac{1 }{1-\bd} \bigl\|
 P_{r_{k},T}^\star(\xi)\bigr\|_2 \, ,
\end{align}
where $\cS^{n}_k := \sum_{i=0}^n\Delta^i_k$.
Inserting the previous estimate into \eqref{eq temp init} and changing $\bd$ into $2 \bd$, we obtain
\begin{align}\label{eq temp paused}
\Delta^j_{k} \le \bd^j \Delta^0_{k}  
+ \frac{ e^{\bd} }{1-\bd}\theta_{k+1}
\bigl(1+
\bigl\|P_{r_{k},T}^\star(\xi)\bigr\|_2
\bigr) 
+   \theta_{k+1} \frac{ \bd e^{\bd}}{1-\bd} \cS^{j-1}_k \, .
%+ e^{\bd} \theta_{k+1}(\bd)^j |\xi|
\end{align} 
% We define
% \begin{align*}
% \alpha = \frac{2 \bd e^{\bd} }{1-\bd}\;,
% \end{align*}
% that will be set later on but we already impose that $\alpha \le 1$.
We note that $\Delta^0_k \le \Lambda(1+ \|P_{r_{k},T}^\star(\xi) \|_{2})$. Recalling $\gamma$ in \eqref{eq de beta}, equation
\eqref{eq temp paused} leads to
\begin{align}\label{eq temp cntd}
\Delta^j_k \le  a_j
%(\Lambda \bd^j  
%+ \gamma \theta_{k+1})(1+\norm{\xi}) 
+   \gamma \theta_{k+1} \bd  \cS^{j-1}_k  \;.
%+ \theta_{k+1}\bd^j e^{\bd}  |\xi|
\end{align} 
where we set
$ 
a_j := (\Lambda \bd^j +\gamma \theta_{k+1}) (1+  \|
P_{r_{k},T}^\star(\xi) \|_{2}). 
$
We  have
\begin{align*}
 \cS^j_k - \cS^{j-1}_k = \Delta^j_k \le 
a_j
+   \gamma \theta_{k+1} \bd \cS^{j-1}_k \, , 
\end{align*}
 and then 
 \begin{align} \label{eq cSj}
\cS^j_k \le e^{\gamma \theta_{k+1}\bd j}\cS^0_k +  \sum_{\ell = 1}^{j} e^{ \gamma \theta_{k+1} \bd (j-\ell)} a_{\ell}\;.
\end{align}
We  compute
\begin{align*}
\sum_{\ell = 1}^j a_\ell \le \Bigl(j\gamma \theta_{k+1} 
+ \frac{\Lambda\bd}{1-\bd} \Bigr)
\bigl(1+\norm{P^\star_{r_{k},T}(\xi)} \bigr) \, ,
\end{align*}
which combined with the properties \eqref{eq control theta k+1} and \eqref{eq cSj} leads to, for all $j \le J-1$,
\begin{align*}
\cS^j_k \le e^{  \gamma \bd } \left(  \gamma  + \frac{ \Lambda}{1-\bd}\right)\bigl(1+
\norm{
P^\star_{r_{k},T}(\xi)}
\bigr) \, ,
\end{align*}
where we recall that $\cS^0_k = \Delta^0_k \le \Lambda(1+\|P^\star_{r_{k},T}(\xi)\|_{2})$.
We insert the previous inequality into \eqref{eq temp cntd} for $j=J$ and get
\begin{align*}
 \Delta^J_k \le  \left( \Lambda \bd^J + \left(\gamma +  \gamma \bd e^{ \gamma \bd } (  \gamma  + \frac{\Lambda}{1-\bd}) \right)\theta_{k+1}  \right) \left(1+  \norm{P^\star_{r_{k},T}(\xi)}\right) \, .
\end{align*}
Using \eqref{eq de beta}, this rewrites
\begin{align*}
  \Delta^J_k \le \left(  \Lambda \bd^J +  e^{\beta \bd} \theta_{k+1}  \right) \left(1+  \norm{P^\star_{r_{k},T}(\xi)}\right) \, ,
\end{align*}
and validates \eqref{eq de theta_k} and thus \eqref{eq bound theta}. 
We then obviously have that \eqref{eq proof hyp rec} holds true. 
%\\Equation \eqref{eq control picard only}
%is simply deduced from the upper bound \eqref{eq bound theta}.
\eproof

\section{Convergence Analysis}
\label{se algo MKVFBSDE analysis}
\subsection{Error analysis in the generic case}

We now study the convergence of a generic implementable solver 
$\sol{}{}$, based upon the local solver $\tsol{}{}{}$ as described above, 
as long as %We consider a discrete measure $\bar{\P}$ approximating the Wiener measure (and denote $\desp{\cdot}$ the associated
%expectation operator).\\
%Importantly, $\tsol{}$ does not allow a perfect computation of the solution to \eqref{eq de fbsde tsol}. 
the output of the local solver  $\tsol{k}{}{}$ satisfies
some conditions, which are shown to be true for Example 
\ref{ex binomial tree}. 

In order to define the required assumption, we use the same letters $\Lambda$ and $\alpha$ as in \HYP{0} and \HYP{1},
except that, without any loss of generality, 
 we assume 
 that  
$\alpha$ is greater than 1. 
For the same coefficients as in the
equation 
\eqref{eq fbsde}, and in particular for the same driver $f$, 
we then ask $\tsol{k}{}{}{}$ to satisfy the following three conditions. 
\begin{align*}
(A1)& \quad \sup_{t \in \pi^k} \normalpha{\cU(t,\Xd_t,\law{\Xd_t})-\Yd_t} \le e^{\Lambda \delta} \normalpha{\cU(r_{k+1},\Xd_{r_{k+1}},\law{\Xd_{r_{k+1}}})-\Yd_{r_{k+1}}} 
\\
&\quad\quad
\quad\quad\quad\quad\quad\quad
\quad\quad\quad 
+ \Lambda \max_{j_k \le i < j_{k+1}}\normalpha{\Xd_{t_i} - \bar{\mathfrak{X}}_{t_i}}
+\mathcal{D}^1(|\pi|) +\mathcal{D}^2(|\pi|)\bigl(1+
 \normalphas{\xi}^{\alpha}\bigr)\,,
\\
(A2)& \quad 
 \sup_{t \in \pi^k} \normalpha{\Xd_t-\Xd'_t} \le \Lambda \delta \sup_{t \in \pi^k}\normalpha{\Yd_t-\Yd'_t}\,,
\\
(A3)& \quad  \normalpha{\cU(r_{k+1},\Xd_{r_{k+1}},\law{\Xd_{r_{k+1}}})-\Yd_{r_{k+1}}}^{\alpha} 
\le \Lambda  \normalpha{\cU(r_{k+1},\Xd_{r_{k+1}},\law{\Xd_{r_{k+1}}})-\Yd_{r_{k+1}}}\, ,
\end{align*}
where $(\Xd{},\Yd{},\Zd{}):=\tsol{k}{\bar{\mathfrak{X}}}{\eta}{f}$, for $f$ as before, and $(\Xd{'},\Yd{'},\Zd{'}):=\tsol{k}{\bar{\mathfrak{X}}'}{\eta'}{f'}$, for another $f'$ either equal to $f$ or $0$, are two output values of  $\tsol{}{}{}$ associated to two input processes $\bar{\mathfrak{X}}$, $\bar{\mathfrak{X}}'$, with the same initial condition $\bar{\mathfrak{X}}_{r_k}=\bar{\mathfrak{X}}'_{r_k}=\xi$, and to two different terminal conditions $\eta$ and $\eta'$. 
For $i \in \set{1,2}$, the function $\mathcal{D}^i:[0,\infty) \rightarrow [0,\infty)$ is a discretization error associated to the use of the grid $\pi$, which satisfies $\lim_{h\downarrow 0} \mathcal{D}^i(h) = 0$. Importantly, both ${\mathcal D}^1$ and ${\mathcal D}^2$ are independent of $\bar{\mathfrak{X}}$, $\bar{\eta}$, $J$ and $N$.

In full analogy with the discussion right below 
Theorem 
\ref{th error picard}, we shall also need some conditions on the 
solver $\tsol{k}{}{0}{0}$ used to initialize the algorithm at each step.
Following the definition of  $(P_{r_{k},t})_{0 \le t \le T}$
introduced in the statement of Theorem 
\ref{th error picard}, we let by induction, for a given $k \in \{0,\cdots,N-1\}$:
%\textcolor{magenta}{
%\begin{equation*}
%\begin{split}
%&{\tt P}_{r_{k},t}(\xi) = \bigl( \tsol{N-1}{\xi}{0}{0} \bigr)_{t}, \quad t \in \pi^{k} \, , \quad \xi \in L^2(\cF_{r_{k}}) \, ,
%\end{split}
%\end{equation*}
%and 
%\begin{equation*}
%\begin{split}
%&{\tt P}_{r_{k},t}(\xi) =  {\tt P}_{r_{k+1},t}
%\bigl( (\tsol{k}{\xi}{0}{0} )_{r_{k+1}} \bigr) \ , 
%\quad t \in \pi, \ t \geq r_{k+1},
%\end{split}
%\end{equation*}
%}
\begin{equation*}
\begin{split}
&{\tt P}_{r_{k},t}(\xi) = \bigl( \tsol{k}{\xi}{0}{0} \bigr)^1_{t}, \quad t \in \pi^{k} \, , \quad \xi \in L^2(\cF_{r_{k}}) \, ,
\end{split}
\end{equation*}
where we 
recall that $\bigl( \tsol{k}{\xi}{0}{0} \bigr)^1$ is the forward component of the algorithm's output,
and, for $k \le N-2$,
\begin{equation*}
\begin{split}
&{\tt P}_{r_{k},t}(\xi) =  {\tt P}_{r_{\ell},t}
\bigl( {\tt P}_{r_{k},r_\ell}(\xi) \bigr), 
\quad t \in {\pi^\ell}, \quad  k< \ell \le N-1 ,%\ t \geq r_{k+1},
\end{split}
\end{equation*}
and then ${\tt P}_{r_{k},T}^\star(\xi)=\max_{s \in \pi, s \in [r_{k},T]} \vert {\tt P}_{r_k,s}(\xi)\vert$, for $\xi \in L^2(\cF_{r_{k}})$.  
It then makes sense to assume  
\begin{align*}
(A4)& \quad 
\normalpha{{\tt P}_{r_{k},T}^\star(\xi) - 
{\tt P}_{r_{k},T}^\star(\xi') } \leq \Lambda 
\normalpha{\xi - \xi'} \,
\hspace{120pt} { \ }
\\
(A5)& \quad 
\normalpha{{\tt P}_{r_{k},T}^\star(\xi) } \leq \Lambda 
\bigl( 1 +
\normalpha{\xi} \bigr) \,
\hspace{120pt} { \ }
\end{align*}
where $\xi,\xi' \in L^{2\alpha}(\cF_{r_{k}})$ and $
k \in \{0,\cdots,N-1\}$.

\begin{remark}
 The main challenging assumption (and maybe the most surprising one) is 
$(A3)$. It is obviously satisfied when $\alpha=1$ as long as $\Lambda$ is assumed to be greater than 1. We refer to \cite{CardaliaguetDelarueLasryLions} and 
\cite[Chap. 12]{CarmonaDelarue_book_I}
for sets of conditions under which this is indeed true. 
When $\alpha >1$, 
Assumption 
$(A3)$
is checked provided we have 
an \textit{a priori bound}
on
$\| \cU(r_{k+1},\Xd_{r_{k+1}},\law{\Xd_{r_{k+1}}})-\Yd_{r_{k+1}}\|_{2\alpha}$, {see Lemma \ref{le:A3}}. This permits to invoke the result
proven in our previous paper 
\cite{chacri15}, which holds true in a weaker setting than the 
solvability results obtained in  
\cite{CardaliaguetDelarueLasryLions} and 
\cite[Chap. 12]{CarmonaDelarue_book_I}. 
\end{remark}

%\textcolor{red}{$\mathcal{D}$ must not depend on $\mathfrak{X}$.}

\begin{theorem}
\label{th propagation error II}
We can find two constants $C_\Lambda,c_{\Lambda}>0$
and a continuous non-decreasing function $\mathfrak{B} : \RR_{+} \rightarrow \RR_{+}$ matching $0$ in $0$, only depending on $\Lambda$,  
 such that, for $\bd:=C_{\Lambda} \delta < \min(c_{\Lambda},1)$
 and $\beta \geq \mathfrak{B}(\bd)$ satisfying
\begin{align}\label{eq main ass disc}
%(J-1)N\left(\Lambda \bd^J+ \mathcal{D}^2(|\pi|)\right){e^{\beta C_\Lambda T}} \le 1\, ,
(J-1) \left( \Lambda \bd^J+ e^{\beta \bd} \mathcal{D}^2(|\pi|)\right)\frac{e^{\beta C_\Lambda T}}{e^{\beta \bd}-1} \le 1\, ,
\end{align}
where  $J$ is the number of Picard iterations in a period, it holds, for any period $k\in \{0,\cdots, N\}$ and $\xi \in L^2(\cF_{r_k})$, 
  {\rm
 \begin{align*} %\label{eq control generic disc}
  \normalpha{\sol{k}{\xi} - \cU(r_k,\xi,\law{\xi})} &\le C 
  \left(\bd^{J-1}+ (N-k)\cD^2(|\pi|)\right)\bigl( 1+ \normalphas{\xi}^\alpha
  \bigr) + C(N-k)\cD^1(|\pi|) \; ,
 \end{align*}
 }
 for a constant $C$ independent of the discretization parameters. 
% then
%  \begin{align*}
%   \max_{k} \bar{\cE}(k) \le C\left(  \delta^{J-1}\max_k \bar{\cI}(k) +N \mathcal{D}(|\pi|) \right)
%  \end{align*}
\end{theorem}

\proof The proof will follow closely the proof of Theorem \ref{th error picard} but
we now have to take into account  the discretization error.
We will first show that for all $k = \{0,\cdots,N\}$,
\begin{align} \label{eq proof hyp rec disc}
 \normalpha{\epsilon^k(\xi)}  \le {\theta}_{k}\bigl( 1+ \normalpha{{\tt P}_{r_{k},T}^\star(\xi)}^\alpha \bigr) + \vartheta_k \cD^1(|\pi|) \; ,
\end{align}
where
\begin{equation*}
\epsilon^k(\xi) =  \sol{k}{\xi} - \cU(r_k,\xi,\law{\xi}) \, ,
\end{equation*}
and $({\theta}_k,\vartheta_k)_{k = 0, \cdots, N}$ is defined by the following backward induction: $({\theta}_N,\vartheta_N) := (0,0)$, recall \eqref{eq error solver very last level solver}, and for $k \in \{0,\cdots,N-1\}$,
\begin{align} \label{eq de theta_k and vartheta_k}
\theta_k :=  \Lambda \bd^J + e^{{\beta} \bd} \set{\theta_{k+1} + \mathcal{D}^2(|\pi|)}
\;\text{ and }\;
\vartheta_k := e^{{\beta}\bd}(\vartheta_{k+1}+1) \; ,
% \theta_k := \Lambda \bd^J + \gamma \theta_{k+1} + 2\gamma \theta_{k+1} %\bd e^{2 \gamma \theta_{k+1}(J-1) \bd}\left( \gamma(J-1) \theta_{k+1} + %\frac{\Lambda}{1-\bd}\right)
\end{align}
$\beta$ being defined as in equation \eqref{eq de beta:disc}. 
%Note that $(\theta_k)$
%is defined as in Theorem \ref{th error picard}.
\vspace{5pt}
\\
Assume for a while that thids holds true. Then, we have, for all $k = \{0,\cdots,N-1\}$, %\textcolor{red}{unifier ecriture}
\begin{align} \label{eq bound theta disc}
\theta_k \le 
\bigl(\Lambda \bd^J  + e^{\beta \bd} \mathcal{D}^2(| \pi |) \bigr)
\frac{e^{\beta \bd (N-k)}-1}{e^{\beta \bd}-1}
%(N-k)\left(\Lambda \bd^J+ \mathcal{D}^2(|\pi|)\right){e^{\beta (N-k) \bd}}
\;\text{ and }\;
\vartheta_k \le e^{\beta  \bd} \frac{e^{\beta (N-k) \bd}-1}{e^{\beta \bd}-1} \; .
%\vartheta_k \le (N-k)e^{\beta (N-k) \bd}.
\end{align}
Recalling that $\bd N = C_{\Lambda} T$, we get the announced inequality. 
\vskip 5pt

We now prove 
\eqref{eq proof hyp rec disc}. Obviously, it holds true for the last step $N$. Assume now that it is true at step $k+1$, for $k<N$ and that \eqref{eq bound theta disc} holds for $\theta_{k+1}$ and $\vartheta_{k+1}$.
\\
In particular, using \eqref{eq main ass disc}, we observe that
\begin{align} \label{eq control theta k+1 disc}
 \theta_{k+1}  j \le 1, \text{ for all } j \le J-1\;.
\end{align}
\textit{First Step.}
 For $j \in \{0,\dots,J\}$, 
 let 
 $$\bar{\Delta}^j_k :=\sup_{t \in \pi^k} \normalpha{ \cU(t,\Xd^{k,j}_t,\law{\Xd^{k,j}_t})-\Yd^{k,j}_t} \, .$$
Under $(A1)-(A2)$, we will prove in this first step
an upper bound for $\bar{\Delta}^j_k$, for $j = 1,\cdots,J$, similar to the one obtained
in Proposition \ref{pr loc error picard}.

Using $(A1)$ and \HYP{1} and the fact that 
$$\Yd^{k,j}_{r_{k+1}} = \cU\bigl(r_{k+1}, \Xd^{k,j-1}_{r_{k+1}},\law{\Xd^{k,j-1}_{r_{k+1}}}\bigr) + \epsilon^{k+1}(\Xd^{k,j-1}_{r_{k+1}})\,,$$ 
we observe that
\begin{align}
\bar{\Delta}^j_k &\le e^{\Lambda \delta} \left[ \normalpha{ \cU\bigl(r_{k+1},\Xd^{k,j}_{r_{k+1}},\law{\Xd^{k,j}_{r_{k+1}}}\bigr)
-\cU\bigl(r_{k+1},\Xd^{k,j-1}_{r_{k+1}},\law{\Xd^{k,j-1}_{r_{k+1}}}\bigr) }  \right.
 \nonumber
\\
&\hspace{15pt} \left.  + \normalpha{ 
\epsilon^{k+1}(\Xd^{k,j-1}_{r_{k+1}}) } \right]
+ \Lambda \max_{j_k \le i < j_{k+1}}\normalpha{\Xd^{k,j}_{t_i} - \Xd^{k,j-1}_{t_i} } 
 + \mathcal{D}^1(|\pi|) + 
\mathcal{D}^2(|\pi|)\bigl(1+  
\normalpha{\xi}^\alpha  \bigr)
\nonumber
\\
&\le C_\Lambda \max_{t \in \pi^k}\normalpha{
\Xd^{k,j}_{t} - \Xd^{k,j-1}_{t} } + 
           e^{\Lambda \delta}
           \normalpha{ 
\epsilon^{k+1}(\Xd^{k,j-1}_{r_{k+1}}) } 
    \label{eq temp analysis II}
\\
&\hspace{15pt}
 + \mathcal{D}^1(|\pi|) +\mathcal{D}^2(|\pi|)\bigl(1+ 
 \normalpha{{\tt P}_{r_{k},T}^\star(\xi)}^\alpha \bigr)\, .  
       \nonumber
\end{align}
Using $(A2)$, we also have
\begin{align}
\sup_{t \in \pi^k} \normalpha{\Xd^{k,j}_{t}-\Xd^{k,j-1}_{t}
}
& \le \Lambda \delta \sup_{t \in \pi^k}
\left[ 
\normalpha{\Yd^{k,j}_t- \cU(t,\Xd^{k,j}_t,\law{\Xd^{k,j}_t})}
+
\Lambda \normalpha{\Xd^{k,j}_t - \Xd^{k,j-1}_t} \right.
\nonumber
\\
  &     \left.         \quad\quad\quad\quad\quad\quad\quad  
+ \normalpha{ \cU(t,\Xd^{k,j-1}_t,\law{\Xd^{k,j-1}_t}) - \Yd^{k,j-1}_t}
\right] 
\nonumber
\\
& \le C_\Lambda \delta \left( \bar{\Delta}^j_k + \bar{\Delta}^{j-1}_k \right) \, , \nonumber
\end{align}
for $\delta$ small enough.
Inserting the previous inequality in \eqref{eq temp analysis II},
we get
\begin{align*}
\bar{\Delta}^j_k 
&\le C_\Lambda \delta \bar{\Delta}_k^{j-1} 
+ e^{C_\Lambda  \delta}\normalpha{\epsilon^{k+1}(\Xd^{k,j-1}_{r_{k+1}})} + \mathcal{D}^1(|\pi|) + \mathcal{D}^2(|\pi|)
\bigl( 1+\normalpha{{\tt P}_{r_{k},T}^\star(\xi)}^\alpha \bigr)\,,
%\\
%& \le C\delta \Delta^{j-1} + e^{C_\Lambda \delta}\bar{\cE}(k) + C\mathcal{D}(|\pi|) \,,
\\
& \le \bd^j \bar{\Delta}^0_k + e^{\bd}\sum_{\ell=0}^{j-1} \bd^\ell \normalpha{\epsilon^{k+1}(\Xd^{k,j-1-\ell}_{r_{k+1}})} + \frac{\mathcal{D}^1(|\pi|)}{1-\bd} 
+ \frac{\mathcal{D}^2(|\pi|)}{1 - \bd}\bigl( 1+
\normalpha{{\tt P}_{r_{k},T}^\star(\xi)}^\alpha  \bigr) \, ,
\end{align*}
%leading to
%\begin{align*}
%\max_k \bar{\cE}(k) &\le C\left(\delta^{J-1} \max_k \bar{\cI}(k) + N\mathcal{D}(|\pi|) \right) \,.
%\end{align*}
%($\eta^Y \sim \delta |\pi|$)
with $\bd:=C_{\Lambda} \delta$. 
We note that compared to \eqref{eq pr local error}, there is a new term, namely 
 $({\mathcal{D}^1(|\pi|)} + {\mathcal{D}^2(|\pi|)} (1+ 
 \normalphas{{\tt P}_{r_{k},T}^\star(\xi)}^\alpha)/(
1 - \bd)$, which is due to the discretization. 
\vspace{5pt}
\\
\textit{Second Step.} Using \eqref{eq proof hyp rec disc} at the previous step $k+1$ and noting that $\bar{\Delta}^0_k \le \Lambda(1 + 
\normalphas{{\tt P}_{r_{k},T}^\star(\xi)}) \leq 2 \Lambda ( 1+ \normalphas{{\tt P}_{r_{k},T}^\star(\xi)}^{\alpha})$, we claim that
\begin{align}
\label{eq temp init disc}
\bar{\Delta}^j_k \le & \left(2 \Lambda  \bd^j 
+\gamma \mathcal{D}^2(|\pi|)\right)\bigl(1 +  \normalpha{{{\tt P}_{r_{k},T}^\star(\xi)}}^{\alpha}\bigr) 
+ \gamma(\vartheta_{k+1}+1)\mathcal{D}^1(|\pi|) \nonumber
\\
&+  e^{\bd} \theta_{k+1} \sum_{\ell=0}^{j-1} \bd^{j-1-\ell}
\left( 1+
\normalpha{{\tt P}_{r_{k+1},T}^\star(\bar{X}^{k,\ell}_{r_{k+1}})}^{\alpha}\right) \;,
\end{align} 
where $\gamma := e^{\bd} / (1-\bd)$.

This corresponds to equation \eqref{eq temp init} adapted to our context.
%\textcolor{magenta}{
%By $(A2)$ and $(A4)$, we have, for $\ell \leq J-1$,
%\begin{equation*}
%\begin{split}
%\normalpha{{\tt P}_{r_{k+1},T}^\star(\Xd^{k,\ell}_{r_{k+1}})
%-
%{\tt P}_{r_{k+1},T}^\star(
%\Xd^{k,0}_{r_{k+1}})}
%&\leq  
%C_{\Lambda}
%\normalpha{\Xd^{k,\ell}_{r_{k+1}}
%-
%\Xd^{k,0}_{r_{k+1}}}
%\\
%&\leq C_{\Lambda} \delta \bar{\Delta}_{k}^{\ell} + C_{\Lambda} 
%\delta
%\bigl( 1+
% \normalphas{\xi}\bigr)
%\\
%&\leq   C_{\Lambda} \delta \bar{\Delta}_{k}^{\ell} + C_{\Lambda} 
%\delta
%\bigl( 1+
% \normalpha{{\tt P}_{r_{k},T}^\star(
%\xi)}\bigr)
% \ ,
%\end{split}
%\end{equation*}
%}
By $(A2)$, we have, for $\ell \leq J-1$,
\begin{align} \label{eq temp 1}
%\begin{split}
\normalpha{{\tt P}_{r_{k+1},T}^\star(\Xd^{k,\ell}_{r_{k+1}})
-
{\tt P}_{r_{k+1},T}^\star(
\Xd^{k,0}_{r_{k+1}})}
\leq  
C_{\Lambda}
\sup_{t \in \pi^k}\normalpha{\Xd^{k,\ell}_{t}
-
\Xd^{k,0}_{t}} \;.
%\\
%&\leq C_{\Lambda} \delta \bar{\Delta}_{k}^{\ell} + C_{\Lambda} 
%\delta
% \bigl( 1+
%  \normalphas{\xi}\bigr)
% \\
% &\leq   C_{\Lambda} \delta \bar{\Delta}_{k}^{\ell} + C_{\Lambda} 
% \delta
% \bigl( 1+
%  \normalpha{{\tt P}_{r_{k},T}^\star(
% \xi)}\bigr)
%  \ ,
% \end{split}
\end{align}
Using $(A4)$, we then compute, recalling that $\Yd^{k,0}=0$,
\begin{align*}
 \sup_{t \in \pi^k}\normalpha{\Xd^{k,\ell}_{t}
-
\Xd^{k,0}_{t}}
&\le 
\Lambda \delta \sup_{t \in \pi^k} \normalpha{\Yd^{k,\ell}_t}
\\
&\le 
\Lambda \delta \left(
\bar{\Delta}_{k}^{\ell} +
\Lambda \sup_{t \in \pi^k} \normalpha{\Xd^{k,\ell}_{t}
-
\Xd^{k,0}_{t}} + \Lambda \bigl(1 + \normalpha{\xi} \bigr)
\right)
\\
& \le
C_{\Lambda} \delta \bar{\Delta}_{k}^{\ell} + C_{\Lambda} 
\delta
\bigl( 1+
 \normalpha{\xi}\bigr)\,,
\end{align*}
where for the last inequality we used the fact that $\delta$ is small enough. Observing that $\normalpha{\xi} \le \normalpha{{\tt P}_{r_{k},T}^\star(\xi)}$ and combining the previous inequality with \eqref{eq temp 1}, we obtain
\begin{align*}
 \normalpha{{\tt P}_{r_{k+1},T}^\star(\Xd^{k,\ell}_{r_{k+1}})
-
{\tt P}_{r_{k+1},T}^\star(
\Xd^{k,0}_{r_{k+1}})}
\leq
C_{\Lambda} \delta \bar{\Delta}_{k}^{\ell} + C_{\Lambda} 
\delta
\bigl( 1+
 \normalpha{{\tt P}_{r_{k},T}^\star(\xi)}\bigr)\;.
\end{align*}
So that, 
by using the fact 
that 
${\tt P}_{r_{k+1},T}^\star(
\Xd^{k,0}_{r_{k+1}}) \leq 
{\tt P}_{r_{k},T}^\star(
\xi)$ together with a
convexity argument, 
\begin{equation*}
\begin{split}
\normalpha{
{\tt P}_{r_{k+1},T}^\star(
\Xd^{k,\ell}_{r_{k+1}})}^{\alpha}
&\leq  \Bigl( 
C_{\Lambda} \delta \bar{\Delta}_{k}^{\ell} + 
\bigl( 1 +
C_{\Lambda} 
\delta
\bigr)
\bigl( 1+
 \normalpha{{\tt P}_{r_{k},T}^\star(
\xi)}\bigr)
\Bigr)^{\alpha},
\\
&\leq   \bigl( 1 + 2 C_{\Lambda} \delta  \bigr)^{\alpha-1} 
\Bigl(
  C_{\Lambda} \delta \bigl(  {\bar \Delta}_{k}^{\ell} \bigr)^{\alpha} +
 \bigl(1+  C_{\Lambda} \delta \bigr)
  \normalpha{{\tt P}_{r_{k},T}^\star(
\xi)}^{\alpha}
\Bigr) \ ,
\end{split}
\end{equation*}
Appealing to 
$(A3)$ and redefining  $\bd$, we get 
\begin{equation*}
\begin{split}
\normalpha{
{\tt P}_{r_{k+1},T}^\star(
\Xd^{k,\ell}_{r_{k+1}})}^{\alpha}
&\leq  \bd  \bar{\Delta}_{k}^{\ell} + e^{\bd}
 \bigl( 1+ 
 \normalpha{{\tt P}_{r_{k},T}^\star(
\xi)}^{\alpha}
\bigr) \ ,
\end{split}
\end{equation*}
which may be rewritten as
\begin{align*}
 \sum_{\ell=0}^{j-1} \bd^{j-1-\ell} 
 \normalpha{
{\tt P}_{r_{k+1},T}^\star(
\Xd^{k,\ell}_{r_{k+1}})}^{\alpha} 
  \le \bd 
 \sum_{\ell=1}^{j-1} \bd^{j-1-\ell}
 \bar{\Delta}_{k}^{\ell}
 +  \frac{e^{\bd}}{1-\bd} \bigl( 1+ 
  \normalpha{{\tt P}_{r_{k},T}^\star(
\xi)}^{\alpha}
\bigr)   \, .
\end{align*} 
%Using a straightforward adaptation of \eqref{eq analysis I dX} in the proof of Proposition \ref{pr loc error picard} together with $(A3)$, we obtain
%\begin{align*}
% \sum_{\ell=0}^{j-1} \bd^{j-1-\ell} \bigl\|\tilde{X}^{k,\ell}_{r_{k+1}}
% \bigr\|_{2\alpha}^{\alpha} 
% \le 
% J^{\alpha-1}
% \bd  \sum_{i=1}^{j-1}  (\bar \Delta^i_k + \bar \Delta^{i-1}_k )\sum_{\ell=i}^{j-1} \bd^{j-1-\ell} 
% + \sum_{\ell=0}^{j-1}\bd^\ell \normalpha{\xi}^{\alpha} \, . 
%\end{align*}
%Then,
%\begin{align*}
% \sum_{\ell=0}^{j-1} \bd^{j-\ell} 
% \left( 1+\normalpha{\tilde{X}^{k,\ell}_{r_{k+1}}}^{\alpha} 
%\right)
% \le 
% \frac{2 J^{\alpha-1} \bd }{1-\bd} \bar \cS^{j-1}_k + \frac{1 }{1-\bd} \normalphas{\xi}^{\alpha} \, ,
%\end{align*}
Recalling the notation 
$\gamma= e^{\bd}/(1-\bd)$ and letting $\bar \cS^{n}_k := \sum_{i=0}^n \bd^{n - i}\bar \Delta^i_k$,
we obtain a new version of \eqref{eq temp cntd}, namely
\begin{align}\label{eq disc init}
\bar{\Delta}^j_k \le  \Lambda \bd^j
\bigl(\tfrac12 + \normalpha{{\tt P}_{r_{k},T}^\star(
\xi)}^{\alpha}
\bigr) 
 + 
 \bar{a}
+  {\theta}_{k+1} \gamma \bd  \bar{\cS}^{j-1}_k 
 \;,
\end{align} 
where we changed
the constant 
$2 \Lambda$ in 
\eqref{eq temp init disc}
into $\tfrac12\Lambda$ as we changed the value of $\bd$, 
and where we put
\begin{align*}
\bar{a} =& \bigl(  \gamma^2 \theta_{k+1}+\gamma \mathcal{D}^2(|\pi|)\bigr)\bigl(1 +  \normalpha{{\tt P}_{r_{k},T}^\star(
\xi)}^{\alpha}
\bigr) 
%\\
%&
+ \gamma(\vartheta_{k+1}+1)\mathcal{D}^1(|\pi|) \, .
\end{align*}
We straightforwardly deduce that 
\begin{equation*}
\begin{split}
\bar{\cS}^{j}_k = \bar{\Delta}^j_{k} +
\bd 
\bar{\cS}^{j-1}_k &\leq  \Lambda \bd^j
\bigl(1 +
 \normalpha{{\tt P}_{r_{k},T}^\star(
\xi)}^{\alpha}
\bigr)  +  \bar{a} + 
\bigl(1+  \gamma   {\theta}_{k+1} \bigr) \bd  \bar{\cS}^{j-1}_k
\\
&\leq e^{ \gamma \theta_{k+1} j} \bd^j 
 \bar{\cS}^{0}_k + 
 \sum_{\ell=0}^{j-1}
 e^{\gamma \theta_{k+1} \ell} 
 \bd^\ell
\Bigl(  \Lambda  \bd^{j - \ell}
\bigl(1 + 
 \normalpha{{\tt P}_{r_{k},T}^\star(
\xi)}^{\alpha}
\bigr)
+
  \bar a  \Bigr) \, ,
\end{split}
\end{equation*}
which yields
\begin{equation*}
\begin{split} 
\bar{\cS}^{j}_k
 &\leq \Lambda  (j +2) \bd^j e^{\gamma \theta_{k+1} (j-1)}
 \bigl( 1 +  \normalpha{{\tt P}_{r_{k},T}^\star(
\xi)}^{\alpha}
 \bigr) + \frac{\bar{a}}{1- e^{\gamma \theta_{k+1}} \bd} \, ,
\end{split}
\end{equation*}
where we used $
\bar{\cS}^{0}_k \leq 2
\Lambda (1+  \normalphas{{\tt P}_{r_{k},T}^\star(
\xi)}^{\alpha})
$.
Thanks to 
\eqref{eq disc init}, we get
\begin{equation*}
\begin{split}
\bar{\Delta}^J_k \le  \Lambda \bd^J
\Bigl(
\tfrac12+ \bd \gamma (J+ 2) \theta_{k+1}
e^{\gamma \theta_{k+1} (J-1)}
\Bigr)
\bigl(1 + \normalpha{{\tt P}_{r_{k},T}^\star(
\xi)}^{\alpha}\bigr) 
 + 
\frac{ \bar{a}}{ 1- e^{\gamma \theta_{k+1}}\bd}
 \;.
\end{split}
\end{equation*}
Recalling that $(J-1) \theta_{k+1} \leq 1$, we deduce that, for $\bd$ small enough,  
\begin{align*}
 \bar{\Delta}^J_k &\le \bigl( \Lambda  \bd^J + e^{\beta \bd}\set{ \theta_{k+1} + \mathcal{D}^2(|\pi|)}\bigr)\bigl(1 + \normalpha{{\tt P}_{r_{k},T}^\star(
\xi)}^\alpha\bigr) + e^{\beta \bd}(\vartheta_{k+1} + 1) \cD^1(|\pi|) \;,
\end{align*}
provided that $\beta$ satisfies
\begin{equation}
\label{eq de beta:disc}
\frac{\gamma^2}{1- e^{\gamma \theta_{k+1}}\bd}
\leq e^{\beta \bd} \, .
\end{equation}
This validates \eqref{eq de theta_k and vartheta_k} and concludes the proof.
\eproof

%\input{errimplemented}

%!TEX root = mainmkv.tex
\subsection{Convergence error for the implemented scheme}

% We denote $\Pi=\cup_k \pi^k$, which the discrete time grid use for the approximation of the Brownian Motion. This approximation
% is made by using a binomial tree.\\
% On each interval $J_k=[r_{k+1}=t_{j_k},r_{k}=t_{j_{k+1}}]$, given an initial condition $\bar{\xi}$, $\widehat{\cF}_{r_{k+1}}$-measurable random variable, and $\bar{\chi}$ $\widehat{\cF}_{r_{k}}$-measurable  random variable, we thus implement the following scheme 
% \begin{align*}
% \tX{}_{\ti{+1}} & = \tX{}_{\ti{}} +b(\tY{}_{\ti{}})h+ \Delta \wh{W}_i \,,\; t_{j_k} \le \ti{} < t_{j_{k+1}}\,, \text{ and } \tX{}_{t_{j_{k}}} = \bar{\xi}\,,
% \\
% \tY{}_{\ti{}} & =  \hEc{\ti{}}{\bar{\eta}} \,,\; t_{j_k} \le \ti{} \le t_{j_{k+1}}\,.
% \end{align*}

We now analyse the global error of our method when the numerical algorithm is given by our benchmark Example \ref{ex binomial tree}, see Section \ref{subse approx U(t,x,mu)}.

\begin{lemma} \label{le stab}  (Scheme stability) Condition $(A2)$ holds true for the scheme given in Example  \ref{ex binomial tree}.
\end{lemma}

\proof  For $k\le N-1$, we consider $(\Xd,\Yd,\Zd):=\tsol{k}{\bar{\mathfrak{X}}}{\eta}{f}$ and $(\Xd',\Yd',\Zd'):=\tsol{k}{\bar{\mathfrak{X}}'}{\eta'}{f'}$ with
$\bar{\mathfrak{X}}_{r_k}=\bar{\mathfrak{X}}'_{r_k}=\xi$.
Letting $\Delta X_{i} = \Xd_{\ti{}}-\Xd'_{\ti{}}$ and $\Delta Y_{i} = \Yd_{\ti{}}-\Yd'_{\ti{}}$, we observe
\begin{align*}
|\Delta X_{i+1}| \le  \biggl\vert \sum_{\ell=j_k}^i (t_{\ell+1}-t_\ell)\Delta b_\ell \biggr\vert +
\biggl\vert 
\sum_{\ell=j_k}^i \Delta \sigma_\ell \Delta \bar{W}_\ell \biggr\vert \, ,
\end{align*}
for $i \in \{j_{k},\cdots,j_{k+1}\}$, where $\Delta b_\ell := b(\Xd_{t_\ell},\Yd_{t_\ell},\law{\Xd_{t_\ell},\Yd_{t_\ell}})-b(\Xd'_{t_\ell},\Yd'_{t_\ell},\law{\Xd'_{t_\ell},\Yd'_{t_\ell}})$ and, similarly, $\Delta \sigma_\ell := \sigma(\Xd_{t_\ell},\law{\Xd_{t_\ell}})-\sigma(\Xd'_{t_\ell},\law{\Xd'_{t_\ell}})$.
\vspace{5pt}
\\
Invoking Cauchy-Schwartz inequality for the first term and the
B\"urkholder-Davis-Gundy inequality for discrete martingales for the second term 
and appealing to the Lipschitz property of $b$ and $\sigma$, we get
\begin{align*}
\normalpha{\Delta X_{i+1}} 
&\le C \delta \max_{\ell=j_{k},\cdots,i} 
\left(
 \normalpha{\Delta Y_\ell} + \normalpha{\Delta X_\ell}\right) + C 
\biggl\| \sum_{\ell = j_k}^i 
\vert 
\Delta \sigma_\ell \vert^2 \cdot \vert \Delta \hat{W}_\ell \vert^2
\biggr\|_{\alpha}^{\frac12}
\, 
\\
& 
\leq C \delta \max_{\ell=j_{k},\cdots,i} 
\left(
 \normalpha{\Delta Y_\ell} + \normalpha{\Delta X_\ell}\right) + C 
 \biggl( \sum_{\ell = j_k}^i 
(t_{\ell+1}-t_\ell)
 \normalpha{\Delta X_\ell}^2
\biggr)^{\frac12}
\, 
\\
& 
\leq 
C \delta \max_{\ell=j_{k},\cdots,i} 
\left(
 \normalpha{\Delta Y_\ell}+ \normalpha{\Delta X_\ell} \right) + 
C \delta^{1/2} \max_{\ell=j_{k},\cdots,i}  
\left( \normalpha{\Delta X_\ell}\right) \, ,
\end{align*}
where we used the identity $t_{\ell+1}-t_\ell =\delta / (j_{k+1}-j_k)$. For $\delta$ small enough (taking the sup in the sum), we then obtain
\begin{align}
\max_{ j_k \le i \le j_{k+1}} \normalpha{\Delta X_{i}} \le C \delta \max_{ j_k \le i \le j_{k+1}} \normalpha{\Delta Y_{i}} \;,
\end{align}
which concludes the proof.
\eproof

\vspace{15pt}
We now turn to the study of the approximation error.

\begin{lemma} \label{le approx}
Assume that \HYP{0}-\HYP{1} are in force. Then, condition $(A1)$ holds true for the scheme given in Example  \ref{ex binomial tree} with
$$\cD^1(|\pi|) \le C \sqrt{|\pi|} \text{ and }  \cD^2(|\pi|) \le C \sqrt{|\pi|}  .$$
%\textcolor{red}{with extra condition on $\alpha$ and $\sigma$}
\end{lemma}

%\begin{bluetext}
\proof \textit{First Step.} Given the scheme defined in Example \ref{ex binomial tree},
we introduce its piecewise continuous version, which we denote by $(\Xd_s)_{0 \le s \le T}$. For $i<n$, $t_i < s < t_{i+1}$, 
\begin{align*}
 \Xd_{s} := \Xd_{t_{i}} + b_i (s - t_i) 
                  +  \sigma_i \sqrt{s-t_{i}} \varpi_{i},
                  \quad \varpi_{i}:=\frac{1}{\sqrt{t_{i+1}-t_{i}}} \Delta \bar{W}_i \; ,
 \end{align*}
with $(b_i,\sigma_i) := (b(\Xd_{t_{i}},\Yd_{t_i},\law{\Xd_{t_i},\Yd_{t_i}}),\sigma(\Xd_{t_i},\law{\Xd_{t_i}}))$. In preparation for the proof, we also introduce a piecewise c\`ad-l\`ag version, denoted by  $(\Xd_s^{(\lambda)})_{0 \le s \le T}$, where $\lambda$ is a parameter in $[0,1)$.   
 For $i<n$, $t_i < s < t_{i+1}$, 
 \begin{align*}
 \Xd_{s}^{(\lambda)} := \Xd_{t_{i}} + b_i (s - t_i) + 
\lambda \sigma_i \sqrt{s-t_{i}} \varpi_{i}
  \; .
 \end{align*}
For the reader's convenience, we also set 
\begin{align*}
\Ud_s &:= \cU\bigl(s,\Xd_s,\law{\Xd_s}\bigr) \,,\;
\\
\Vd_s^x &:= \partial_x \cU\bigl(s, \Xd_s,\law{\Xd_s}\bigr) \,,\; 
\Vd^{\mu}_s := \partial_\mu \cU\bigl(s, \Xd_s,\law{\Xd_s}\bigr)(\cc{\Xd_s}) \; ,
\\
\Vd^{x,0}_s & := \partial_x \cU\bigl(s, \Xd^{(0)}_s,\law{\Xd_s}\bigr) \;.
\end{align*}
Applying the discrete It\^o formula given in Proposition \ref{pr disc ito formula}, and using the PDE solved by $\cU$, recall \eqref{eq pde decoupling field}, we compute
\begin{align*}
\Ud_{t_{i+1}} &= \Ud_{t_i} 
+ \int_{t_i}^{t_{i+1}}\Vd_s^x
\cdot \left\{
b\bigl(\Xd_{t_{i}},\Yd_{t_i},\law{\Xd_{t_i},\Yd_{t_i}}\bigr) - b\bigl(\Xd_{t_{i}},\Ud_{t_i}, \law{\Xd_{t_i},\Ud_{t_i}}\bigr) 
\right\} \ud s 
\\
&\hspace{15pt} + \int_{t_i}^{t_{i+1}}
\ccesp{
\Vd_s^\mu \cdot \bigl\{ \cc{b\bigl(\Xd_{t_i},\Yd_{t_i},\law{\Xd_{t_i},\Yd_{t_i}}\bigr) - b\bigl(\Xd_{t_i},\Ud_{t_i}, \law{\Xd_{t_i},\Ud_{t_i}}\bigr) \bigr\}}
}
\ud s
\\
&\hspace{15pt} -(t_{i+1} - t_i) f\Bigl(\bar{\mathfrak X}_{t_{i}},\Ud_{t_{i}},\sigma^\dagger \bigl(\Xd_{t_i},\law{\Xd_{t_{i}}}\bigr) \Vd^x_{t_i},\law{\bar{\mathfrak{X}}_{t_i},\Ud_{t_i}}\Bigr)
\\
&\hspace{15pt}+ \Vd_{t_i}^x \cdot \Bigl( \sqrt{t_{i+1}-t_i} \sigma\bigl(\Xd_{t_i},\law{\Xd_{t_i}}\bigr) \varpi_i
\Bigr)
\\
&\hspace{15pt} + \cR^{w}_i + \cR^f_i  +\cR^{bx}_i+\cR^{b\mu}_i +\cR^{\sigma x}_i
+\cR^{\sigma \mu}_i
+ \delta {\mathcal M}(t_{i},t_{i+1})+ 
\delta {\mathcal T}(t_{i},t_{i+1}) \, ,
\end{align*}
with
\begin{align*}
 \cR^w_i := &\int_{t_i}^{t_{i+1}} ( \Vd^{x,0}_s - \Vd_{t_i}^{x,0} ) \cdot \frac{\sigma(\Xd_{t_i}^{(0)},\law{\Xd_{t_i}}) \varpi_i}{2 \sqrt{s - t_i}} \ud s\,, \;
 \\
 \cR^f_i := &\int_{t_i}^{t_{i+1}} \left \{ f\Bigl(\Xd_{s},\Ud_{s},
  \sigma^\dagger\bigl(\Xd_s,\law{\Xd_{s}}\bigr) \Vd_s^x,
  \law{\Xd_s,\Ud_s}\Bigr) \right.
  \\
&\hspace{100pt}
\left. - f\Bigl(\bar{\mathfrak{X}}_{t_{i}},\Ud_{t_{i}},
 \sigma^\dagger\bigl(\Xd_{t_i},\law{\Xd_{t_{i}}}\bigr) \Vd_{t_i}^x ,\law{\bar{\mathfrak{X}}_{t_i},\Ud_{t_i}}\Bigr) \right \} \ud s\,,
  \\
 \cR^{bx}_i := & \int_{t_i}^{t_{i+1}}\Vd_s^x
 \cdot
 \left \{
 b\bigl(\Xd_{t_{i}},\Ud_{t_i},\law{\Xd_{t_i},\Ud_{t_i}}\bigr) - 
 b\bigl(\Xd_{s},\Ud_s,\law{\Xd_s,\Ud_s}\bigr)
 \right \}
 \ud s \, ,
 \\
  \cR^{b \mu}_i  := &\int_{t_i}^{t_{i+1}}
\ccesp{
\Vd_s^\mu\cdot \bigl\{ \cc{b\bigl(\Xd_{t_{i}},\Ud_{t_i}, 
\law{\Xd_{t_i},\Ud_{t_i}}\bigr)
-b\bigl(\Xd_{s},\Ud_{s}, \law{\Xd_{s},\Ud_{s}}\bigr) }\bigr\}
}
\ud s \, , 
\end{align*}
%and
%\begin{align*}
% R^w_i & = \int_{t_i}^{t_{i+1}} ( \Vd^0_s - \Vd_{t_i} ) \frac{\sigma \eta_i}{2 \sqrt{t_{i+1} - t_i}} \ud s
% \\
% R^f_i &= \int_{t_i}^{t_{i+1}} \left \{ f(\Vd_s \sigma)- f(\Vd_{t_i} \sigma) \right \} \ud s
%\end{align*}
and
\begin{align}
 &\cR^{\sigma x}_i  =  \frac12 \int_{t_i}^{t_{i+1}}\!\! \int_0^1 \!  
 \Delta^x (s,\lambda)\ud \lambda \ud s
 \,,\; \cR^{\sigma \mu}_i  =  \frac12 \int_{t_i}^{t_{i+1}} \!\! \int_0^1 \! 
 \Delta^\mu(s,\lambda)
 \ud \lambda \ud s \, ,  \label{eq de Rsigmai}
 \end{align}
 where
 \begin{align}
 \Delta^x(s,\lambda) & :=
 \textrm{\rm Tr}
 \Bigl\{ 
 \partial_{xx}^2 \cU\bigl(s,\Xd^{(\lambda)}_s,\law{\Xd_s}\bigr)
 a(\Xd_{t_i},\law{\Xd_{t_i}})
 - \partial_{xx}^2 \cU\bigl(s,\Xd_s,\law{\Xd_s}\bigr)a(\Xd_{s},\law{\Xd_{s}})
\Bigr\} \nonumber
\\
\Delta^\mu(s,\lambda) & :=
\!\hat{\EE} \Bigl[
\textrm{\rm Tr}
\Bigl\{ 
\partial_v \partial_\mu \cU\bigl(s,\Xd_s,\law{\Xd_s}\bigr)(\cc{\Xd^{(\lambda)}_s})\cc{a(\Xd_{t_i},\law{\Xd_{t_i}})} 
\nonumber
\\
&\hspace{100pt}-
\partial_v \partial_\mu \cU\bigl(s,\Xd_s,\law{\Xd_s}\bigr)(\cc{\Xd_s}) \cc{a(\Xd_{s},\law{\Xd_{s}})} \Bigr\}
 \Bigr] \nonumber
\;. 
 \end{align}
Also, $\delta {\mathcal M}(t_{i},t_{i+1})$ is a martingale increment satisfying $\esp{ \vert \delta {\mathcal M}(t_{i},t_{i+1}) \vert^{2\alpha} \, \vert \, \cF_{t_{i}}}^{1/(2\alpha)}
 \leq C h_{i}$ and 
 $\normalphas{\delta {\mathcal T}(t_{i},t_{i+1})} \leq C_{\Lambda} h_{i}^{\frac32}$, recall Proposition \ref{pr disc ito formula}.
%  
% \begin{align}
%  &2\cR^\sigma_i  =  \label{eq de Rsigmai}
%  \\
%  & \!\!\int_{t_i}^{t_{i+1}} \!\!\!\int_0^1\!
%  \left \{ 
%  \partial_{xx}^2 U(s,\Xd^{(\lambda)}_s,\law{\Xd^{(1)}_s})\sigma(\Xd_{t_i})^2\eta_i^2
%  - \partial_{xx}^2 U(s,\Xd^{(1)}_s,\law{\Xd^{(1)}_s})\sigma(\Xd^{(1)}_{s})^2
% \right \}
% \ud \lambda
% \ud s \nonumber
% \\
% &\quad\quad\quad\quad\quad\quad+ \nonumber
% \\
% &\!\!\int_{t_i}^{t_{i+1}}\!\!\! \int_0^1
% \!\ccesp{\partial_v \partial_\mu U(s,\Xd^{(1)}_s,\law{\Xd^{(1)}_s})(\cc{\Xd^{(\lambda)}_s})\cc{\sigma(\Xd_{t_i})^2\eta_i^2} \right. \nonumber
% \\
% &\quad\quad\quad\quad\quad\quad
% \left. -
% \partial_v \partial_\mu U(s,\Xd^{(1)}_s,\law{\Xd^{(1)}_s})(\cc{\Xd^{(1)}_s})\cc{\sigma(\Xd^{(1)}_{s})^2} \nonumber
% }
% \ud \lambda
% \ud s \;. \nonumber
% %%
%  \end{align}
\vskip 5pt

\noindent \textit{Second Step.} 
Denoting $\urpi_i := \varpi_i/\sqrt{t_{i+1} - t_i}$ and
\begin{align*}
\delta b_i &:= \frac1{h_i}\int_{t_i}^{t_{i+1}}\Vd_s^x \cdot \left\{
b\bigl(\Xd_{t_{i}},\Yd_{t_i},\law{\Xd_{t_i},\Yd_{t_i}}\bigr) - b\bigl(\Xd_{t_{i}},\Ud_{t_i}, \law{\Xd_{t_i},\Ud_{t_i}}\bigr) 
\right\} \ud s 
\\
&\hspace{15pt} + \frac1{h_i} \int_{t_i}^{t_{i+1}}
\ccesp{
\Vd_s^\mu \cdot \bigl\{ \cc{b\bigl(\Xd_{t_{i}},\Yd_{t_i},\law{\Xd_{t_i},\Yd_{t_i}}\bigr) - b\bigl(\Xd_{t_{i}},\Ud_{t_i}, \law{\Xd_{t_i},\Ud_{t_i}}\bigr)}\bigr\}
}
\ud s \, ,
\end{align*}
%Introducing  
%$\beta_i = \frac1{t_{i+1}-t_i}\int_{t_i}^{t_{i+1}} \Vd_s \ud s $,
%$\beta^\mu_i = \frac1{t_{i+1}-t_i}\int_{t_i}^{t_{i+1}} \Vd^\mu_s \ud s $, 
%$\delta b_i = b(\Yd_{t_i},\law{\Xd_{t_i},\Yd_{t_i}}) - b(\Ud_{t_i}, \law{\Xd_{t_i},\Ud_{t_i}}) $,
%$\Hd_i = \frac{\eta_i}{\sqrt{t_{i+1} - t_i}}$
the previous equation reads
\begin{equation}
\label{eq bar U disc}
\begin{split}
 \Ud_{t_{i+1}} &= \Ud_{t_i} + \zeta_i
 \\ 
&+ h_i \Bigl[ \delta b_i  
  - f\Bigl(\bar{\mathfrak{X}}_{t_{i}},\Ud_{t_{i}}, \sigma^\dagger\bigl(\Xd_{t_i},\law{\Xd_{t_i}}\bigr) \Vd_{t_i}^x,\law{\bar{\mathfrak{X}}_{t_i},\Ud_{t_i}}\Bigr)
 + \Vd_{t_i}^x \cdot \bigl( \sigma(\Xd_{t_i},\law{\Xd_{t_i}}) \urpi_i 
 \bigr) 
  \Bigr] \;,%\nonumber
\end{split}
\end{equation}
where 
\begin{equation} \label{eq de zeta}
\begin{split}
\zeta_i  
&:= \cR^{w}_i + \cR^f_i
+
\cR^{b x}_i+\cR^{b \mu}_i  + \cR^{\sigma x}_i + \cR^{\sigma \mu}_i
+
\delta {\mathcal M}(t_{i},t_{i+1})
+
\delta {\mathcal T}(t_{i},t_{i+1})   \; . \nonumber
\end{split}
\end{equation}
On the other hand, the scheme can be rewritten as
\begin{align} \label{eq scheme expli}
\Yd_{t_i} = \Yd_{t_{i+1}} 
+ h_i f\bigl(\bar{\mathfrak{X}}_{t_i},\Yd_{t_{i}},\Zd_{t_i},\law{\bar{\mathfrak{X}}_{t_i},\Yd_{t_i}}\bigr) 
- h_i \Zd_{t_i} \cdot \urpi_i  - \Delta M_i  \; ,
\end{align}
where $\Delta M_i$ satisfies 
\begin{align} \label{eq prop Delta M}
\EFp{\ti{}}{\Delta M_i} = 0\, , \; \EFp{\ti{}}{\urpi_i \cdot \Delta M_i} = 0 \text{ and } \esp{|\Delta M_i|^2} < \infty \,.
\end{align}
%\textcolor{red}{rem: actually in the case of the bin tree, $\Delta M = 0$.}
Denoting $\Delta \Yd_i = \Yd_{t_i} - \Ud_{t_i}$, $\Delta \Zd_i = \Zd_{t_i} -  \sigma^\dagger(\Xd_{t_i},\law{\Xd_{t_i}}) \Vd_{t_i}^x$,
and adding \eqref{eq bar U disc} and \eqref{eq scheme expli}, we get
\begin{align}\label{eq bckwd scheme}
 \Delta \Yd_{i} = \Delta \Yd_{i+1} + h_i\left(\delta b_i + \delta	 f_i \right) + \zeta_i  - h_i \Delta \Zd_i \cdot \urpi_i - \Delta M_i \; ,
\end{align}
where 
\begin{equation*}
\begin{split}
\delta f_i 
= f\bigl(\bar{\mathfrak{X}}_{t_i},\Yd_{t_{i}},
\Zd_{t_i},\law{\bar{\mathfrak{X}}_{t_i},\Yd_{t_i}}\bigr)- 
 f\Bigl(\bar{\mathfrak{X}}_{t_i},\Ud_{t_{i}},
\sigma^\dagger\bigl(\Xd_{t_i},\law{\Xd_{t_i}}\bigr) \Vd^x_{t_i},\law{\bar{\mathfrak{X}}_{t_i},\Ud_{t_i}}\Bigr) \;.
\end{split}
\end{equation*}
%$$\gamma^z_i = \frac{f(\Zd_{t_i},\law{\bar{\mathfrak{X}}_{t_i},\Yd_{t_i}})-f(\Vd_{t_i}\sigma(\Xd_{t_i}),\law{\bar{\mathfrak{X}}_{t_i},\Yd_{t_i}})}{\Zd_{t_i} - \Vd_{t_i}\sigma(\Xd_{t_i})}\1_{\set{\Pd_{i}  \neq 0}}$$ and
%$$\gamma^y_i = \frac{f(\Vd_{t_i}\sigma(\Xd_{t_i}),\law{\bar{\mathfrak{X}}_{t_i},\Yd_{t_i}})
%-f(\Vd_{t_i}\sigma(\Xd_{t_i}),\law{\bar{\mathfrak{X}}_{t_i},\Ud_{t_i}})}{\|\Yd_{t_i} - \Ud_{t_i} \|_2}\1_{\set{\|\Gd_{i}\|_2  \neq 0}}.$$
For later use, we observe that
\begin{align}\label{eq de bound db+df}
|\delta b_i| + |\delta f_i| \le C_{\Lambda} \bigl(|\Delta \Yd_i| + \norm{\Delta \Yd_i} + |\Delta \Zd_i| \bigr)\;.
\end{align}
Summing the  equation \eqref{eq bckwd scheme} from $i$ to $j_{k+1}-1$, we obtain
\begin{align*}
\Delta \Yd_i + \sum_{\ell = i}^{j_{k+1}-1}\set{h_{\ell} \Delta \Zd_{\ell} 
\cdot \urpi_{\ell}  +\Delta M_{\ell} }
=
\Delta \Yd_{j_{k+1}} +  \sum_{\ell = i}^{j_{k+1}-1} h_{\ell}\left(\delta b_{\ell} + \delta	 f_{\ell} \right)
- \sum_{\ell = i}^{j_{k+1}-1} \zeta_\ell\,.
\end{align*}
Squaring both sides and taking expectation, we compute, using \eqref{eq prop Delta M} for the left side and Young's and conditional Cauchy-Schwarz inequality for the right side,
\begin{equation*}
\begin{split}
&\EFp{t_q{}}{|\Delta \Yd_i|^2 } + \sum_{\ell = i}^{j_{k+1}-1}\!\!\! h_{\ell} \EFp{t_q{}}{|\Delta \Zd_{\ell}|^2 }
\\
&\hspace{15pt}
\le \EFp{t_q{}}{(1+C\delta)|\Delta \Yd_{j_{k+1}}|^2
+  C\sum_{\ell = i}^{j_{k+1}-1} h_{\ell}|\delta b_{\ell} + \delta	 f_{\ell} |^2
+ \frac{C}\delta \biggl(\sum_{\ell = i}^{j_{k+1}-1} \zeta_\ell\biggr)^2 }\, ,
\end{split}
\end{equation*}
for $i \geq q \geq j_{k}$.
Combining \eqref{eq de bound db+df} and Young's inequality, this leads to
\begin{equation*}
\begin{split}
&\EFp{t_q{}}{|\Delta \Yd_i|^2 } + \frac12\sum_{\ell = i}^{j_{k+1}-1}\!\!\! h_{\ell} \EFp{t_q{}}{|\Zd_{\ell}|^2}
\le \EFp{t_q{}}{e^{C\delta}|\Delta \Yd_{j_{k+1}}|^2
+  C\sum_{\ell = i}^{j_{k+1}-1} h_{\ell}|\Delta \Yd_{\ell}|^2
+ \frac{C}{\delta} \biggl(\sum_{\ell = i}^{j_{k+1}-1} \zeta_\ell
\biggr)^2  }\,.
\end{split}
\end{equation*}
Using the discrete version of Gronwall's lemma and recalling that 
$\sum_{\ell = j_{k}}^{j_{k+1}-1} h_{\ell} = \delta$, we obtain,
for $i=q$,
\begin{equation*}
\begin{split}
|\Delta \Yd_i|^2 
 \le
\EFp{t_i{}}{e^{C \delta}|\Delta \Yd_{j_{k+1}}|^2
+ \frac{C}{\delta} \max_{j_{k} \le i \le j_{k+1}-1}\biggl(\sum_{\ell = i}^{j_{k+1}-1} \zeta_\ell \biggr)^2  } \; ,
\end{split}
\end{equation*}
and then,
\begin{equation}
\label{eq stability}
\begin{split}
\Delta_Y^2 &:= \max_{j_k\le i \le j_{k+1}} \normalpha{\Delta \Yd_i}^2
 \le
e^{C \delta}\normalpha{\Delta \Yd_{j_{k+1}}}^2
+ \frac{C}{\delta} \biggl\| \max_{j_{k} \le i \le j_{k+1}-1}\biggl(\sum_{\ell = i}^{j_{k+1}-1} \zeta_\ell \biggr) \biggr\|_{2\alpha}^2 \;.
\end{split}
\end{equation}
\vskip 5pt

\noindent \textit{Third Step.} To conclude, we need an upper bound for the
error
$\normalphas{\max_{j_{k} \le i \le j_{k+1}-1} (\sum_{\ell = i}^{j_{k+1}-1} \zeta_\ell)}^2$
where $\zeta_\ell$ is defined in \eqref{eq de zeta}. To do so, we
study each term in \eqref{eq de zeta} separately.
We also define
%\begin{align*}
$
 \Delta_X := \max_{t \in \pi^k}\normalpha{\Xd_t - \bar{\mathfrak{X}}_t}
$ 
and we
 %\text{ and } 
 %\Delta_Y := \max_{t \in \pi^k}\norm{U(t,\Xd_t,\law{\Xd_t}) - \Yd_t}\;.
%\end{align*}
recall that $\Xd_{r_k} = \xi$.
\vskip 5pt

\noindent \textit{Third Step a.} We first study the contribution of $\cR^f_i$ to the global error term and note that
\begin{align} \label{eq decomp de base}
\biggl\|
\max_{j_{k} \le i \le j_{k+1}}
\biggl(
\sum_{\ell = i}^{j_{k+1}-1} \cR^f_{\ell}\biggr) \biggr\|^2_{2\alpha} \le 
C \frac{\delta}{\vert \pi \vert} \sum_{\ell = j_k}^{j_{k+1}-1} \normalpha{\cR^f_{\ell}}^2 \, .
 %\esp{ \frac{|\EFp{t_i}{\cR^f_i}|^2}{h_i} + |\cR^f_i|^2 } \le \frac{2}{h_i} \esp{|\cR^f_i|^2}\;.
\end{align}
We will upper bound this last term.
\vspace{5pt}
\\
Let us first observe, that, for $t_i \le s \le t_{i+1}$, 
\begin{align*}
|\Vd_s^x - \Vd_{t_i}^x|  & \le \bigl|\partial_x \cU\bigl(s, \Xd_s,\law{\Xd_s}\bigr) - \partial_x \cU\bigl(t_{i}, \Xd_{t_i},\law{\Xd_{t_{i}}}\bigr)\bigr| 
%\\
%&\hspace{15pt} +   \bigl|\partial_x U\bigl(s, \Xd_{t_i},\law{\Xd_s}\bigr) - \partial_x \cU\bigl(s, \Xd_{t_i},\law{X^{t_i,\Xd_{t_i}}_s}\bigr)\bigr| \label{eq disc temp 2}
%\\
%&\hspace{15pt}+
%\bigl|\partial_x \cU\bigl(s, \Xd_{t_i},\law{X^{t_i,\Xd_{t_i}}_s}\bigr) - \partial_x \cU\bigl(s, \Xd_{t_i},\law{\Xd_{t_i}}\bigr)\bigr| 
%\label{eq disc temp 3}
\\
&\le C\left(|\Xd_s-\Xd_{t_i}|+\cW_2(\law{\Xd_s},\law{\Xd_{t_i}})
 + h_i^\frac12
\bigl( 1+ \vert \Xd_{t_{i}} \vert 
 + 
 \| \Xd_{t_{i}} \|_{2} \bigr)
 \right) \, ,
% \label{eq disc temp V final}
\end{align*}
where
%, to upper bound the term in the RHS of \eqref{eq disc temp 1}, 
we used the Lipschitz property of $\partial_x \cU$
given in \HYP{1}, together with 
\eqref{eq cont grad mu}
and 
\eqref{eq cont grad t}.
% Note that
% \begin{align*}
%  \cW_2(\law{\Xd^{(1)}_s},\law{X^{t_i,\Xd_{t_i}}_s}) \le \left( 
%  \norm{\Xd^{(1)}_s-\Xd_{t_i}} 
%  + \norm{X^{t_i,\Xd_{t_i}}_s-\Xd_{t_i}}
%  \right)\;.
% \end{align*}
% Combining the previous inequality with \eqref{eq disc temp V final}, we obtain
 Hence,
 \begin{align}\label{eq interm Vd}
  \normalpha{\Vd_s^x - \Vd_{t_i}^x}^2 \le 
  C\left(  \normalpha{\Xd_s-\Xd_{t_i}}^2
%  C\left( h_i + \esp{|\Xd^{(1)}_s-\Xd_{t_i}|^2
%  + 
%  |X^{t_i,\Xd_{t_i}}_s-\Xd_{t_i}|^2
%  }
  +h_i
\bigl( 1+ 
 \| \Xd_{t_{i}} \|_{2 \alpha}^2 \bigr) 
  \right)\;.
 \end{align}
% which leads, using Lemma \ref{le basic estimates for X&Y}, to
% \begin{align}
%  \esp{|\Vd_s - \Vd_{t_i}|^2} \le C\left( h_i + h_i^2\norm{\eta}^2
%   + 
%   \esp{
%   |X^{t_i,\Xd_{t_i}}_s-\Xd_{t_i}|^2
%   }
%   \right) \label{eq disc temp esp(V)}
% \end{align}
%For the last term in the RHS of the previous inequality, we first observe
%that from the  Lipschitz property \eqref{eq de U lipschitz} of $U$, one obtains 
%classicaly 
%\begin{align*}
% \esp{|X^{t_i,\Xd_{t_i}}_s|^2} \le C(1+\esp{|\Xd_{t_i}|^2} )
%\end{align*}
%and then
%\begin{align}\label{eq control true process one time step}
% \esp{
%  |X^{t_i,\Xd_{t_i}}_s-\Xd_{t_i}|^2
%  } \le C(h_i + h_i^2 \esp{|\Xd_{t_i}|^2})\;.
%\end{align}
From the boundedness of $\sigma$ and the Lipschitz property of $b$ and $\cU$, we
compute
\begin{align}
\label{eq bound dX}
 \normalpha{\Xd_s-\Xd_{t_i}}^2
 \le
 C_{\Lambda}\left(h_i + h_i^2\normalpha{\bar{U}_{t_i}-Y_{t_i}}^2 +
h_{i}^2 \normalpha{\Xd_{t_i}}^2 \right) \; .
\end{align}
%Inserting the previous inequality and \eqref{eq control true process one time step} into \eqref{eq interm Vd}, we obtain, 
Using Lemma \ref{le basic estimates for X&Y} from the appendix below,
we obtain
%
%Inserting the previous inequality into \eqref{eq disc temp esp(V)} and using %Lemma \ref{le basic estimates for X&Y}, we get
\begin{align*}
  \normalpha{\Vd_s^x - \Vd_{t_i}^x}^2 \le C\Bigl( h_i 
 \bigl( 1
  + 
 \normalphas{\xi}^2\bigr)
  + h_i^2 \Delta_Y^2
  \Bigr) \, .%\label{eq disc temp esp(V) final}
\end{align*}
From the boundedness of $\partial_x \cU$, $\sigma$ and the lipschitz property
of $\sigma$, we 
obtain
\begin{align*}
  \normalpha{
 \sigma^\dagger\bigl(\Xd_s,\law{\Xd_s}\bigr)  \Vd_s^x  
 - \sigma^\dagger \bigl(\Xd_{t_i},\law{\Xd_{t_{i}}}\bigr) \Vd_{t_i}^x}^2 \le C\left( h_i 
 \bigl( 1+ 
 \normalphas{\xi}^2 \bigr) 
 + h_i^2 \Delta_Y^2 
\right) \, ,%\label{eq disc temp esp(Vsigma) final}
\end{align*}
where we used the same argument as above to handle the 
difference between the two $\sigma$ terms.
%We observe that $\norm{\eta} \le \norm{U(r_k,\Xd_{r_k},\law{\Xd_{r_k}}) - %\eta} + 2L(1+ \norm{\Xd_{r_k}})$ from the linear growth property of $U$.
%
Combining the previous inequality  with the Lipschitz property of $f$ and replicating the analysis to handle the difference between the $\bar{U}$ terms, we  deduce
\begin{align}
\normalpha{\cR^f_i}^2 \le Ch_i^2
\Bigl( \Delta_X^2
+
h_i 
\bigl( 1+ \normalphas{\xi}^2 \bigr)
+ h_i^2 \Delta_Y^2 
\Bigr)\;. \label{eq control Rfi}
\end{align}

\noindent \textit{Third Step b.} Combining the Lipschitz property of $b$, 
the fact that $|\Vd_s^x|^2 + \hat{\EE}[ \vert \Vd_s^\mu|^2] \le C$ and Cauchy-Schwarz inequality, we get
\begin{align} \label{eq interm Rbi}
\normalpha{\cR^{b x}_i}^2+ \normalpha{\cR^{b \mu}_i}^2 &\le Ch_i^2  \normalpha{\Ud_s - \Ud_{t_i}}^2 + \normalpha{\Xd_s - \Xd_{t_i}}^2 \, .
%\\
%&\le C h_i^2\esp{ |\Xd_s - \Xd_{t_i}|^2}
\end{align}
%where we used for the last inequality the Lipschitz property of $U$ given in
%\eqref{eq de U lipschitz}.
%From the boundedness of $\sigma$ and the Lipschitz property of $b$ and $U$, we
%also compute that,
%\begin{align} \label{eq bound dX}
% \esp{ |\Xd^{(1)}_s - \Xd_{t_i}|^2} \le C_L\left(h_i + h_i^2\esp{|\bar{U}_{t_i}-Y_{t_i}|^2 + |\Xd_{t_i}|^2} \right) \;.
%\end{align}
%Inserting the previous inequality into \eqref{eq interm Rbi} and using lemma
%\ref{le basic estimates for X&Y}, we obtain
Arguing as in the previous step, we easily get
\begin{align}
\normalpha{\cR^b_i}^2 &\le Ch_i^2\Bigl(h_i
\bigl( 1
+ \normalphas{\xi}^2\bigr)
+h_i^2 \Delta_Y^2 \Bigr) \;.\label{eq control Rbi}
\end{align}
%which leads to, combining Corollary \ref{co control bar U} and \textcolor{red}{basic estimates} %for $\Xd$,
%\begin{align}
%\desp{|\cR^b_i|^2} \le Ch_i^3 \;. \label{eq control Rbi}
%\end{align}
\textit{Third Step c.}
 We now study the contribution of the terms $\cR^w_i$ to the global error.
From the independance property of $(\varpi_i)_{i=0,\cdots,n-1}$,
we may regard each $\cR^w_{\ell}$ as a martingale increment. 
By Burkholder-Davies-Gundy inequalities for discrete martingales, we first compute, using the fact that each $\urpi_{i}$ is uniformly bounded,
\begin{align*}
\biggl\| \max_{j_{k} \le i \le j_{k-1}}
 \biggl(\sum_{\ell = i}^{j_{k+1}-1} \cR^w_\ell \biggr)
 \biggr\|^2_{2\alpha}
&\leq C \left\| \sum_{\ell = j_k}^{j_{k+1}-1}\left| \int_{t_i}^{t_{i+1}} 
\sigma^\dagger(\Xd_{t_i}^{(0)},\law{\Xd_{t_i}})
 \frac{\bar{V}^{0,x}_s -\bar{V}^{0,x}_{t_i}}{\sqrt{s-t_i}}\ud s\right|^2 \right\|_{\alpha}
 \\
 & \le C\sum_{\ell = j_k}^{j_{k+1}-1} h_i 
 \biggl( h_{i} \bigl( 1 + \normalphas{\xi}^2 \bigr)
+ 
 \normalphadeuxB{\sup_{s \in [t_i,t_{i+1}] } |\Xd_{s}^{(0)} - \Xd_{t_i}|^2 }\biggr) \; .
\end{align*}
Since $|\Xd_{s}^{(0)} - \Xd_{t_i}| \le h_i|b(\Xd_{t_i},\Yd_{t_i},\law{\Xd_{t_i},\Yd_{t_i}})|$, for $s \in [t_i,t_{i+1}]$,
so that 
$\| \Xd_{s}^{(0)} - \Xd_{t_i} \|_{2\alpha} \leq
C_{\Lambda} h_{i} ( 1+ \|\Xd_{t_{i}} \|_{2\alpha}+ 
\| \Yd_{t_i} \|_{2 \alpha} )
\leq C_{\Lambda} h_{i} ( 1+ 
\|\Xd_{t_{i}} \|_{2\alpha} 
+ 
\Delta^2_{Y})$, 
 the previous inequality, together with 
Lemma \ref{le basic estimates for X&Y}, 
 leads to
\begin{align*}
 \biggl\| \max_{j_{k} \le i \le j_{k-1}}
 \biggl(\sum_{\ell = i}^{j_{k+1}-1} \cR^w_\ell \biggr)
 \biggr\|^2_{2\alpha} \le C \delta \vert \pi \vert  \Bigl(1 + \normalphas{\xi}^2
 + \vert \pi \vert \Delta_Y^2 
  \Bigr)\;.
\end{align*}
Similarly, 
\begin{align*}
 \biggl\| \max_{j_{k} \le i \le j_{k-1}}
 \biggl(\sum_{\ell = i}^{j_{k+1}-1} 
 \delta {\mathcal M}(t_{\ell},t_{\ell+1})
 \biggr)
 \biggr\|^2_{2\alpha} 
 &\leq C\, \sum_{\ell = j_{k}}^{j_{k+1}-1}
\bigl\| 
\bigl\vert \delta {\mathcal M}(t_{\ell},t_{\ell+1})
\bigr\vert^2
\bigr\|_{\alpha}
 \\
& \le C \, \delta \vert \pi \vert  \;.
\end{align*}
Hence,
\begin{equation}
\label{eq control Rwi}
\begin{split}
& \biggl\| \max_{j_{k} \le i \le j_{k-1}}
 \biggl(\sum_{\ell = i}^{j_{k+1}-1} \cR^w_\ell \biggr)
 \biggr\|^2_{2\alpha}
 + \biggl\| \max_{j_{k} \le i \le j_{k-1}}
 \biggl(\sum_{\ell = i}^{j_{k+1}-1} 
 \delta {\mathcal M}(t_{\ell},t_{\ell+1})
 \biggr)
 \biggr\|^2_{2\alpha}
\\ 
&\hspace{15pt}  \le C \delta \vert \pi \vert \bigl( 1
+ \normalphas{\xi}^2 \bigr)
+ C \delta \vert \pi \vert^2  \Delta_Y^2 \;.
\end{split}
\end{equation}
% 
% \begin{align*}
% |\cR^w_i| & \le C \sup_{s \in [t_i,t_{i+1}]}|\Xd_s^{(0)} - \Xd_{t_i}| |\sigma(\Xd_{t_i}) \eta_i| \sqrt{h_i}\;,
% \\
% & \le C |b(\Yd_{t_i},\law{\Xd_{t_i},\Yd_{t_i}})|h_i^{\frac32}
% \end{align*}
% where we used the Lipschitz continuity of $\partial_x U$ and the definition of $X^{(0)}$.
% Combining the previous inequality with Lemma \ref{le basic estimates for X&Y}, we obtain
% \begin{align}
% \desp{|\cR^w_i|^2} \le C h_i^3(\norm{U(r_k,\Xd_{r_k},\law{\Xd_{r_k}}) - \eta}^2 + \norm{\xi}^2+ \delta)\;.   \label{eq control Rwi}
% \end{align}
\textit{Third Step d.} (i) We study the contribution of $\cR^{\sigma x}_i$. 
%First, using the
%boundedness of \textcolor{red}{$\partial^2_{xx}U$}, we compute
%\begin{align}\label{eq Rsigmax strong}
% \esp{|\cR^{\sigma x}_i|^2} \le C h_i^2 \;. 
%\end{align}
We observe that % can be decomposed as the sum
\begin{align*}
|\Delta^x(s,\lambda)| \le
  &\bigl|\partial^2_{xx}\cU\bigl(s,\Xd^{(\lambda)}_s,\law{\Xd_s}\bigr)
  -
  \partial^2_{xx}\cU\bigl(s,\Xd_{s},\law{\Xd_s}\bigr)
 \bigr| \cdot \vert a(\Xd_{t_i},\law{\Xd_{t_i}}) \vert  
 \\
 &+ \bigl|
 \partial^2_{xx}\cU(s,\Xd_s,\law{\Xd_s})
 \bigr\vert \cdot \bigl\vert
 a(\Xd_{t_i},\law{\Xd_{t_i}}) 
  -
 a(\Xd_{s},\law{\Xd_{s}})
\bigr| 
  \, ,
\end{align*} 
for $s \in [t_{i},t_{i+1}]$.
Using the boundedness and Lipschitz continuity of $\partial^2_{xx} \cU$ and $\sigma$, we get, from the previous expression,
\begin{align}
 \normalpha{\Delta^x(s,\lambda)}^2  \le 
 C\left( \normalpha{\Xd^{(\lambda)}_s- \Xd_{s}}^2 + \normalpha{\Xd_s- \Xd_{t_i}}^2 \right) \, .
\end{align}
Observing that $\normalpha{\Xd_s^{(\lambda)}-\Xd_s} \le C \sqrt{h_i}$, we obtain
using \eqref{eq bound dX}, for $ t_i \le s \le t_{i+1}$
\begin{align*} 
 \normalpha{\Delta^x(s,\lambda)}^2  \le Ch_i\bigl(1+h_i
 \normalpha{\bar{U}_{t_i}-Y_{t_i}}^2 + 
 h_{i}
 \normalpha{\Xd_{t_i}}^2
 \bigr)\;.
\end{align*}
which leads, using Lemma \ref{le basic estimates for X&Y} again, to
\begin{align}\label{eq control R sigma x}
\normalpha{\cR^{\sigma x}_i}^2 \le C h_i^2\Bigl(h_i+h_i^2 \bigl(\Delta_Y^2 + \normalphas{\xi}^2\bigr)\Bigr) \, .
 %\esp{ \frac{|\EFp{t_i}{\cR^{\sigma x}_i}|^2}{h_i} + |\cR^{\sigma x}_i|^2 } \le C h_i\left(h_i   +h_i^2( \norm{U(r_k,\Xd_{r_k},\law{\Xd_{r_k}}) - \eta}^2 + \norm{\xi}^2)\right)\;.
\end{align}

%
% The first term in the RHS of \eqref{eq de Rsigmai} is controlled by using the Lipshitz continuity of $\partial_{xx} U$ given in \HYP{2a}
% and \textcolor{red}{basic estimates} for $\Xd$. 

(ii) To study $\cR^{\sigma \mu}_i$, we first observe 
that 
\begin{align}
 |\Delta^\mu(s,\lambda)| 
 & \le C 
\ccesp{  
 \bigl|\partial_\upsilon \partial_\mu \cU(s,\Xd_s,\law{\Xd_s})(\cc{\Xd_s^{(\lambda)}})
  -
  \partial_\upsilon \partial_\mu \cU(s,\Xd_s,\law{\Xd_s})(\cc{\Xd_s})\bigr|  
  }
     \label{eq decomp Delta mu}
 \\
 &\hspace{15pt} + 
 \ccesp{ 
 \bigl|\partial_\upsilon \partial_\mu \cU(s,\Xd_s,\law{\Xd_s})(\cc{\Xd_s})
 \bigr\vert \cdot 
 \bigl\vert
 \cc{a(\Xd_{t_i},\law{\Xd_{t_i}}) 
  -
a(\Xd_{s},\law{\Xd_{s}}) 
  } \bigr\vert} \, .
  \nonumber
\end{align}
For the last term, we combine Cauchy-Schwarz inequality \eqref{eq bound second deriv} and boundedness and Lipschitz continuity of $\sigma$ to get
\begin{align*}
 &\ccesp{ \bigl|\partial_\upsilon \partial_\mu \cU(s,\Xd_s,\law{\Xd_s})(\cc{\Xd_s}) \bigr\vert \cdot
 \bigl\vert 
 \cc{a(\Xd_{t_i},\law{\Xd_{t_i}})
  -
  a(\Xd_{s},\law{\Xd_{s}})  } \bigr|}
  \\
 &\hspace{15pt}
 \le C 
 \norm{\Xd_{t_i}-\Xd_{s}}
 \le 
 C 
 \normalpha{\Xd_{t_i}-\Xd_{s}}.
\end{align*}
Recalling
from \eqref{eq bound dX} that
%\begin{align*}
 $
 \normalpha{\Xd_s-\Xd_{t_i}}^2
 \le
 C_{\Lambda}(h_i + h_i^2 (\Delta_{Y}^2  + \|\Xd_{t_i}\|_{2\alpha}^2) )$, %\end{align*}
we obtain, using Lemma \ref{le basic estimates for X&Y}, that
\begin{align}
 &\ccesp{ \bigl|\partial_\upsilon \partial_\mu \cU(s,\Xd_s,\law{\Xd_s})(\cc{\Xd_s}) \bigr\vert
 \cdot
 \bigl\vert 
 \cc{a(\Xd_{t_i},\law{\Xd_{t_i}})  -
  a(\Xd_{s},\law{\Xd_{s}})
  } \bigr|
} \nonumber
  \\
 &\hspace{15pt}
 \le 
 C_{\Lambda} h_i^\frac12\left(1 + h_i^\frac12\set{ \Delta_Y + \normalphas{\xi}} \right) \, .
 \label{eq control Delta mu 1}
 \end{align}
For the first term in \eqref{eq decomp Delta mu},
we  use $\HYP{1}$ equation \eqref{eq control dvdmu} to get
\begin{align*}
 &|\partial_\upsilon \partial_\mu \cU(s,\Xd_s,\law{\Xd_s})(\cc{\Xd_s^{(\lambda)}})
  -
  \partial_\upsilon \partial_\mu \cU(s,\Xd_s,\law{\Xd_s})(\cc{\Xd_s})|
  \\
  & \le
 C\left\{ 1 + |\cc{\Xd_s^{(\lambda)}}|^{2\alpha} 
	      + |\cc{\Xd_s}|^{2 \alpha}
	      +\norm{\Xd_s}^{2 \alpha}
	      \right\}^{\frac12} |\cc{\Xd_s^{(\lambda)}}-\cc{\Xd_s}|\,.
\end{align*}
% 
% we have
% \begin{align*}
%  &\Delta_s := |\partial_\upsilon \partial_\mu U(s,\Xd^{(1)}_s,\law{\Xd^{(1)}_s})(\cc{\Xd^{(\lambda)}_s})
%  - 
%  \partial_\upsilon \partial_\mu U(s,\Xd^{(1)}_s,\law{\Xd^{(1)}_s})(\cc{\Xd^{(1)}_s})
%  | 
%  \\
%  & \le L\left\{ 1 + |\cc{\Xd_s^{(\lambda)}}|^{2\alpha} 
% 	      + |\cc{\Xd_s^{(1)}}|^{2\alpha}
% 	      +\norm{\Xd_s^{(1)}}^{2\alpha}
% 	      \right\}^{\frac12} |\cc{\Xd_s^{(\lambda)}}-\cc{\Xd_s^{(1)}}|\,.
% \end{align*}
By Cauchy Schwarz inequality, we obtain
\begin{align}
&\ccesp{|\set{\partial_\upsilon \partial_\mu \cU(s,\Xd_s,\law{\Xd_s})(\cc{\Xd_s^{(\lambda)}})
  -
  \partial_\upsilon \partial_\mu \cU(s,\Xd_s,\law{\Xd_s})(\cc{\Xd_{s}}) \nonumber
  } \nonumber
  |} 
  \\
 &\le  \label{eq first part temp}
 C \sqrt{h_{i}} \left(1+\| \Xd_s^{(\lambda)} \|_{2 \alpha}^{\alpha} + \| \Xd_s \|_{2 \alpha}^{\alpha} \right) \, .
 \end{align}
We then observe that
\begin{align*}
\normalpha{\Xd_s^{(\lambda)}} + \normalpha{\Xd_s}
& \le
C\left(\normalpha{\Xd_{t_i}}+h_i\normalpha{\Ud_{t_i} - \Yd_{t_i}}+\sqrt{h_i}\right)
\\
& \le C \left(1+ \normalphas{\xi} + \delta \Delta_Y \right) \, ,
\end{align*}
where we used lemma \ref{le basic estimates for X&Y} for the
last inequality. Combining the last inequality with \eqref{eq first part temp} and using also
\eqref{eq control Delta mu 1}, we compute
%\textcolor{red}{
\begin{align*}
 |\cR^{\sigma \mu}_i| \le C h_i^{\frac32}\left(1+ \normalphas{\xi} + \delta \Delta_Y \right) \, ,
\end{align*}
%}
and then
 \begin{align}\label{eq control R sigma mu}
 \biggl\| \sum_{\ell=j_k}^{j_{k+1}-1} \vert \cR^{\sigma \mu}_\ell 
 \vert
 \biggr\|_{2\alpha}^2  &\le C \vert \pi \vert \delta^2 \left(1 + \delta^2\Delta_Y^2 + \normalphas{\xi}^2 \right)\;.
 %\nonumber
 %\\
 %&\le C h_i\left(h_i  +h_i^2( \norm{U(r_k,\Xd_{r_k},\law{\Xd_{r_k}}) - \eta}^2 + \norm{\xi}%^2)\right)\;.
 \end{align}

% Using Lemma \ref{le basic estimates for X&Y}, we then obtain
% \begin{align}
%  \ccesp{|\Delta_s|^2} \le Ch_i(1+ \ccesp{|\Xd_{t_i}|^2} + h_i^2\norm{\eta}^2 )
% \end{align}

% recall that 
% \textcolor{red}{$\desp{|\Xd_s|^2} \le C$}. Using H\"older's inequality, we compute
% \begin{align*}
%  \whEFp{t_i}{\Delta_s^2} \le C
%  \left\{ 1 + \whEFp{t_i}{|\cc{\Xd_s(\lambda)}|^{2\alpha(1+\epsilon)} }^{\frac1{1+\epsilon}}
% 	      + \whEFp{t_i}{|\cc{\Xd_s}|^{2\alpha(1+\epsilon)}}^{\frac1{1+\epsilon}}
% 	      \right\}
% 	      \whEFp{t_i}{|\cc{\Xd_s(\lambda)}-\cc{\Xd_s}|^{2p_\epsilon} }^\frac{1}{ p_\epsilon}\;.
% \end{align*}
% with $p_\epsilon := \frac{1+\epsilon}{\epsilon}$. Noticing that
% \begin{align*}
%  \whEFp{t_i}{|\cc{\Xd_s(\lambda)}-\cc{\Xd_s}|^{2p_\epsilon} }^\frac{1}{ p_\epsilon} \le C \sigma(\Xd_{t_i})^2 h\,,
% \end{align*}
% \textcolor{magenta}{
% OK si vol bornee, et $\alpha(1+\epsilon)<1$ si point de depart de $X$ dans
% $L^2$ mieux sinon...
% }
% %
% %

4. Collecting the estimates \eqref{eq control Rfi}, \eqref{eq control Rbi} and \eqref{eq control R sigma x}, we compute 
\begin{align*}
 {\left(\sum_{\ell = j_k}^{j_{k+1}-1} 
  \normalpha{\cR^f_\ell + \cR^b_\ell +\cR^{\sigma x}_\ell }\right)^2} 
 \le C   \delta^2 \left(\Delta_X^2 + |\pi| \set{1 + \normalphas{\xi}^2} + |\pi|^2	 \Delta_Y^2 \right) \, .
\end{align*}
Observing that 
\begin{align*}
 {\left(\sum_{\ell = j_k}^{j_{k+1}-1} 
  \normalpha{
  \delta {\mathcal T}(t_{i},t_{i+1})   }\right)^2} 
 \le C   \delta^2 \vert \pi \vert \, ,
\end{align*}
and combining the previous inequality with \eqref{eq control R sigma mu}, \eqref{eq control Rwi} and \eqref{eq stability}, we obtain
\begin{align*}
 \Delta_Y^2 \le e^{C\delta}\|\Ud_{r_{k+1}} - \Yd_{r_{k+1}}\|_{2\alpha}^2 
 + C \Bigl(\delta \Delta_X^2 +  \vert \pi \vert
 \bigl( 1 +  \normalphas{\xi}^2\bigr)
  +
 |\pi|\delta \Delta_Y^2  \Bigr) \, ,
\end{align*}
which concludes the proof for $\delta$ small enough.
\eproof
%\end{bluetext}
\vspace{5pt}

\begin{lemma}
\label{le:A3}
Assume that $g$ and $f(\cdot,0,0,\law{\cdot,0})$ are bounded. Then
$(A3)$ is satisfied whatever the value of $\alpha$. 
\end{lemma}

\proof
It suffices to prove that $\cU$ is bounded on the whole space and that $\bar{Y}$ 
is bounded independently of the discretization parameters. 

We refer to \cite{chacri15} for the proof of the boundedness of $\cU$. 

The bound for $\bar{Y}$ may obtained by squaring 
\eqref{eq scheme expli} and then by taking the conditional expectation exactly as done in the second step of the proof of 
Lemma \ref{le approx}.
\eproof

Assumptions $(A4)$ and $(A5)$ are easily checked. It suffices to observe that $({\tt P}_{r_{k},t}(\xi))_{t \in \pi, t \geq r_{k}}$ coincides
with the solution of the discrete Euler scheme:
\begin{equation*}
\bar{X}_{t_{i+1}}^0
=\bar{X}_{t_{i}}^0
+ (t_{i+1}-t_{i}) b\bigl(X_{t_{i}}^0,0,\law{X_{t_{i}}^0,0}\bigr)
+ \sqrt{t_{i+1}-t_{i}}
\sigma\bigl(X_{t_{i}}^0,\law{X_{t_{i}}^0}\bigr)
\varpi_{i},
\end{equation*}
with $\bar{X}_{r_{k}}^0=\xi$ as initial condition. 

Combining Lemma \ref{le approx},  Lemma \ref{le stab} 
and Lemma 
\ref{le:A3}
with Theorem \ref{th propagation error II}, we have
the following result.
\begin{corollary}
\label{cor:implemented}
Under \HYP{1}-\HYP{0}, assuming \eqref{eq main ass disc}, the following holds
 {\rm
 \begin{align*} %\label{eq main theorem}
 & \normalpha{\sol{k}{\xi} - \cU(r_k,\xi,\law{\xi})} 
\le C 
  \left( (C \delta)^{J-1}+ |\pi|^{\frac12}\delta^{-1} (1+\normalphas{\xi}) 
  \right) \; ,
%  + C|\pi|^{\frac12}\delta^{-1} \right)\;,
 \end{align*}
 }
 for $\bd$ small enough.
\end{corollary}
The first term in the right hand side is connected with the local Picard iterations on a
step of length $\delta$. As expected, it decreases geometrically fast with the number of iterations. 
The second term is due to the propagation of the error along the mesh. 
The leading term $\vert \pi \vert^{\frac12}$ is consistent with that observed for classical forward-backward systems, see for instance 
\cite{delmen06,delmen08}. The normalization by $\delta$ is due to the propagation of the error through the successive local solvers.

%is not so regular as it could be. We here pay for the fact that $\cU$ and its derivative are just assumed to be $1/2$-H\"older in time, locally in space, see
%\eqref{eq cont grad t},
%and that the second order derivatives in the measure argument may be of polynomial growth, see 
%\eqref{eq control dvdmu}. We are pretty sure that the rate would be stronger under stronger regularity assumptions. For consistency, we here chose assumptions in accordance with our previous result in \cite{chacri15}, which holds true under pretty general conditions. We felt that it was a good framework to prove the convergence of our algorithm, as the polynomial growth observed in bounds of the form 
%\eqref{eq control dvdmu} is a typical feature of the mean field structure of the problem.  

%\input{numerics}

%!TEX root = mainmkv.tex
%\newpage

\section{Numerical applications}
\label{se numerics}

In practice, we would like to approximate the value of $\cU(0,\cdot)$ at some point
$(x,\mu) \in \R^d\times\cP_2(\R^d)$. In the first section below, we explain how to retrieve such approximation using the approximation
of $\cU(0,\xi,\law{\xi})$ given by the algorithm $\sol{0}{}$, for some $\xi \sim \mu$. 
In a second part, we discuss the numerical results obtained by implementing $\sol{0}{}$ with two levels, i.e. $N=2$. In particular, we show that it
is more efficient than an algorithm based simply on Picard Iterations.

% In practise, we need to discretise the initial condition as well.
 %\textcolor{red}{easy???}

\subsection{Approximation of $\cU(0,x,\mu)$}
\label{subse approx U(t,x,mu)}
% \textcolor{red}{in fact describe two algorithms to do the computation and observe that very often we don't need to do the second one... it will be clearer}

The goal of this section is to show how to obtain an approximation of $\cU(0,x,\law{\xi})$ with $\xi \sim \mu$ and $x \in \mathrm{ supp}( \mu)$. We will assume that we thus have at hand a discrete valued random variable $\xi^{|\pi|} \sim \mu^{|\pi|}=\sum_{\ell=1}^M p_\ell \delta_{x\ell} $ such that $\mu^{\vert \pi\vert}$ is a good approximation of $\mu$ for the Wasserstein distance. For instance, such an approximation can be constructd by using quantization techniques. 
Then, we can use $\sol{0}{\xi^{|\pi|}}$ to obtain an approximation of $\cU(0,\xi^{|\pi|},\law{\xi^{|\pi|}})$.

%The algorithm based on Example \ref{ex binomial tree} is based on %iteration on trees: It is important to observe that, at time $0$, %$x_\ell$, $1 \le \ell \le M$, are the roots of $M$ trees.

Note that $\sol{0}{\xi^{|\pi|}}$ is a discrete random variable as the algorithm is initialised by a discrete random variable as well. In practice, this means that  each point  $x^\ell$ will be the root of a tree and will be associated to an output value $y^\ell=\cU(0,x^\ell, \law{\xi^{|\pi|}})$ and then
$\sol{0}{\xi^{|\pi|}} \sim \sum_{\ell=1}^M p_\ell \delta_{y^\ell} $.
It is important to remark that the computations on the trees are connected via the McKean-Vlasov interaction. 
%and $y^\ell$ is then an approximation of $U(0,x^\ell, %\law{\xi^{|\pi|}})$.

Using the Lipschitz continuity of $\cU$, one easily obtains
\begin{equation}
\label{eq first error}
\begin{split}
|\cU(0,x,\mu) - \cU(0,x^{\bar{\ell}},\mu^{|\pi|})| &\le C \left( \min_{y \in \mathrm{ supp}( \law{\xi^{|\pi|}}) }|y-x| + \cW_2(\mu^{|\pi|},\mu) \right) 
\\
&=: \cE_1(|\pi|,\xi)\,,
\end{split}
\end{equation}
where $x^{\bar{\ell}}$ is a point in the support of $\mu^{|\pi|}$ realising the minimum in the first line.

\begin{remark}
In many cases, it will be easy to have $x \in \mathrm{ supp}( \mu^{|\pi|})$ and thus reduce the above error to the term $\cW_2(\mu^{|\pi|},\mu)$. This is obviously the case if $\xi$ is deterministic.
\end{remark}

As mentioned above, the approximation of $\cU(0,x^{\bar{\ell}},\mu^{|\pi|})$ is obtained by running $\sol{0}{\xi^{|\pi|}}$ and by taking its value on the tree initiated at $x^{\bar{\ell}}$, precisely we have $\cU(0,x^{\bar{\ell}},\mu^{|\pi|}) = y^{\bar{\ell}}$. 
The corresponding pointwise error is given by
\begin{equation}
\label{eq second error}
\cE_2(|\pi|,\delta,\xi) :=
|y^{\bar{\ell}} - \cU(0,x^{\bar{\ell}},\law{\xi^{|\pi|}})| \, . 
\end{equation}
Of a course, this might be estimated by 
\begin{equation*}
\cE_2(|\pi|,\delta,\xi)
\le \frac1{p_{\bar{\ell}}}\norm{\cU(0,\xi^{|\pi|},\law{\xi^{|\pi|}}) - \sol{0}{\xi^{|\pi|}}} \; ,
\end{equation*}
but this is very poor when the initial distribution $\mu$ is diffuse and accordingly
when
$\mu^{\vert \pi \vert}$ has a large support, in which case $p_{\bar{\ell}}$ 
is expected to be small. 

To bypass this difficulty, we must regard $\cE_2(|\pi|,\delta,\xi)$ as a conditional error. 
Somehow, it is the error of the numerical scheme conditional on the initial root of the tree. It requires a new analysis, but it should not be so challenging: 
Now that we have investigated the error for the McKean-Vlasov component, we can easily revisit the proof 
of Theorem \ref{th propagation error II} in order to derive 
a bound for this conditional error. 

Instead of revisiting the whole proof, we can argue by doubling the variables.
For $\xi$ and $x$ as above, we can regard the 
four equations
\eqref{eq Xtxmu}, 
\eqref{eq Ytxmu}, 
\eqref{eq Xtxi}
and 
\eqref{eq Ytxi}
as a single forward-backward system of the McKean-Vlasov type. 
The forward component of such a doubled system is
${\mathbb X} = (X^{0,x,\mu},X^{0,\xi})$ and the backward components are ${\mathbb Y}=(Y^{0,x,\mu},Y^{0,\xi})$ and 
${\mathbb Z}=(Z^{0,x,\mu},Z^{0,\xi})$. Except for the fact that the dimension of ${\mathbb X}$ is no longer equal to the dimension of the noise, which we assumed to be true for convenience only,
and for the fact that ${\mathbb Y}$ takes values in $\RR^2$,
the setting is exactly the same as before, namely 
 $({\mathbb X},{\mathbb Y},{\mathbb Z})$ can be regarded 
 as the solution of a McKean-Vlasov  forward-backward SDE in
 which the mean field component reduces to the marginal law of 
 $(X^{0,\xi},Y^{0,\xi})$. We observe in particular that 
 \begin{equation*}
 Y^{0,x,\mu}_{t} = \cU(t,X^{0,x,\mu}_{t},\law{X_{t}^{0,\xi}}),
 \quad 
 Y^{0,\xi}_{t} = \cU(t,X_{t}^{0,\xi},\law{X_{t}^{0,\xi}}), \quad t \in [0,T],
 \end{equation*}
 with similar relationships for $Z^{0,x,\mu}$ and $Z^{0,\xi}$. Hence, 
${\mathbb Y}_{t}$ (and ${\mathbb Z}_{t}$) can be represented as a function of ${\mathbb X}$, which was the key assumption in our analysis. 
For sure, the fact that ${\mathbb Y}$ takes values in dimension 2
is not a limitation for duplicating the arguments used to prove Theorem 
\ref{th propagation error II}. 

Numerically speaking, the tree initiated at root $x^{\bar{\ell}}$
under the initial distribution $\mu^{\vert \pi\vert}$ provides an approximation of $\cU(0,x^{\bar{\ell}},\law{\xi^{\vert \pi \vert}})$, which is equal to $Y^{0,x^{\bar{\ell}},\law{\xi^{\vert \pi\vert}}}$. So our numerical (implemented) scheme is in fact a numerical for the whole process 
$({\mathbb X},{\mathbb Y},{\mathbb Z})$.

This leads us to the following result.
\begin{theorem}
 Let $y^{\bar{\ell}}$ be the approximation of $\cU(0,x,\mu)$ obtained by calling  {\rm $\sol{0}{\xi^{\vert \pi \vert}}$}, where $\bar{\ell}$ is defined in \eqref{eq first error}. Then, the following holds
 \begin{align*}
  |\cU(0,x,\mu)- y^{\bar{\ell}}| \le \cE_1(|\pi|,\xi) + \cE_2(|\pi|,\delta,\xi)\, ,
 \end{align*}
where 
$\cE_2(|\pi|,\delta,\xi)$ can be estimated by 
Corollary 
\ref{cor:implemented}, with 
$(1+\normalphas{\xi})$ replaced by 
$(1+\vert x^{\bar{\ell}} \vert +\normalphas{\xi})$.
%recall \eqref{eq first error} and \eqref{eq second error}.
\end{theorem}

%Obviously, the second term is rather poor when $p_{\bar{\ell}}$ is small, which is for instance the case when the initial distribution of $\xi$ is diffuse. 
%This is well-expected: The global error of the algorithm, see for instance
%Theorem \ref{th propagation error II},
% is in $L^p$-norm for some $p \geq 2$; some price has to be paid to recover a pointwise estimate.

\subsection{Numerical illustration}

In this section, we will prove empirically the convergence of the 
approximation obtained by the solver $\sol{}{}$. In particular, we will compare the output of our algorithm $\sol{}{}$, when implemented with two levels, i.e. $N=2$ (we simply call it \emph{two-level algorithm}), with the output of a basic algorithm based only on Picard iterations, which can be seen as a solver $\sol{}{}$, but with only one level, i.e. $N=1$ (we simply call it \emph{one-level algorithm}). In both cases, we use 
Example \ref{ex binomial tree} as discretization scheme, with a standard Bernoulli quantization of the normal distribution, $d$ being equal to $1$. In the numerical studies below, we show that  the \emph{two-level algorithm} converges in case when the \emph{one-level algorithm} fails.

%Let us note that in the graphs below $Y0$ refers to $U(0,x,\delta_x)$. 

\subsubsection{The example of a linear model}
In this part, we compare the output of both algorithms for the following linear model where a closed-form solution is available:
 \begin{align*}
  \ud X_t &= - \rho \esp{Y}_t \ud t + \sigma \ud W_t\;, \quad \; X_0 = x\,,
  \\
  \ud Y_t &= -a Y_t \ud t  + Z_t \ud W_t \, , \quad \;\text{ and } Y_T = X_T\, ,
 \end{align*}
 for $\rho,a>0$,
 and the true solution for $\esp{X_0} = m_0$ is given by
\begin{align*}
 Y_0 = \frac{m_0 e^{aT}}{ 1 + \frac{\rho}{a} (e^{aT}-1) }\;.
\end{align*}
%The parameters $\rho$ and $a$ are fixed to \textcolor{magenta}{?check the code?}. 
The errors for various time steps and for both algorithms are shown on the  log-log error plot of Figure $1$.  %The value \textcolor{red}{$Y0$ stands for the approximation of $\cU(0,x,\delta_x)$}. 
The parameters are fixed as follows: $\rho = 0.1$, $a = 0.25$, $\sigma=1$, $T = 1$ and $x=2$. Moreover, the \emph{two-level algorithm} uses $5$ Picard Iterations {per level}, and the \emph{one-level algorithm} computes $25$ Picard Iterations.

%\begin{center}
\begin{figure}[h] \label{fig safety check}
\includegraphics[width=0.8\textwidth]{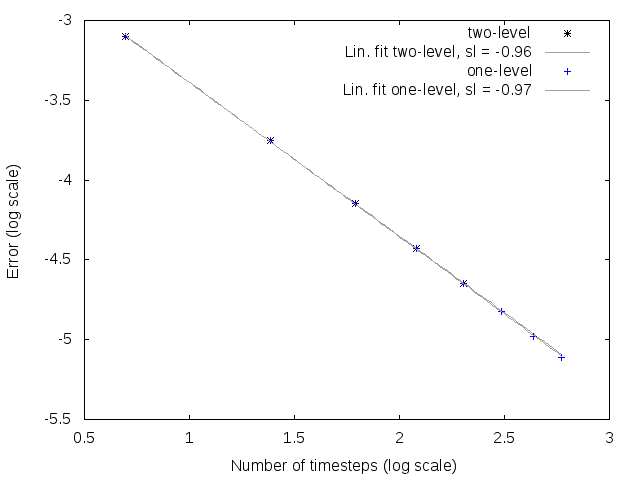}%{./img/safety_check.png}
\caption{Convergence of the algorithms: log-log error plot {for the same data as in the text}. We can observe that both algorithms return the same value which is close to the true value. This validates the convergence of both methods in this simple linear setting.}
\end{figure}
%\end{center}

%On Figure $1$,  

\vspace{5pt} 
% The linear model is the following
% \begin{align*}
%  \ud X_t &= - \epsilon \esp{Y_t} + \sigma \ud W_t
%  \text{ and }
%  Y_t = \EFp{t}{X_T + \int_t^T aY_s \ud s}
% \end{align*}
% and the true solution for $\esp{X_0} = m_0$ is given by
% \begin{align*}
%  Y_0 = \frac{m_0 e^{aT}}{ 1 - \frac{\epsilon}{a} (1 - e^{aT}) }\;.
% \end{align*}

% 
% \textcolor{magenta}{check the code for the data
% %\begin{itemize}
% %  \item The linear model:
% %  \begin{align*}
% %   \ud X_t &= - \rho \esp{Y}_t \ud t + \sigma \ud W_t\;,\; X_0 = x\,,
% %   \\
% %   \ud Y_t &= -a Y_s \ud s  + Z_s \ud W_s \;\text{ and } Y_T = X_T\,.
% %  \end{align*}
% %  and the true solution for $\esp{X_0} = m_0$ is given by
% % \begin{align*}
% %  Y_0 = \frac{m_0 e^{aT}}{ 1 - \frac{\rho}{a} (1 - e^{aT}) }\;.
% % \end{align*}
% % \item The coupling parameter is fixed.
% %\item 
% We study the convergence of the discretisation error for both method
% \begin{enumerate}
%  \item Picard Iteration (25 iterations)
%  \item $\sol{}{}$ with two levels (5 Picard iterations each) 
% \end{enumerate}
% %\end{itemize}
% %
% }

\subsubsection{Efficiency of the $\sol{}{}$ algorithm}

In this section, we compare the \emph{two-level algorithm} and the \emph{one-level algorithm} on two models, for which existence and uniqueness to the master equation (or the FBSDE system)  hold true for any arbitrary terminal time $T$ and Lipschitz constant $L$ of the coefficients function. Nevertheless, as stated in the  theorems above, the convergence of the algorithms is guaranted only for a periods of time which are controlled by $L$ and $T$. Here, we fix the terminal date $T$ and allow $L$ to vary with the use of a coupling parameter $\rho$, see equations \eqref{ex:ex1} (for a case without McKean-Vlasov interaction)
and 
 \eqref{ex:ex2}
  (for a case with McKean-Vlasov interaction). We will see below that, as expected, the \emph{two-level algorithm} converges for a larger range of coupling parameter than the \emph{one-level algorithm}.

% 
% As shown below, the  $\sol{}{}$ algorithm is better than the ``basic''one
% as it allows to compute an approximation of the true solution for larger 
% coupling parameter. This is demonstrated on two models: the first one is a standard FBSDE,
% the second one is a FBSDE with McKean-Vlasov interaction. Both models admit, in theory,
% solution for arbitrary value of the parameter.

%\vspace{5pt}
%On both graphs: x-axis is coupling parameter, y-axis is output value of the %algorithm (plot 5 last iteration value).
\paragraph{An example with no McKean-Vlasov interaction} \textcolor{white}{.}\\

\noindent Here, the model is the following
\begin{equation}
 \label{ex:ex1}
\begin{split}
 \ud X_t &= \rho \cos(Y_t) \ud t + \sigma \ud W_t, \quad \; X_0 = x \, ,
 \\
 Y_t &= \EFp{t}{ \sin(X_T)} \, .
\end{split}
\end{equation}

%\textcolor{magenta}{Redo the graph..}

On Figure $2$, we plot the output of the two-level and one-level algorithm along with a proxy of the true solution computed by usual BSDE approximation method (after a Girsanov transform) and with a very high-level of precision. On the graph, the value $Y0$ stands for the approximation of $\cU(0,x)$: There is no dependence upon the initial measure as there is no MKV interaction in this example. The parameters are fixed as follows:  $\sigma=1$, $T = 1$ and $x=0$. Moreover, the \emph{two-level algorithm} uses $5$ Picard Iterations {per level}, and the \emph{one-level algorithm} computes $25$ Picard Iterations.

\begin{figure}[h]
\includegraphics[width=0.8\textwidth]{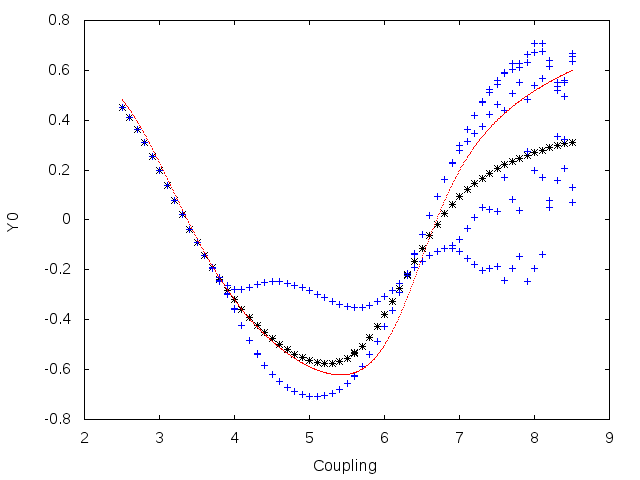}%{./img/fbsde_bifur_nomkv.png}
\caption{
Comparison of algorithms' output for different value of the coupling parameter and {for the same data as in 
Example  \eqref{ex:ex1}}: two-level (black star), one-level (blue cross), 
true value (red line). The \emph{two-level algorithm} converges for larger coupling parameter than the \emph{one-level algorithm}. It is close to the true solution up to parameter $\rho=7$, the discrepancy for large coupling parameter coming most probably from the discrete-time error. Interestingly, the \emph{one-level algorithm} shows bifurcations. 
} 
\end{figure}

%\includegraphics[width=\textwidth]{./img/pres_nomkv_bifur.png}

%\emph{legend.} black:``solver'', blue:``basic'', red: ``true''.
%\textcolor{red}{Peux tu d\'etailler plus ce que tu as repr\'esent\'e ? Abscisses %= couplage ? Ordonn\'ees = valeur $t=0$, $x=?$ ? Bifurcation ?}

\paragraph{An example from large population stochastic control} \textcolor{white}{.}\\

\noindent For this part, the model is given by 
%\begin{itemize}
 %\item The model
 \begin{equation}
 \label{ex:ex2}
 \begin{split}
  \ud X_t & = - \rho Y_t \ud t + \ud W_t\,,\, X_0 = x\,,
  \\
  \ud Y_t & =  \mathrm{atan}(\esp{X_t}) \ud t + Z_t \ud W_t \text{ and } Y_T = G'(X_T):=\mathrm{atan}(X_T).
 \end{split}
 \end{equation}
 coming from Pontryagin principle applied to MFG
 \begin{align*}
  \inf_\alpha \esp{G(X^\alpha_t) + \int_0^T\left(\frac1{2\rho} \alpha^2_t + X_t^\alpha \mathrm{atan}(\esp{X^\alpha_t})\right) \ud t  }
 \end{align*}
 with $\ud X^\alpha_t =  \alpha_t \ud t + \ud W_t$, see e.g. \cite{CarmonaDelarue_sicon}.
%\end{itemize}

We do not know the exact solution for this model and it is not possible to obtain easily an approximation as in the previous example. We plot on Figure $3$, the output value of the \emph{one-level algorithm} and \emph{two-level algorithm}. On the graph, the value $Y0$ stands for the approximation of $\cU(0,x, \delta_x)$. The parameters are fixed as follows:  $\sigma=1$, $T = 1$ and $x=1$. Moreover, the \emph{two-level algorithm} uses $5$ Picard Iterations {per level}, and the \emph{one-level algorithm} computes $25$ Picard Iterations.

%On the graph below, 

\begin{figure}[h]
 \includegraphics[width=0.8\textwidth]{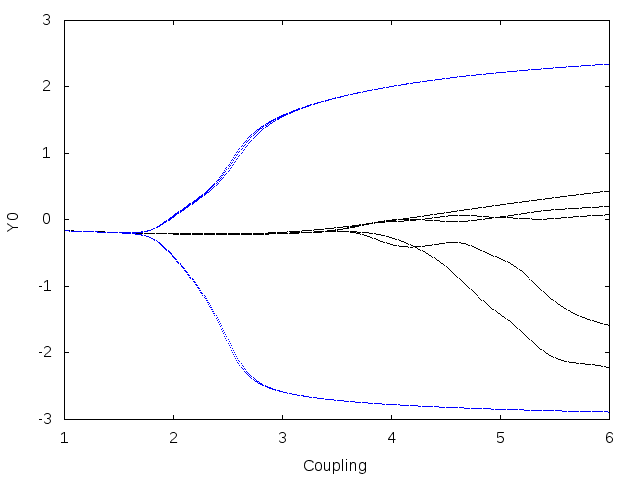}%{./img/mkvfbsde_bifur_francois.png}
\caption{Algorithms' output for {the same data as in Example 
 \eqref{ex:ex2}}: \emph{one-level algorithm} (blue line), \emph{two-level algorithm} (black line). We observe the same phenomenon as in the previous model: The \emph{two-level algorithm} {converges} to a unique value for a larger range of coupling parameter than the \emph{one-level algorithm}, which exhibits a bifurcation. Observe that the \emph{two-level algorithm} fails to converge at some points: One should add a level of computation to shorten the time period $\delta$.}
 \end{figure}
% 
% \vspace{5pt}
% Here, the model is the following
% \begin{align*}
%  \ud X_t &= -\epsilon \esp{\cos(Y_t)} \ud t + \sigma \ud W_t
%  \\
%  Y_t &= \EFp{t}{ \cos(X_T)}
% \end{align*}
% 
% comment: check assumption on model and also very strong coupling, maybe change model... observe jump around $10.5$ (why?)
% 

% 
% 
% \textcolor{magenta}{
% numerics
% \begin{enumerate}
%  \item Picard Iterations (25) - in \emph{blue}
%  \item $\tsol{}{}$ with two levels (5 iterations per level) - in \emph{black}
% \end{enumerate}
% }

%\input{appendix}

%!TEX root = mainmkv.tex
%\newpage
\section{Appendix}

\subsection{A discrete It\^o formula}

We consider the following  Euler scheme on the discrete time grid $\pi$ of the interval $[0,T]$, recall \eqref{eq de disc grid pi},
\begin{align}\label{eq de app euler scheme}
 \Xd_{t_{i+1}} = \Xd_{t_{i}} + b_i (t_{i+1} - t_i) 
                   + \sigma_i \sqrt{t_{i+1} - t_i} \varpi_i \;, 
\end{align}
where $(\varpi_i)_{i_ \le n}$ are i.i.d. centered $\RR^d$-valued random variables
such that the covariance matrix $\EE[ \varpi_{i} \varpi_{i}^\dagger]$ is the identity matrix and $\normalphas{\varpi_{i}} ^2
\leq \Lambda h_{i}$, and $(b_i,\sigma_i) \in L^2(\cF_{t_i})$, for all $i \le n$.
\vskip 5pt

\noindent We also introduce a  piecewise continuous version of the previous scheme, for $i<n$, $t_i \le s < t_{i+1}$
% \begin{align}\label{eq de app euler scheme cont}
%  \Xd_{s} = \Xd_{t_{i}} + b_i (s - t_i) 
%                    + \sigma_i \sqrt{s - t_i} \eta_i \;
% \end{align}
and  $\lambda \in [0,1]$, the process $(\Xd_t^{(\lambda)})_{0 \le t \le T}$,
\begin{align}\label{eq de app euler scheme cont}
 \Xd_{s}^{(\lambda)} = \Xd_{t_{i}} + b_i (s - t_i) 
                  + \sigma_i \lambda \sqrt{s - t_i} \varpi_i \;
 \end{align}
 and $\Xd_{t_n}^{(\lambda)} = \Xd_{t_n}$.
 Following the notation used in the proof of Lemma
 \ref{le approx}, we just write 
 $(\Xd_{s})_{0 \le s \le T}$ for 
 $(\Xd_{s}^{(1)})_{0 \le s \le T}$, which defines a continuous version of the Euler scheme given in \eqref{eq de app euler scheme}.

\begin{proposition}\label{pr disc ito formula}
For any $i \in \{0,\cdots,n-1\}$,
the following holds true:
\begin{align*}
 &\cU(t_{i+1},\Xd_{t_{i+1}},\law{\Xd_{t_{i+1}} })
 = \cU(t_i,\Xd_{t_i},\law{\Xd_{t_i} })
 + \int_{t_i }^{t_{i+1}} \partial_t \cU(s,\Xd_s,\law{\Xd_s}) \ud s
 \\
 &+ \int_{t_i}^{t_{i+1}} \left(
  \partial_x \cU(s,\Xd_s,\law{\Xd_s}) \cdot b_i  
  + 
\frac12\int_0^1 \textrm{\rm Tr} \bigl[ \partial^2_{xx} \cU(s,\Xd^{(\lambda)}_s,\law{\Xd_s}) a_i \bigr] \ud \lambda
   \right)
   \ud s
\\
&+ 
 \int_{t_i}^{t_{i+1}}  
\ccesp{\partial_\mu \cU(s,\Xd_s,\law{\Xd_s})(\cc{\Xd_s})
\cdot \cc{b_i}} ds
\\
&
  +\frac12\int_0^1 
\ccesp{  \textrm{\rm Tr} \bigl[
  \partial_\upsilon \partial_\mu \cU(s,\Xd_s,\law{\Xd_s})(\cc{\Xd^{(\lambda)}_s})
 \cc{a_i}\bigr] \ud \lambda 
 }
  \ud s
\\
&+ \int_{t_i}^{t_{i+1}}\partial_x \cU(s, \Xd^{(0)}_s, \law{\Xd_s})
  \frac{\sigma_i \varpi_i}{2 \sqrt{s-t_i}} \ud s  
  + \delta {\mathcal M}(t_{i},t_{i+1}) + \delta {\mathcal T}(t_{i},t_{i+1})
  \;,
 \end{align*}
 where $a_{i}$ is here equal to $\sigma_{i} \sigma_{i}^\dagger$, 
 and $\delta {\mathcal M}(t_{i},t_{i+1})$ is a martingale increment satisfying $\normalphas{\delta {\mathcal M}(t_{i},t_{i+1})} \leq C_{\Lambda} h_{i}^{2}$ and 
 $\normalphas{\delta {\mathcal T}(t_{i},t_{i+1})} \leq C_{\Lambda} h_{i}^{\frac32}$.
 \end{proposition}

\proof
%We consider $0 < h < t_{i+1}-t_i$. 
By writing
\begin{equation*}
\Xd_{t_{i+1}} = \Xd_{t_{i}}
+ \int_{t_{i}}^{t_{i+1}} 
\bigl(b_i +\frac{\sigma_i \varpi_i}{2\sqrt{s-t_i}}\bigr) ds,
\end{equation*}
and by using the standard chain rule 
for continuously differentiable functions on a Hilbert space, 
we get
\begin{align*}
 \cU(t_{i+1},\Xd_{t_{i+1}},\law{\Xd_{t_{i+1}} })
 &= \cU(t_{i},\Xd_{t_{i}},\law{\Xd_{t_{i}} })
 + \int_{t_i}^{t_{i+1}} \partial_t \cU(s,\Xd_s,\law{\Xd_s}) \ud s
 \\
 &\hspace{15pt}+ \int_{t_i}^{t_{i+1}}
 \left( \partial_x \cU(s,\Xd_s,\law{\Xd_s}) \cdot \bigl(b_i +\frac{\sigma_i \varpi_i}{2\sqrt{s-t_i}}\bigr) \ud s
\right.
\\
&\hspace{30pt} + \left. \ccesp{\partial_\mu \cU(s,\Xd_s,\law{\Xd_s})(\cc{\Xd_s})
\cdot 
             \cc{b_i +\frac{\sigma_i \varpi_i}{2\sqrt{s-t_i}}}} \right)
  \ud s\;.
 \end{align*}
Now we observe that,
\begin{align*}
 \partial_x \cU(s,\Xd_s,\law{\Xd_s}) &= \partial_x \cU(s,\Xd^{(0)}_s,\law{\Xd_s})
+\sqrt{s-t_i} \int_0^1 \partial^2_{xx} \cU(s,\Xd^{(\lambda)}_s,\law{\Xd_s})
\sigma_i\varpi_i \ud \lambda
\\
&= \partial_x \cU(s,\Xd^{(0)}_s,\law{\Xd_s})
+ \sqrt{s-t_i}\, \partial^2_{xx} \cU(s,\Xd^{(0)}_s,\law{\Xd_s}) \sigma_i\varpi_i   + \sqrt{s-t_{i}}{\mathcal T}_{1}(s) \, ,
 \end{align*}
 where 
 ${\mathcal T}_{1}(s)$ is a random variable 
 defined on $(\Omega,\cF,\PP)$ such that 
 $\| {\mathcal T}_{1}(s) \|_{2 \alpha} \leq 
 C h_{i}^{\frac12}$,
and
\begin{align*}
 \partial_\mu \cU(s,\Xd_s,\law{\Xd_s})(\cc{\Xd_s}) &=
 \partial_\mu \cU(s,\Xd_s,\law{\Xd_s})(\cc{\Xd^{(0)}_s})
 \\
 &\hspace{15pt}+\sqrt{s-t_i} \int_0^1 \partial_\upsilon \partial_\mu \cU(s,\Xd_s,\law{\Xd_s})(\cc{\Xd^{(\lambda)}_s})
  \cc{\sigma_i  \varpi_i}\ud \lambda
 \\
 &= \partial_{\mu} \cU(s,\Xd_s,\law{\Xd_s})(\cc{\Xd^{(0)}_s})
 \\
 &\hspace{15pt}
+ \sqrt{s-t_i} \, \partial_{v} \partial_{\mu} \cU(s,\Xd_s,\law{\Xd_s})(\cc{\Xd_{s}
^{(0)}})
 \cc{\sigma_i\varpi_i} 
+ \sqrt{s-t_{i}} {\mathcal T}_{2}(s) \, ,
\end{align*}
where ${\mathcal T}_{2}(s)$ is a random variable on the enlarged space $(\Omega \times \hat{\Omega},\cF \otimes \hat{\cF},\PP \otimes \hat{\PP})$ such that 
$\ccesp{\vert {\mathcal T}_{2}(s) \vert^{2\alpha}}^{1/(2\alpha)} \leq 
C h_{i}^{\frac12}$. 

We insert these expansions back into the identity we obtained 
for the term
$\cU(t_{i+1},\Xd_{t_{i+1}},\law{\Xd_{t_{i+1}} })$.
We let
\begin{equation*}
\begin{split}
\delta {\mathcal M}(t_{i},t_{i+1}) &= 
 \frac12
\int_{t_{i}}^{t_{i+1}} 
\Bigl[
\partial^2_{xx} \cU(s,\Xd^{(0)}_{s},\law{\Xd_{s}}) 
\sigma_{i} \varpi_{i} \cdot 
\bigl( 
\sigma_{i} \varpi_{i} \bigr)
\\
&\hspace{50pt}
-
\EFp{t_i}{\partial^2_{xx} \cU(s,\Xd^{(0)}_{s},\law{\Xd_{s}}) 
\sigma_{i} \varpi_{i} \cdot 
\bigl( 
\sigma_{i} \varpi_{i} \bigr)}
\Bigr]
 \ud s \, ,
\\
\delta {\mathcal T}(t_{i},t_{i+1}) &= 
\frac12
\int_{t_{i}}^{t_{i+1}} 
\bigl( {\mathcal T}_{1}(s) +  {\mathcal T}_{2}(s) \bigr)
\cdot \sigma_{i} \varpi_{i}
 \ud s \, . 
\end{split} 
\end{equation*}
It defines a martingale increment satisfying 
$\EFp{t_i}{ \vert \delta {\mathcal M}(t_{i},t_{i+1}) \vert^{2\alpha} }^{1/(2\alpha)}
 \leq C h_{i}$. 
Observing that for $ t_i \le s \le t_{i+1}$,
\begin{align*}
& \ccesp{
 \partial_\mu \cU(s,\Xd_s,\law{\Xd_s})(\cc{\Xd_s^{(0)}}) \cdot \cc{\sigma_{i}\varpi_i}} = 0 \, ,
 \\
 &\EFp{t_i}{\partial^2_{xx} \cU(s,\Xd^{(0)}_s,\law{\Xd_s})
 \sigma_{i} \varpi_{i}
 \cdot
\bigl( \sigma_{i} \varpi_{i}\bigr) }
= \EFp{t_i}{\textrm{\rm Tr} \bigl(\partial^2_{xx} \cU(s,\Xd^{(0)}_s,\law{\Xd_s})
a_{i} \bigr) } \, ,
 \\
 &\EFp{t_i}{\partial^2_{xx} \cU(s,\Xd^{(0)}_s,\law{\Xd_s})
 \sigma_{i} \varpi_{i}
 \cdot
\bigl( \sigma_{i} \varpi_{i}\bigr) }
= \EFp{t_i}{\textrm{\rm Tr} \bigl(\partial^2_{xx} \cU(s,\Xd^{(0)}_s,\law{\Xd_s})
a_{i} \bigr) } \, ,
\\
&\ccesp{
\partial_{v} \partial_{\mu} \cU(s,\Xd_s,\law{\Xd_s})
(\cc{\Xd_s^{(0)}})
\sqrt{s-t_i} 
\cc{\sigma_i\varpi_i}
\cdot
\cc{\sigma_i\varpi_i}
}  = \ccesp{\textrm{\rm Tr} \bigl(\partial_{v} 
\partial_{\mu} \cU(s,\Xd_s,\law{\Xd_s})(\cc{\Xd_s^{(0)}})
\cc{a_{i}} \bigr) } \, ,
\end{align*}
we complete the proof.
%obtain the proof of the Lemma letting $h$ goes to zero and
%noting that due to the bound on $\partial_x \cU$, the last integral is well defined (at least in a generalised sense)
%.
% taking conditional expectation
% on both side, we get
% \begin{align*}
%  &\EFp{t_i}{U(t_{i+1},\Xd_{t_{i+1}},\law{\Xd_{t_{i+1}} })}
%  = \EFp{t_i}{U(t_{i}+h,\Xd_{t_{i}+h},\law{\Xd_{t_{i}+h} })
%  + \int_{t_i + h}^{t_{i+1}} \partial_t U(s,\Xd_s,\law{\Xd_s}) \ud s}
%  \\
%  &+ \int_{t_i + h}^{t_{i+1}}\EFp{t_i}{
%   \partial_xU(s,\Xd_s,\law{\Xd_s})b_i  
%   + 
% \frac12\int_0^1 \partial^2_{xx}U(s,\Xd_s(\lambda),\law{\Xd_s})\sigma_i^2 \eta_i^2 \ud \lambda
%   } \ud s
% \\
% &+ 
% \int_{t_i + h}^{t_{i+1}}\EFp{t_i}{
% \ccesp{\partial_\mu U(s,\Xd_s,\law{\Xd_s})(\cc{\Xd_s})
%              \cc{b_i}} 
%   +\int_0^1 \partial_\upsilon \partial_\mu U(s,\Xd_s,\law{\Xd_s})(\cc{\Xd_s(\lambda)})
%  \cc{\sigma^2_i  \eta^2_i} \ud \lambda
%   }\ud s\;.
%  \end{align*}
% Letting $h$ goes to zero concludes the proof of the Lemma.
\eproof

\subsection{Estimates for the scheme given in Example \ref{ex binomial tree} }

\begin{lemma}\label{le basic estimates for X&Y}
Under \HYP{0}-\HYP{1}, the following holds for the forward component of the scheme given in Example \ref{ex binomial tree} and
its continuous version,
\begin{align}
\max_{t \in \pi^k} \normalpha{ \Xd_{t} } 
&\le C_\Lambda\left(1 + \normalpha{\Xd_{r_k}} 
+ \delta \max_{t \in \pi^k} \normalpha{\cU(t,\Xd_{t},\law{\Xd_{t}}) - Y_{t}}\right)\,,
\label{eq control Xd}
%\\
%\norm{\Xd_s^{(\lambda)}-\Xd_{t_i}} & \le C_\Lambda\left(\sqrt{t_{i+1} - t_i} + (t_{i+1}-t_i) \norm{\eta} %\right)\;.
\end{align}
%Moreover, the following holds true for its backward component
%\begin{align}
% \max_{j_k \le i <j_{k+1}} \norm{\Yd_{t_i}} &\le e^{C_\Lambda\delta}(\norm{\eta} + C_\Lambda \sqrt{\delta})\;.
%\end{align}
\end{lemma}

\proof
%1. We now study the estimate \eqref{eq control Xd}.
We introduce $\mathfrak{d}_i := |\cU(t_i,\Xd_{t_i},\law{\Xd_{t_i}})  - \Yd_{t_i}|$
%\begin{align*} 
%\mathfrak{d}_k := 
%\max_{j_k \le i < j_{k+1}} |U(t_i,\Xd_{t_i},\law{\Xd_{t_i}})  - \Yd_{t_i}|
%\end{align*}
and observe from the Lipschitz property of $b$ and $\cU$ that
\begin{align} \label{eq up bound b}
 \bigl|b\bigl(\Xd_{\ti{}},\Yd_{\ti{}},\law{\Xd_{\ti{}},\Yd_{\ti{}}}\bigr)\bigr| \le C_\Lambda
 \bigl(1 + |\Xd_{t_i}| +
 \normalpha{\Xd_{t_i}} + \mathfrak{d}_i + \normalpha{\mathfrak{d}_i}
 \bigr)\;.
\end{align}
Recall that the scheme for the forward component reads
\begin{align*}
 \Xd_{\ti{+1}}  = \Xd_{r_k} + \sum_{\ell = j_k}^{i}b
 \bigl(\Xd_{t_\ell},\Yd_{t_\ell},\law{\Xd_{t_\ell},\Yd_{t_\ell}}\bigr)(t_{\ell+1}-t_\ell) +  \sum_{\ell = j_k}^{i} \sigma\bigl(\Xd_{t_\ell},
\law{\Xd_{t_\ell}}
\bigr)\Delta \bar{W}_\ell\;.
\end{align*}
Squaring the previous inequality, using Cauchy-Schwarz inequality for the first sum and the martingale property for the second sum, we obtain
\begin{align*}
 \normalpha{\Xd_{\ti{+1}}}^2
 &\le C
 \normalpha{\Xd_{r_{k}}}^2
 + C \sum_{\ell = j_k}^{i}
  h_\ell \left( \delta 
  \normalpha{b(\Xd_{t_{\ell}},\Yd_{t_\ell},\law{\Xd_{t_\ell},\Yd_{t_\ell}})}^2  +  \|\sigma(\Xd_{t_\ell},
  \law{\Xd_{t_\ell}} 
  )\|^2_{2 \alpha} \right) \; \ ,
\end{align*}
where we used again B\"urkholder-Davis-Gundy inequality for discrete martingales. 

% then using the identity $|y|^2 = |x|^2 + 2x(y-x) + |y-x|^2$ with $y=\Xd_{\ti{+1}}$ and $x =\Xd_{\ti{}}$, we compute
% \begin{align*}
%  \esp{|\Xd_{\ti{+1}}|^2} \le& \esp{|\Xd_{\ti{}}|^2}
%  + 2 \esp{|\Xd_{\ti{}}b(\Yd_{\ti{}},\law{\Xd_{\ti{}},\Yd_{\ti{}}})|}h_i
%  \\
%  & + h_i^2\esp{|b(\Yd_{\ti{}},\law{\Xd_{\ti{}},\Yd_{\ti{}}})|^2
%  + |\sigma(\Xd_{\ti{}})\Delta \bar{W}_i|^2
%  }
% \end{align*}
Combining \eqref{eq up bound b} with the boundedness of $\sigma$, we then have
\begin{align*}
 \normalpha{\Xd_{\ti{+1}}}^2 
 \le C
 \Bigl( 
 \normalpha{\Xd_{r_{k}}}^2 
+\delta + \delta^2\max_{j_k \le i < j_{k+1}} \normalphas{\mathfrak{d}_i}^2 + C\delta\sum_{\ell=j_k}^{i}
h_{\ell}
\normalpha{\Xd_{t_\ell}}^2
 \Bigr) \; .
\end{align*}
Using the discrete version of Gronwall's lemma, 
the result easily follows. 
%to
%\begin{align*}
% \max_{j_k \le i \le j_{k+1}} \esp{|\Xd_{t_i}|^2}
% \le
% C\left(\delta + \esp{|\Xd_{r_k}|^2} 
% + \delta^2\max_{j_k \le i < j_{k+1}}\esp{|\mathfrak{d}_i|^2} \right)
%\end{align*}
%and conclude the proof for this step.
% 
% 
% 
% 
% 2. For the second estimate, we compute
% \begin{align*}
% \norm{\Xd_s^{(\lambda)}-\Xd_{t_i}^{(\lambda)}} & \le C_\Lambda(\sqrt{t_{i+1} - t_i} + (t_{i+1}-t_i)\max_{j_k \le i <j_{k+1}} \norm{\Yd_{t_i}} )\;.
%\end{align*} 
\eproof

\end{document}